\documentclass[amsppt,11pt]{amsart}
\usepackage{curves}
\usepackage{amsmath}
\def\Line#1{\hbox to\textwidth{#1}}%

\def\orig{{\bf 0}}
\def\LA{{\mathcal A}}
\def\Ba{{B}}
\def\bw{{\bold w}}
\def\bu{{\bold u}}
\def\bv{{\bold v}}

\def\bx{{\bold x}}

\def\gas{{\gamma^\star}}
\def\ba{{\bold a}}

\def\N{{\hbox{I\kern-.15em\hbox{N}}}}
\def\R{{\hbox{I\kern-.15em\hbox{R}}}}
\def\Zd{\mathbb Z^d}
\def\Zt{{\hbox{Z\kern-.4em\hbox{Z}}^2}}

\def\part{\partial_{\rm ext}{\mathcal C}(v_0)}

\def\piab1{\prod \ [a_i - 1, b_i + 1]}

\newcommand{\wcon}{\leftrightarrow}

\newtheorem {theorem}{Theorem}[section]
\newtheorem {prop}[theorem]{Proposition}
\newtheorem {lemma}[theorem]{Lemma}
\newtheorem {corollary}[theorem]{Corollary}
\newcounter{conjecture}
\newcounter{remark}

\newcommand{\eqnsection}{
   \renewcommand{\theequation}{\thesection.\arabic{equation}}
   \makeatletter
   \csname @addtoreset\endcsname{equation}{section}
   \makeatother}
\newcommand{\req}[1]{(\ref{#1})}
\newtheorem {defin}{Definition}[section]

\def \bt{\begin{theorem}}
\def \et{\end{theorem}}
\def \bs{\begin{split}}
\def \es{\end{split}}
\def \bea{\begin{eqnarray}}
\def \eea{\end{eqnarray}}
\def \bas{\begin{eqnarray*}}
\def \eas{\end{eqnarray*}}
\newcommand{\comment}[1]{}

\def \al{\alpha}

\def \ga{\gamma}
\def \Ga{\Gamma}

\def \eps{\epsilon}

\def \la{\lambda}

\def \R{{\bf R}}

\def \BB{{\mathcal B}}

\def \LL{{\mathcal L}}

\def \WW{{\mathcal W}}

\newcommand{\PPP}{{\mathbb P}}
\newcommand{\beq}{\begin{equation}}
\newcommand{\eeq}{\end{equation}}

\begin{document}
\bibliographystyle{amsplain}
\eqnsection

\title[Greedy lattice
animals]{Greedy lattice animals: Geometry and Criticality \\ (with an Appendix)}
\author{Alan Hammond}\thanks{This work was supported in part by NSF
grant DMS-0071448}
\date{21st November, 2004}
\maketitle

\begin{abstract}
Assign to each site of the integer lattice $\Zd$ a real score,
sampled according to the same distribution $F$, independently of
the choices made at all other sites. A lattice animal is a
finite connected set of sites, with its weight being the sum of the scores
at its sites. Let $N_n$ be the maximal weight of those
lattice animals of size $n$ that contain the origin. Denote by $N$ the
almost sure finite constant limit of $n^{-1} N_n$, which exists
under a mild condition on the positive tail of $F$.
We study certain geometrical aspects of the lattice animal
with maximal weight among those contained in an $n$-box where $n$ is large,
both in the supercritical phase where $N > 0$, and in the
critical case where $N = 0$.
\end{abstract}

\section{Introduction}

In this paper, a {\it lattice animal} is a connected set $\xi$ of
sites in $\Zd$, where $d \geq 2$. To each $\xi$ in the set $\LA$
of all lattice animals assign the random {\it weight}
$S(\xi) = \sum_{\bv \in \xi} X_\bv$, where $\{ X_\bv : \bv \in \Zd\}$ are
independent random variables, each
having a common distribution $F$. With $|\xi|$
denoting the number of sites in a set $\xi \subset \Zd$, a
{\it greedy lattice animal of size $n$} is a set $\xi \in \LA$
that contains the origin $\orig$ with $|\xi|=n$ and whose weight $N_n$
is maximal among all such sets.
The study of greedy lattice animals was begun by \cite{GLAone,GLAtwo}.
The authors of \cite{GLAone}
present some optimization problems that motivate the definition of
greedy lattice animals.

It is shown in \cite[Theorem 1]{GLAtwo}
that $n^{-1} N_n$ converges almost surely and in $L^1$ to a
non-random finite constant $N$, in the case that the quantities
$X_\bv$ are non-negative and
\begin{equation}
\label{eqntwopointone}
\int_1^\infty x^d (\log x)^{c} dF(x) < \infty \;
\textrm{for some $c > d$}.
\end{equation}
The same conclusion is derived in
\cite[Theorem 1.1]{Martin} for non-negative $X_\bv$ under the
slightly weaker condition that
\begin{equation}\label{martins}
\int_0^\infty (1-F(x))^{1/d} dx < \infty
\end{equation}
(which in particular holds whenever $c>d-1$ in \req{eqntwopointone}).
By a subtle and involved argument, \cite[Theorem 2.1]{GLAthree} extends the almost sure convergence
$n^{-1} N_n \to N$ to any real-valued $X_\bv$ whose distribution
satisfies \req{eqntwopointone}. 
The condition \req{martins} is almost optimal,
as $\limsup_n n^{-1} N_n = \infty$ whenever $\int_0^\infty x^d dF(x) = \infty$
(see \cite[Theorem 2.2]{GLAthree}). 

We note in passing that
a related object, the greedy lattice path of size $n$, in which the
space $\LA$ is replaced by the collection of finite self-avoiding
paths, is studied in
\cite{GLAtwo,Martin}. These papers prove that
the corresponding normalized weights $n^{-1} M_n$ converge
to a non-random finite constant $M \leq N$
(subject to the same non-negativity and tail conditions on $X_\bv$).
It is further shown in
\cite{Lee} that $M=N$ only when $X_\orig$ is of bounded support and the
probability that $X_\orig$ equals the right end point of its support
is at least the critical probability $p_c$ for site percolation on $\Zd$.

One of the central themes in the study of greedy lattice animals is
the phase transition that the model undergoes as the constant $N$
changes from being negative to positive. It is true that
independent site percolation, obtained by taking
$X_\bv \in \{ -\infty, 1 \}$, is excluded from the theory. 
(Indeed, in the supercritical phase, that is, when $\PPP(X_\orig =1)>p_c$,
the limit $N$ of $n^{-1} N_n$ is 
a non-degenerate random variable on $\{-\infty,1\}$,
with $\{N=1\}$ being the event that $\orig$ is in the infinite
open percolation cluster.) It is helpful however to think of the
natural objects of study in the theory of greedy lattice animals as
counterparts of well-known objects in percolation. We will mention
some of these parallels, as well as comparing and contrasting
percolation with the current framework, throughout this Introduction.

Let $\LA_{C}$ denote the collection of lattice animals contained in
a set $C \subset \Zd$. For positive integers $n$, let
$\Ba_{\bv,n}= \bv+\{0,\ldots,n-1\}^d$ denote
the $n$-box shifted by the vector $\bv \in \Zd$, with $\Ba_n=\Ba_{\orig,n}$.
The paper \cite{GLAthree} studies the
limiting growth of the weight of the {\it greedy lattice animal in the
$n$-box},
\[
G_n := \max \{ S(\xi) \, : \, \xi \in \LA_{\Ba_n} \} \,,
\]
and its size
\[
L_n := \min \{ |\xi| \,:\, \xi \in \LA_{\Ba_n} \; {\rm and} \; S(\xi) = G_n  \},
\]
the minimum being taken to break ties in the case where $F$ is not atomless.

In the percolation model ($X_\bv \in \{ -\infty,1\}$),
there is a transition for the quantity $G_n=L_n$
from $O(\log n)$ at $p<p_c$ to $O(n^d)$ at $p>p_c$.
This transition is analogous to the emergence of a
giant component in the random graph model $G(n,p)$
for $p=c/n$ at $c=1$. It is shown in
\cite[Theorems 3.1 and 3.2]{GLAthree} that a similar transition
occurs for any proper distribution $F$ satisfying \req{martins},
that is, $n^{-1} G_n$ is almost surely bounded in $n$ if $N<0$,
while for some constant $c \in (0,1)$ almost surely
$n^{-d} G_n \in (c,c^{-1})$ for all $n$ large enough.

Our main goals here are to understand more fully the transition of the
weight and size of the greedy lattice animal in the $n$-box, its shape
and the behaviour at criticality, that is, when $N=0$.

Our first result sheds some light on the geometry of the greedy lattice
animal in a large $n$-box in the supercritical case. We define
an $\ell$-box percolation process of parameter $p \in [0,1]$
to be the random collection of disjoint $\ell$-boxes
$\big\{ \Ba_{\ell \ba,\ell}: \ba \in P \big\}$, 
where $P$ is the collection of open sites for an independent
site percolation in $\Zd$ where each site is open with probability $p$.
\begin{theorem}\label{thmjkl}
Let $F$ be a distribution
satisfying \req{martins} for which $N > 0$. For any
$\eps > 0$ there exist $C,\ell \in \mathbb{N}$ and
an $\ell$-box percolation 
$\big\{ \Ba_{\ell \bold a,\ell} : \bold a \in P \big\}$
of parameter at least $1 - \eps$ such that for all $n$ sufficiently large,
each greedy lattice animal in the $n$-box $\Ba_n$ intersects every $\ell$-box
from the largest connected component of
$\big\{ \Ba_{\ell \bold a,\ell} : \bold a \in P, \, \Ba_{\ell \bold a,\ell} \subseteq \{ C \ell, \ldots, n -
1 - C\ell \}^d \big\}$.
\end{theorem}
Theorem \ref{thmjkl} 
implies that $L_n$ exceeds the number of $\ell$-boxes in the largest
connected component to which the theorem refers. 
Applying some well-known facts about the supercritical phase of
percolation, that will later be stated in Lemma \ref{lempva},
the limiting fraction of the $n$-box occupied by
the corresponding greedy lattice animal,
\begin{equation}\label{ldef}
L = \liminf_{n \to \infty} n^{-d} L_n \, ,
\end{equation}
is therefore bounded away from zero when $N > 0$. In this way, 
the theorem 
removes the restriction of exponentially decaying positive tail
of $X_{\orig}$ under which this is proved in \cite[Theorem
4.4]{GLAthree}.
Theorem \ref{thmjkl} shows how pervasive any greedy lattice animal in
a large box must be: it reaches into all but a small fraction of the
array of $\ell$-boxes.  
The limiting density $\lim_{n}{n^{-d}{\vert \gamma_n \vert}}$ 
of the largest cluster $\ga_n$ in $\Ba_n$ of a
supercritical percolation is the density of the unique infinite cluster, $\theta(p) =
\PPP(\vert C(\orig) \vert = \infty)$, by 
\cite[Lemma 7.89]{Gr}. For this reason, 
the counterpart in the
framework of greedy lattice animals of the density of the infinite
cluster 
is the limiting
fraction $L$. At least in principle, this quantity may be random.
Our next result
advances the treatment of
\cite[Theorem 3.2]{GLAthree} by resolving the corresponding question
for the quantity $G$.
\begin{theorem}\label{pthmthree}
For any distribution $F$ that satisfies condition
\req{martins}, there exists a non-random finite constant $G$
such that almost surely
\[
G= \lim_{n \to \infty} n^{-d} G_n.
\]
\end{theorem}
In common with \cite{GLAthree}, the assumption that $F$ may have an
arbitrary negative tail has created the need for more intricate
techniques in the proof of Theorem \ref{thmjkl}, 
and still more so in that of Theorem
 \ref{pthmthree}, than those
that would work were these results to suppose conditions on the
negative tail.  
We next prove a relation between $L$ and the constants $G$ and $N$.
\begin{theorem}\label{pthmsix}
If $F$ is a ditribution satisfying condition \req{martins} for which
$N>0$, 
then the inequality $G \leq L N$ holds almost surely.
\end{theorem}

Given Theorem \ref{thmjkl}, 
the proof of Theorem \ref{pthmsix} is simple. 
Let $\xi_n$ be
a greedy lattice animal in the box $\Ba_n$ for which $\vert \xi_n
\vert = L_n$. If $\xi_n$ happens to contain $\orig$, then 
\begin{equation}\label{gkk}
N_{L_n} \geq S(\xi_n) = G_n. 
\end{equation}
The quantity $N_{L_n}$
behaves like $N L_n$ for high values of $n$ by  \cite[Theorem
2.1]{GLAthree}, from which the inequality $G \leq L N$ follows. 
In general, of
course, the origin may not lie in $\xi_n$. However, it follows
from Theorem \ref{thmjkl} that $\xi_n$ reaches to within a distance
of $\orig$ that is bounded above, uniformly in $n$. This means that
\req{gkk} holds up to a constant term, implying Theorem \ref{pthmsix}.

Consider a distribution $F$ for which $N=0$. We may apply  
Theorem \ref{pthmsix} in the supercritical case of the law $F_{\eps}$,
defined as the distribution 
of the random variable  $X_{\orig} + \eps$, where $X_{\orig}$ has the law
$F$, and where $\eps>0$ is arbitrarily small. It readily follows that 
$G=0$ when $N=0$. 
The authors of \cite{GLAthree} comment that they do not address the
critical case.
The next theorem does so, providing a more precise estimate on the growth of
$G_n$ at criticality than the statement that $G=0$ in this case.
\begin{theorem}\label{thm-crit}
Suppose $F$ satisfies \req{martins} and is such that $N=0$.
Then, for any $c > d/(d-1)$, we have that almost surely
\begin{equation}\label{crit-bd}
\lim_{n \to \infty} (\log n)^{-c} n^{-1} G_n = 0 \,.
\end{equation}
\end{theorem}
The result stands in contrast to that valid for critical percolation
in $\mathbb{Z}^2$, for which \cite{keszha} proves that 
the size of the largest open cluster in the
box $\Ba_n$ grows at a rate exceeding $n^{1 + \delta}$, for some
$\delta > 0$. Theorem \ref{thm-crit} is an optimal result up to
logarthimic corrections, because for choices of $F$ that come close to
violating the positive tail condition \req{martins}, the maximum
weight of sites in $\Ba_n$ is typically of power order $n$. Indeed,
for any $\al > 0$, $F(x) = 1 - x^{-d}(\log x)^{-d(1+\al)}$ satisfies
\req{martins}, with 
\begin{equation}\label{fgb}
 \lim_{n \to \infty}{\PPP \Big( \max_{\bold v \in \Ba_n}{X_{\bold v}}
 > n (\log n)^{-1-\al} \Big)} = 1,
\end{equation}
where $X_{\orig}$ has law $F$.
The random variables given by $Y_{\bold v}(\la) = X_{\bold v} 1 \! \! 1 \{
X_{\bold v} > \la \}$ satisfy \req{fgb} for each $\la \in
\mathbb{R}$. Provided that $\la$ is high enough, the corresponding
value of $N$ may be zero or even negative. Thus, the bound in Theorem
\ref{thm-crit} cannot be improved by more than a logarthmic correction. 
Nonetheless, it may be that for some less contrived choices of $F$ for
which $N=0$, the growth rate of $G_n$ is sublinear. 
 
Finally, we present some results that arise from considering `greedy'
lattice animals that are constrained to occupy a given fraction of the
sites of a large box.

For $n \in \mathbb{N}$ and $\alpha \in (0,1)$, let
\[
\tilde{G}_n (\al) := \max \{ S(\xi) \, : \, \xi \in \LA_{\Ba_n}, \,
| \xi | = \lfloor n^d \alpha \rfloor \} \,,
\]
denote the maximal weight among lattice animals of specified size
that are also contained inside the $n$-box. 
The constant $N$ can be obtained as the limit of the
weight per site of low density maximal weight animals in a large box:  
\begin{prop}\label{thm-conc}
Suppose that $F$ satisfies \req{martins}. Then
\[
\tilde{G} (\al) = \lim_{n \to \infty} n^{-d} \tilde{G}_n(\al)  
\]
exists almost surely for each $\al \in (0,1)$, with
$\tilde{G}: (0,1) \to \mathbb{R}$ a concave, non-random function.
If further $X_\orig$ is bounded below then
$\al^{-1} \tilde{G}(\al) \to N$ as $\al \downarrow 0$.
\end{prop}
Proposition \ref{thm-conc} has the following consequence.
\begin{corollary}\label{webcor}
Suppose that $F$ satisifes \req{martins}, and is such that $X_{\orig}$
is bounded below. Then $G \leq LN$.
\end{corollary}
Indeed, Amir Dembo proposed the method of proof of Proposition
\ref{thm-conc} as a means of deriving the inequality $G \leq
LN$. Of course, we have derived this result by other means.
We note however that, unlike Theorem \ref{pthmsix}, 
Corollary \ref{webcor} does not require that $N>0$. 
We know from \cite{GLAthree} that $G=0$ when $N<0$ and $G>0$ when
$N>0$. The upshot of Corollary \ref{webcor} for the case when $N <
0$ is therefore that $L=0$, provided that
$X_\orig$ is bounded below. 
(In the case where  $X_\orig$ has exponentially decaying
positive tail and $N<0$,  
\cite[Theorems 4.3 and 4.4]{GLAthree} prove that $G_n = O(\log
n)$ and $L_n=O(\log n)$, similarly to the case of percolation.)

If the distribution $F$ is chosen so that $N=0$ and the function
$\tilde{G}:(0,1) \to \mathbb{R}$ is strictly concave, this function attains
its unique maximum at $0$, so that $\tilde{G}(0)=G = \lim_n n^{-d} G_n =
\lim_n n^{-d} \tilde{G}_n(n^{-d}L_n)$ implies that $L=0$. The corresponding
statement for percolation is that $\theta(p_c)=0$, which amounts to
the absence of any infinite cluster at the critical value. 
One analogue of continuity of the percolation probability is trivially
false: if $X_{\orig}$ is equal to $\eps$ almost surely, then $L$ is
the almost sure constant $1$, whatever the value of $\eps>0$. However,
the map $F \mapsto G[F]$ does not have a jump discontinuity as the law
$F$ is increased through choices for which $N =0$. This follows from
\cite{Leetwo}, which proves, under uniform
stochastic dominance and moment assumptions, that $F \mapsto N[F]$  is 
continuous with respect to weak convergence of measures, and the upper
bound $G \leq L N \leq N$ asserted by Theorem \ref{pthmsix}.

\noindent{\bf The global geometry of the greedy lattice animal: two
examples.}

Theorem \ref{thmjkl} prompts the question: if $N > 0$,
how closely does the geometry of a greedy lattice animal in a large
box resemble that of the infinite component of a supercritical site
percolation? Two examples illustrate how the answer depends on the
choice of the distribution $F$. 

In the first example, 
$X_{\bv} \in \{-\la,1\}$ with $\PPP(X_{\orig}=1) = 1 -
\PPP(X_{\orig} = -\la ) = p$ for a pair $(p,\la) \in
(0,1)\times[0,\infty)$ for which $N > 0$. If the parameter $p>p_c$ is
supercritical for site percolation, then a greedy lattice animal in
$\Ba_n$ will, for $n$ large enough, contain the largest connected
component $\Gamma$ of $\{\bv \in \Ba_n: X_{\bv} = 1\}$. Moreover, a smaller
cluster $\gamma$ of one-valued sites in $\Ba_n$ would lie in the
greedy animal by forming a path into
$\Gamma$ unless it is isolated from $\Gamma$ by a region of
$-\la$-valued sites which requires a path of length about $\vert \ga
\vert/\la$ sites to cross. Connected sets of one-valued sites are
therefore much more prone to be part of the greedy lattice animal  in
a large box than connected sets of open sites are liable to 
form part of the largest connected component of open
sites in the same box for a supercritical percolation. This means that, in this case, 
the global
geometry of a greedy lattice animal is at least as connected as that
of the largest component of a percolation.

This choice of law for $X_{\orig}$ has in fact been studied previously.   
 The behaviour of the greedy lattice path in this
model is
closely related (by setting $\la = \rho/(1 + \rho)$) 
to the model of $\rho$-percolation, introduced in 
\cite{man}. Taking $\la > 0$ fixed, \cite{Leeone} explores the
 behaviour of $N[p]$ as
$p \downarrow 0$ and that of $1-N[p]$ as $p \uparrow p_c$.

As a contrast, consider the case where  $X_{\bv} \in \{-1,\la\}$
with $\PPP(X_{\orig} = -1) = 1 - \PPP(X_{\orig} = \la) = p$, with
$\la$ high and $p$ close to one. Supposing that for a greedy lattice
animal $\ga$ in a large box $\Ba_n$, the collection $\ga \cap \{\bv \in
\Ba_n:X_{\bv} = \la \}$ is some given set $\phi \subseteq \Ba_n$,
then $\ga$ will minimize the size of the set of connection costs 
$\{\bv \in \Ba_n:X_{\bv} = - 1 \}$ subject to joining together the
sites of $\phi$ by paths in $\Ba_n$. The choice for the set $\phi$
is made by then optimizing over subsets of $\{\bv \in
\Ba_n:X_{\bv} = \la \}$. 
Thinking of $\la$-valued sites as cities that benefit from travel between them
but must pay some given cost per unit distance of connecting road (or
$-1$-valued site), the greedy lattice animal is the network of cities
that renders the greatest benefit over road cost.  
As such, the greedy animal in this
case is closely related to the 
tree with minimal total edge length subject
to the constraint that its vertices contain some fixed fraction of the
points of a constant rate Poisson process in a large box in $\mathbb{R}^d$,
with the edges being line segments between such points.
The global geometry of this object would seem to be starkly different
from that of the largest connected component of a supercritical
percolation in a large box, having a much higher graphical distance
between a typical pair of distant points.

We mention that, in the current case, we may interpret Theorems \ref{thm-crit}
and \ref{thmjkl} as
asserting that if the benefit $\la$ per city is high enough that there exists a network of cities in $\Ba_n$ whose collective
benefit exceeds road cost by an amount that is super-linear in $n$,
then the optimal network in fact comprises a positive proportion of all the
cities in $\Ba_n$. 

We remark also that the greedy lattice path for this choice of law $F$ 
corresponds to the variant of the travelling salesman problem for
the Poisson  process of points, where the salesman need only visit a high
but fixed fraction of the points in a large box.

\noindent{\bf Organization.}

In Section \ref{sectwo}, we firstly define notations and prove some
lemmas that will be of use throughout
the paper, and then give the proof of Theorem \ref{thmjkl}. The key
to this proof is Lemma \ref{lemmainr}, which shows that it is probable that
large boxes, of sidelength $\ell$, contain weighty lattice animals that may readily be joined
to moderately sized animals in their surroundings. Any greedy lattice
animal in a much larger box $\Gamma$ fails to avoid most such animals in the
array of $\ell$-boxes that lie in $\Gamma$, because
each of these animals is liable to increase the weight of any animal
that runs nearby by joining up with it.
 
In Section \ref{secthree}, we prove 
 Theorem \ref{pthmthree}, beginning with an outline of the argument.
The proof of  Theorem \ref{pthmsix} is given in Section \ref{secold}.

In Section \ref{secfour}, we prove Theorem \ref{thm-crit} by
showing that the negation of \req{crit-bd} implies that $N > 0$:
animals that witness the violation of the bound \req{crit-bd} at
finite values of $n$ are concatenated
by reasonably short paths the weight of whose sites is uniformly
bounded below. All sufficiently large boxes are therefore highly
likely to contain animals whose weight is a high multiple of the
sidelength of the box. A further argument which involves concatenating
these animals establishes the conclusion that $N>0$.   

In the Appendix, we prove 
Proposition \ref{thm-conc} 
and then apply it to deduce Corollary \ref{webcor}.

\noindent{{\bf Remark.}}
 
The results of this paper have been stated with the positive tail
condition \req{martins} being used. 
We will make use of some results of
\cite{GLAthree}, which are stated using the slightly stronger condition
\req{eqntwopointone}. Each of these proofs is valid when   
\req{eqntwopointone} is replaced by  \req{martins}, as explained in
the note on page 207 of 
\cite{GLAthree}.

\section{The geometry of the maximal weight animal: proof of Theorem \ref{thmjkl}}\label{sectwo}

We shall examine in this section some aspects of
the geometry of the greedy lattice animals in the $n$-box $\Ba_n$ 
for large $n$,
in the case where $N>0$. We show below that such animals
intersect well the largest connected component in the $n$-box for
supercritical $\ell$-box percolation, provided that $\ell$ is some
fixed large value and $n \geq n(\ell)$ is high enough.
The proof of Theorem \ref{thmjkl}, in common with several later proofs, will require numerous
definitions and lemmas, which we now state and prove. 
\begin{defin}\label{defnfour}
Let $G_{\Ba}$ and $L_{\Ba}$ denote the weight and size
of a greedy lattice animal of minimal size in a given $m$-box
$\Ba=\Ba_{\bold x,m}$. For
$\lambda \in \mathbb{R}$, say that a site $\bv \in \mathbb{Z}^d$
is $\lambda$-white (or white, for short) if $X_\bv \geq - \lambda$.
The set of $\la$-white paths is a percolation: we will define numerous
site percolations on $\Zd$, so that each percolation process is an
independent site
percolation, unless stated otherwise, and so that $p_c =p_c(d)$ denotes the
critical value for site percolation in $\Zd$. 
The minimal length among all $\lambda$-white paths
in $\Zd$ from some white site in $B$ to some white site
in $A$ is denoted by $D(B,A)$ (or by $D(\bv,A)$ or $D(\bv,\bu)$, in
the case that 
$B=\{\bv\}$ and possibly $A=\{\bu\}$).
We write $B \wcon A$ (in $C$) in case such a path (in $C$) exists. 
We also write $\ell_{\infty}(\bv,A) = \inf_{\bu \in A}{\| \bv
- \bu \|}$ for the minimal sup-norm distance from $\bv$ to a
site in $A$.
Further, for an $m$-box $\Ba=\Ba_{\bold x,m}$ and for any $q \in \mathbb{N}$, set
\begin{displaymath}
\Ba [q] = \bigcup_{0 \leq \| \ba \| \leq q} \, \Ba_{\bold x + m \ba,m},
\end{displaymath}
noting that $\Ba \subseteq \Ba[q]$. 
\end{defin}
\begin{defin}\label{defng}
The boundary of $A \subset \Zd$, denoted $\partial A$, is the collection
of sites $\bu \notin A$ adjacent in $\Zd$ to some $\bv \in A$. For $B
\subseteq \Zd$, the $B$-boundary $\partial_{B} A$ of $A \subseteq B$
is given by $B \cap \partial A$. 
The {\it white cluster} $\WW (\bv)$ of
$\bv \in \Zd$ is the collection of sites $\bu$ such that
$\bu \wcon \bv$ (in particular $\WW(\bv)$ is empty in case $\bv$ is
black). In an analogous
manner, we define the {\it black cluster} of $\bv$, denoted by $\BB
(\bv)$. 

Let $\LL$ be the graph with vertex
set $\Zd$ and with an edge between any pair $\bu \neq \bv \in \Zd$ with
$\| \bu - \bv \| = \max_{1 \leq i \leq d} | u(i) - v(i) | = 1$.
We use the notation 
$\WW_{\LL}(\bv)$ and $\BB_{\LL} (\bv)$ to denote 
the white and black clusters of the site
$\bv$ with respect to the graph $\LL$. 
For $\bv \in \Ba_n$, we write
$\WW_n (\bv)$ or $\BB_n(\bv)$ for the white or black clusters of $\bv$ in
the graph $\Ba_n$, and $\WW_{n,\LL}(\bv)$ or $\BB_{n,\LL}(\bv)$ for these
clusters in the induced subgraph of $\LL$ with vertex set $\Ba_n$.
By `path' or `$\LL$-path', we mean a finite self-avoiding path in
the nearest neighbour or $\LL$-topology.
We will occasionally write `$\Zd$-path' to emphasise that the
nearest-neighbour topology is being used.      

The set $C$ separates $A,B \subseteq \mathbb{Z}^d$ (in $D$) if any
path (in $D$) from
$A$ to $B$ intersects $C$. If each such path intersects $C$ at a
location not lying in $A \cup B$, then we say that $C$ properly
separates $A$ and $B$. The set $C$ separates $A$ from infinity if
any infinite path from $A$ intersects $C$.

For $C \subseteq \Zd$, the exterior boundary of $C$ is given by
\begin{eqnarray}
 \partial_{ext}(C) & = &\Big\{ \bold v \in \Zd : \bold v \, \, \textrm{is adjacent in
 $\LL$ to some $\bold w \in C$,} \label{extbdy} \\
  & &  \qquad \textrm{$\exists$ a path from $\infty$ to $\bold v$ disjoint
  from $C$} \Big\}. \nonumber
\end{eqnarray}
For $C \subseteq \Ba_n$ and $\bold x \in \Ba_n \setminus C$, 
the boundary of $C$ visible from $\bold x$ in
$\Ba_n$ is given by
\begin{eqnarray}
 \partial_{vis(\bold x),n}(C) & = &\Big\{ \bold v \in \Ba_n : \bold v
 \,\, 
\textrm{is adjacent in
 $\LL$ to some $\bold w \in C$,} \label{partbdy} \\
  & &  \qquad \textrm{$\exists$ a path in $\Ba_n$ from $\bold x$ to
  $\bold v$ disjoint
  from $C$} \Big\}. \nonumber
\end{eqnarray}
\end{defin}
\begin{lemma}\label{lemtop}
If $C \subseteq \Zd$ 
is finite and $\LL$-connected, then $\partial_{ext}(C)$ is
$\Zd$-connected. If $C \subseteq \Ba_n$ is $\LL$-connected, 
with $\bold x \in \Ba_n$ and 
$\bold x \not\in
C$, then  $\partial_{vis(\bold x),n}(C)$ is $\Zd$-connected in $\Ba_n$.
\end{lemma}
\noindent{\bf Proof:}  
 The first part of the lemma is \cite[Lemma
2.23]{Kesten}. We prove the second part by altering Kesten's proof. 
The topological setting is that of  a closed $d$-ball instead of
$\mathbb{R}^d$, making the changes more than
merely notational. We write 
\[
\overline{U} = \Big\{ x \in \mathbb{R}^d: \vert x(i) \vert \leq 1/2,
\, i \in \big\{1,\ldots,d\big\} \Big\}
\]
and 
\[
U^{\eps} = \Big\{ x \in \mathbb{R}^d: \vert x(i) \vert < 1/2 + \eps, \, i \in \big\{1,\ldots,d\big\} \Big\}
\]
for some $\eps \in(0,1/8]$, and
we set $N = \bigcup_{\bv \in C}(\bv + U^{\eps}) \, \cap \, [0,n-1]^d$.
Kesten first proves, by invoking Alexander and Poincar\'{e} duality, that, if $N' \subseteq \mathbb{R}^d$ is bounded and
connected with its topological boundary $\partial N'$ a topological
$(d-1)$-manifold, then the set
\begin{eqnarray}
 \partial_{ext}(N') & := &  \Big\{ y \in \partial{N'} : \, \exists \, \,
 \textrm{a continuous path $p:(0,\infty)\to \mathbb{R}^d \setminus
 N'$}, \label{kestset} \\ 
 & & \qquad \quad  p(t) \to y  \, \, \textrm{as $t \to 0$}, \, \,
 \vert p(t) \vert \to \infty \, \, \textrm{as $t \to \infty$} \Big\} \nonumber
\end{eqnarray}
is path connected. 
We define for any sets $F,O \subseteq \mathbb{R}^d$ with $O$ open, 
and for each point $y \in F \setminus {\rm clos}(O)$,
\begin{eqnarray}
 \partial_{vis(x),F}(O) & = & \big\{ y \in F \cap \partial O : \, \exists \, \,
 \textrm{continuous path $p:[0,T]\to F$:} \nonumber \\
 & & \qquad \qquad p(0)=x, p(T)=y, \, p[0,T] \cap O = \emptyset \big\} . \nonumber
\end{eqnarray} 
 (Note that the use of the
$\partial_{ext}$ and $\partial_{vis}$ notations causes no conflict
with the discrete case \req{extbdy} or \req{partbdy}, 
because we are considering $O \subseteq
\mathbb{R}^d$).
The analogue of the path-connectedness of the set \req{kestset} that we require is that, for
each $x \in [0,n-1]^d \setminus {\rm clos}(N)$,
\begin{equation}\label{hereset}
 \textrm{the set} \, \, \partial_{vis(x),[0,n-1]^d}(N) \, \, 
\textrm{is $[0,n-1]^d$-path connected.} 
\end{equation}
To establish \req{hereset}, let $x \in  [0,n-1]^d \setminus {\rm
clos}(N)$ be fixed.
We assume in the first instance 
that  $N \cap \partial([0,n-1]^d) \not= \emptyset$,
and will
prove \req{hereset} in this instance 
by reducing to the $\mathbb{R}^d$ case by a
doubling trick. 
We let closed $d$-balls $B_d(1)$ and $B_d(2)$
 denote two
homeomorphic images of $[0,n-1]^d$, writing $x$ and $N$ 
for the images in $B_d(1)$ of $x$ and $N$ 
by a harmless abuse of notation. 
We glue the two balls by identifying
their boundaries to form $\Gamma = B_d(1) \cup {}_{S^{d-1}} B_d(2)$,
which is a homeomorphic image of the sphere $S^d$. 
Let $\phi: \Gamma \to \Gamma$ map each point in $B_d(i)$ to the
corresponding one in $B_d(2-i)$ for $i \in \{1,2\}$ (so that $\phi$  fixes the common boundary of
$B_d(1)$ and $B _d(2)$), and let $\hat{\phi}:\Gamma \to
\Gamma$ be given by 
\begin{eqnarray}
 \hat{\phi}(y) & =  & y \, \, \textrm{if $y \in B_d(1)$} \nonumber \\
 & = & \phi(y) \, \, \textrm{if $y \not\in B_d(1)$.} \nonumber 
\end{eqnarray} 
Define  $\hat{N} = N \cup_Z \phi(N)$, where $Z = N \cap \partial
B_d(1)$ is the part of the boundary of $B_d(1)$ along which $N$ and
its mirror $\phi(N)$ are glued.
If $\chi: \Gamma \setminus \{x\} \to \mathbb{R}^d$ is a homeomorphism, 
then $\partial_{vis(x),\Gamma}(\hat{N})=
\chi^{-1}\big(\partial_{ext}(\chi(\hat{N}))\big)$.
We now check that the conditions on  $\chi(\hat{N})$ that permit us to
take apply \req{kestset} with the choice $ N' =  \chi(\hat{N})$.
Note that $\chi(\hat{N}) \subseteq \mathbb{R}^d$ is bounded, because
$x \not\in {\rm clos}(\hat{N})$, and that $\chi(\hat{N})$ is connected because
the fact that $N \cap \partial [0,n-1]^d \not= \emptyset$ 
implies that $\hat{N}$ is connected. 
We must also show that  $\partial_{\mathbb{R}^d}(\chi(\hat{N}))$
is a topological $(d-1)$-manifold.
Kesten proved that
\begin{equation}\label{topdman}
\textrm{$\partial_{\mathbb{R}^d}  \bigcup_{\bv \in C}(\bv + U^{\eps})$
is a topological $(d-1)$-manifold}.
\end{equation}
Recalling that we consider $N$ as a subset of $B_d(1)$,
it follows easily from the definition of $N$ that 
$Z = N \cap \partial B_d(1)$
is a submanifold of $\partial B_d(1)$ with boundary. 
This is a
sufficient condition for the doubled object $\hat{N}$ to be a
topological $d$-manifold and furthermore, for its boundary
in $\Ga$ to be given by 
\begin{equation}\label{vag}
 \partial_{\Ga} ( \hat{N} ) = \big( \partial_{B_d(1)} N \setminus
 \textrm{\.{Z}} \big) \cup_{\partial_{\partial B_d(1)} Z} \big( \partial_{B_d(2)} \phi(N) \setminus \textrm{\.{Z}} \big),
\end{equation}
where $\textrm{\.{Z}} = Z \setminus \partial Z$,
with $\partial Z$ denoting the boundary in $\partial B_d(1)$ of $Z$  
(in \req{vag}, $\textrm{\.{Z}}$ is
respectively embedded in $B_d(1)$ or $B_d(2)$). See \cite[Chapter 8.2]{hirsch}
for a discussion of gluing of such manifolds. We must now check that
the right-hand-side of \req{vag} is a $d-1$-manifold. We will do so by
applying the same sufficient condition already mentioned, which in
this case means that $\partial Z$ is a $d-2$
manifold. We show in fact that each of the finitely many connected
components of $\partial Z$ is  a $d-2$ manifold. To this
end, let $Z'$ denote one of these components. Note that $U
\cap Z'$ is homeomorphic to an open subset of a set of the form
appearing in \req{topdman}, with $d$ replaced by $d-1$, provided that
$U$ is an open set 
lying in a part of $\partial B_d(1)$ corresponding to an open
face of $[0,n-1]^d$. If $U$ is a neighbourhood of a point in the
boundary of at least two such faces, then an initial homeomorphism is
required to flatten $\partial B_d(1) \cap U$. The image of $U \cap Z'$
under this map is then homeomorphic 
to the same type of set as in the other case. We thus learn from
\req{topdman} (applied with $d$ replaced by $d-1$) that 
$Z'$ is indeed a manifold. We find 
that the glued object $\partial_{\Ga} ( \hat{N} )$ is a
$(d-1)$-manifold, so that $\partial_{\mathbb{R}^d} \big(\chi(\hat{N})\big)$ 
is, as well. 

We may apply \req{kestset} to  $\chi(\hat{N})$ as we sought, learning by
doing so that 
$\partial_{ext}(\chi(\hat{N}))$ is path
connected, so that 
$\partial_{vis(x),\Gamma}(\hat{N})$ is $\Gamma$-path connected.  
We now use this fact to verify that the set in \req{hereset} is
$[0,n-1]^d$-path connected. 
To this end, let $y,z \in \partial_{vis(x),B_d(1)}(N)$. As
$N \subseteq B_d(1)$, we may consider $y$ and $z$ as members of
$\partial_{vis(x),\Gamma}(\hat{N})$. We use the $\Gamma$-path
connectedness of this last set
to find a path $p:[0,1] \to \partial_{vis(x),\Gamma}(\hat{N})$ such
that $p(0)=y$ and $p(1)=z$. We claim that the path
$\hat{\phi} \circ p : [0,1] \to B_d(1)$ satisfies
$\hat{\phi} \circ p [0,1] \subseteq \partial_{vis(x),B_d(1)}(N)$. 
Indeed, to any $t \in [0,1]$, let $q:[0,1] \to \Gamma$ denote
a continuous map satisfying $q(0)=x$, $q(1)=p(t)$ and $q[0,1] \subseteq
\Gamma \setminus \hat{N}$ (such a map existing because
$p(t) \in \partial_{vis(x),\Gamma}(\hat{N}$).) The map $\hat{\phi} \circ q
:[0,1] \to B_d(1) \setminus N$
demonstrates that $\hat{\phi} \circ p(t) \in \partial_{vis(x),B_d(1)}(N)$ for each
$t \in [0,1]$. Since $y = (\hat{\phi} \circ p) (0)$ and  
 $z = (\hat{\phi} \circ p) (1)$, we have proved that the set in
 \req{hereset} is
$[0,n-1]^d$-path connected, in the case that  $N \cap
\partial([0,n-1]^d) \not= \emptyset$.

If  $N \cap
\partial([0,n-1]^d) = \emptyset$, we simply define a homeomorphism $\xi:[0,n-1]^d
\setminus \{ x \} \to \mathbb{R}^d$, and note that we may apply
\req{kestset}
with the choice $N' = \xi(N)$ by a similar argument to that which
permitted the choice $N' = \chi(\hat{N})$. We learn that, in this case also, 
the set $\partial_{vis(x),[0,n-1]^d}(N) = \xi^{-1}\big(
\partial_{ext}(\xi(N))\big)$ is $[0,n-1]^d$-path connected, as we
sought.

Secondly, Kesten proves that if $\bold v',\bold v'' \in \Zd$
are connected by a path $\phi:[0,1] \to \mathbb{R}^d \setminus N$, then $\bold v'$ and
$\bold v''$ may be connected by a $\Zd$-path that is disjoint from
$N$, and intersects only such cubes $\bold v + \overline{U}$, $\bold v
\in \Zd$, that also contain a point of $\phi$. Correspondingly, our
pair of sites $\bold v',\bold v''$ lie in $[0,n-1]^d$, and the path
$\phi$ that connects them in $[0,n-1]^d \setminus N$, and we require that the
$\Zd$ path lies in $[0,n-1]^d \setminus N$. For this, the path produced in the
$\mathbb{R}^d$-case does the job.

There are two more steps in Kesten's argument. The use of $[0,n-1]^d$ in
place of $\mathbb{R}^d$ makes no difference to the proofs of these
steps, and we only state the form they take in our case. The third
step asserts that, for any $\bold x \in \Ba_n \setminus C$, and for each
$\bold v
\in \Ba_n$, we have that $\bold v \in \partial_{vis(\bold x),\Ba_n}(C)$
if and only if $\bold v \not\in C$ and $(\bold v + \overline{U}) \cap
\partial_{vis(\bold x),[0,n-1]^d}(N) \not= \emptyset$. The fourth step claims that each pair
$\bold v',\bold v'' \in \partial_{vis(\bold x),\Ba_n}(C)$
can be connected by a path in $[0,n-1]^d$ which lies in $(\bold v' +
\overline{U}) \cap (\bold v'' + \overline{U}) \cup
\partial_{vis(\bold x),[0,n-1]^d}(N)$. Similarly to Kesten, this last
path may be deformed by the procedure in the second step into a path
in $\Ba_n$ that contains only sites of $\partial_{vis(\bold
x),\Ba_n}(C)$.
Thus, $\partial_{vis(\bold
x),\Ba_n}(C)$ is connected in $\Ba_n$, as required. $\Box$     
\begin{lemma}\label{lemmbul}
Suppose that the connected sets $C,D \subseteq \Ba_n$ are disjoint, and that $E
\subseteq \Ba_n$ separates $C$ and $D$ in $\Ba_n$, and is disjoint
from $D$.
Then there exists an $\LL$-connected subset of $E$ that also 
separates $C$ and $D$ in $\Ba_n$. 
Suppose instead that $C$,$D$ and $E$ are subsets of $\Zd$ that satisfy the
same conditions of disjointness, and that $C$ and $D$ are connected.
If $E$ separates $C$ and $D$, then,
similarly, an $\LL$-connected subset of $E$ separates $C$ and $D$.
\end{lemma}
\noindent{\bf Proof:} 
We comment briefly on the content of what is asserted in the lemma,
before proving it. In the case of the first part of the lemma, 
set
\begin{equation}\label{ronem}
 \tilde{C} = \bigcup_{\bv \in C \setminus E}\big\{ \bold x \in \Ba_n :
 \bold x \leftrightarrow \bv  \, \, \textrm{in} \, \, E^c \big\},
\end{equation}
and
\begin{equation}\label{rone}
\hat{C}  = (C \cap E) \, \cup \, \partial_{\Ba_n}{\tilde{C}},
\end{equation}
Note that $\hat{C} \subseteq E$, and 
that $\hat{C}$ separates $C$ and $D$ in $\Ba_n$: indeed, the
first element of $E$ encountered in any path in $\Ba_n$ from $C$ to $D$
lies in $\hat{C}$. From this fact and the disjointness of $D$ and $E$,
it is straightforward that, for 
any $\bold y \in D$, the set
\begin{equation}\label{zas}
 \Big\{ \bold v \in \hat{C} : \bold v
 \, \,  
\textrm{is $\Zd$-adjacent to some $\bold w \in \Ba_n \setminus \hat{C}$,}
 \, \, \textrm{$\exists$ a path in $\Ba_n$ from $\bold y$ to
  $\bold w$ disjoint
  from $\hat{C}$} \Big\}. 
\end{equation}
separates $C$ and $D$ in $\Ba_n$.
A variant of the second part of Lemma \ref{lemtop} asserting that
the set \req{zas} is $\LL$-connected would therefore suffice for our
purpose. We mention that a sketch of an elementary proof of a similar assertion 
appears in the Appendix of \cite{matremy}. The variant might also be obtained
 by re-examining the proof of  \cite[Lemma 2.23]{Kesten}. Instead of
 doing this, we will prove the first assertion of the lemma by finding
 an $\LL$-connected component of $\hat{C}$ that separates $C$ and $D$
 in $\Ba_n$, making direct use of the second part of Lemma
 \ref{lemtop} in the course of the argument .

Let $\big\{ \gamma_1,\ldots, \gamma_k \big\}$
denote the $\LL$-connected components of $\hat{C}$, and let 
\[
\phi_i = \big\{ \bv \in \Ba_n : \gamma_i \, \,
\textrm{separates $\bv$ from $D$ in $\Ba_n$} \big\}.
\]
 We claim that
\begin{equation}\label{rtwo}
 \phi_i \cap \phi_j = \emptyset  \, \, \textrm{implies that either} \,
 \,\phi_i
 \subseteq  \phi_j , \, \, \textrm{or} \, \,  \phi_j
 \subseteq  \phi_i.
\end{equation} 
For $\bold x \in \phi_i \cap \phi_j$, we consider the case that there
exists a path $\tau$ in $\Ba_n$ from $\bold x$ to $D$ that intersects $\gamma_i$
before it intersects $\gamma_j$. We will show that $\phi_i \subseteq
\phi_j$, and, to this end, we let $\bold y \in \phi_i$ be arbitrary. 
Let $\big( \bold
x = \bold x_0 ,
\bold x_1,\ldots, \bold x_{r_1} \big)$ denote the segment of the path $\tau$
until its first intersection with $\gamma_i$ (so that $\bold x_{r_1} \in
\gamma_i$). Let  $\big( \bold y = \bold y_0 ,
\bold y_1,\ldots, \bold y_{r_2} \big)$ denote any path in
$\Ba_n$ from $\bold y$ to
$D$. 
There exists $r_3 \in \big\{0,\ldots,r_2 \big\}$
for which $\bold y_{r_3} \in \gamma_i$, since $\bold y \in
\phi_i$. Note that 
\begin{equation}\label{rthr}
 \Big\{ \bold v \in \Ba_n : \ell_{\infty} (  \bold v , \gamma_i ) \leq 1
 \Big\} \cap \gamma_j = \emptyset \, \, \textrm{for $j \not= i$,}
\end{equation}
because the sets $\gamma_j$ are $\LL$-connected components of a larger
set. There exists an $\LL$-path from $\bold x_{r_1}$ to $\bold y_{r_3}$ in
$\ga_i$. Any consecutive members $\bold u$ and $\bold v$ of this path each lie in
a unit cube contained in $\Ba_n$. We may find a $\Zd$-path in this cube
from $\bold u$ to $\bold v$.
In this way, we may find a path 
$\big( \bold x_{r_1} = \bold z_1,\ldots, \bold z_{r_4} = \bold y_{r_3} \big)$ in $\Ba_n$ 
such that $\ell_{\infty} (
 \bold z_l ,  \gamma_i ) \leq 1$ for  $l \in \big\{ 1, \ldots, r_4 \big\}$, 
and thus, by \req{rthr}, such that $\bold z_l \not\in
\gamma_j$ for these $l$ and for $j \not= i$. Note that
the path $\big( \bold x = \bold x_0,\bold x_1,\ldots \bold x_{r_1} =
\bold z_1,\ldots,\bold z_{r_4} =
\bold y_{r_3},\ldots \bold y_{r_2} = \bold y \big)$ connects $\bold x$ to $D$ in $\Ba_n$. 
Note also that
its initial segment $\big( \bold x_0,\bold x_1,\ldots, \bold x_{r_1} ,
\bold z_2,\ldots \bold z_{r_4}
\big)$ is disjoint from $\ga_j$ by construction. The fact
that $\bold x \in \phi_j$ implies that  $\bold y_r \in \gamma_j$ for some $r \in \big\{ r_3 +
1,\ldots, r_2 \big\}$.
The site $\bold y_r$ lies in the path $\big( \bold y_0,\ldots,\bold y_{r_2}\big)$,
which was chosen to be an arbitrary path from $\bold y$ to $D$.
Thus, $\bold y \in \phi_j$, and $\phi_i \subseteq \phi_j$, as we
claimed. 
The conclusion
$\phi_j \subseteq \phi_i$ arises in the other case, where there exists a path from $\bold x$ to
$D$ that intersects $\gamma_j$ before it intersects
$\gamma_i$. This establishes (\ref{rtwo}).

Note that (\ref{rtwo}) implies that to each $i \in
\big\{1,\ldots,k\big\}$, there is an unique $j \in 
\big\{1,\ldots,k\big\}$ such that $\phi_i \subseteq \phi_j$ and for
which $\phi_j \subseteq \phi_l \implies j = l$. Let $\big\{
j_1,\ldots,j_s \big\}$ denote the collection of $j \in 
\big\{1,\ldots,k\big\}$  arising in this way for some 
 $i \in
\big\{1,\ldots,k\big\}$. Note that, by (\ref{rtwo}), $\phi_{j_i} \cap
 \phi_{j_l} = \emptyset$ for $i,l \in \big\{ 1,\ldots,
s\big\}$ with $i \not= l$. We claim that:
\begin{equation}\label{rfour}
 C \subseteq \bigcup_{i=1}^{s}{\phi_{j_i}}.
\end{equation}
To derive this, suppose that  $\bold x \not\in
\bigcup_{i=1}^{s}{\phi_{j_i}} = \bigcup_{i=1}^{k}{\phi_i}$ for some $\bold x \in \Ba_n$. Let $T = \big(
\bold x = \bold v_0,\ldots \bold v_{r_5}
\big)$ denote an arbitrary path in $\Ba_n$ from $\bold x$ to $D$. If $T \cap
\hat{C} = \emptyset$, then $x \not\in C$, because $\hat{C}$ separates
$C$ and $D$ in $\Ba_n$. We may assume
then that $T \cap \hat{C} \not= \emptyset$, and that $T$
intersects $\gamma_1$ before it intersects any 
$\ga_i$ for $i \in \{ 2,\ldots k\}$. In this case, let $0<r_6 \leq r_7<r_5$ be such that $\bold v_{r_6}$
is the first visit of $T$ to $\ga_1$, while
$\bold v_{r_7}$ is the last such. (That $r_6 > 0$ follows from $\bold x
\not\in \phi_1$, and thus $\bold x \not\in \ga_1$). We next show that 
\begin{equation}\label{rfive}
\bold v_{r_6 - 1}, \bold v_{r_7 +1} \in
 \bigcup_{\bold y \in D} \partial_{vis(\bold y),n}\big( \ga_1 \big): 
\end{equation} 
firstly, each of these two sites is
adjacent to a member of $\ga_1$. Note also that $(\bold v_{r_7
+ 1},\ldots,\bold v_{r_5})$ is a path in $\Ba_n \setminus \ga_1$ 
from $\bold v_{r_7 +
1}$ to $D$. Since $\bold x \not\in
\phi_1$, we may find a path in $\Ba_n$ from $\bold x$ to some site
$\bold y
\in D$ that is disjoint
from $\ga_1$. Concatenating this path to
$\big(\bold v_{r_6 -1},\bold v_{r_6 -2} , \ldots, \bold v_1 \big)$ forms a path in $\Ba_n$
from $\bold v_{r_6 - 1}$ to $\bold y$ that is disjoint from
$\ga_1$. Thus, $\bold v_{r_6 - 1} \in \partial_{vis(\bold y),n}\big(
\ga_1 \big)$, and (\ref{rfive}).

Note that $\partial_{vis(\bold y),n}\big( \ga_1 \big)$ is independent of
$\bold y \in D$ because the connected set $D$ satisfies $D \cap E =
\emptyset$ and thus $D \cap \ga_1 = \emptyset$ 
(because $\ga_1 \subseteq \hat{C}
\subseteq E$). By the second assertion of Lemma \ref{lemtop},
we may choose a path
$\big( \bold v_{r_6 - 1} = \bold u_0,\bold u_1,\ldots,\bold u_{r_8} = \bold v_{r_7 + 1}\big)$ that
lies in $\partial_{vis(\bold y),n}(\gamma_1)$ for any given $\bold
y \in D$.
It follows from $\partial_{vis(\bold y),n}(\ga_1) \subseteq
\big\{ \bold v \in \mathbb{Z}^d : 
 \ell_{\infty} ( \bold v , \ga_1 ) \leq 1 \big\}$ and (\ref{rthr}) that
\[
 \big\{ \bold u_0,\ldots,\bold u_{r_8} \big\} \, \cap \, \bigcup_{i=2}^{k}
 \gamma_i = \emptyset.  
\]
We alter the path $T$ to form
\[
 \big( \bold x = \bold v_0,\ldots,\bold v_{r_6 - 1} = \bold u_0,\ldots,\bold u_{r_8} = \bold v_{r_7 +
 1},\ldots, \bold v_{r_5} \big).
\]
This new path reaches $D$ from $\bold x$, is disjoint from
$\ga_1$ and intersects  $\bigcup_{i=2}^{k}
 \ga_i$ only at points where the path $T$ does. By
 performing alterations that similarly remove the intersections of the
 new path with the other sets 
$\gamma_i$, we produce a path in $\Ba_n$ from $\bold x$ to $D$ that is
disjoint from  $\cup_{i=1}^{k}
 \ga_i = \hat{C}$. Since
 $\hat{C}$ separates $C$ and $D$ as noted after
 (\ref{rone}), $\bold x \not\in C$. We have proved (\ref{rfour}).

We now claim that
\begin{equation}\label{rsix}
  \textrm{there exists $i \in \big\{ 1,\ldots,s\big\}$ for which
  $C \subseteq \phi_{j_i}$.}
\end{equation}
Were this not the case, there would exist by \req{rfour} adjacent sites
$\bold w_1,\bold w_2 \in C$ such that $\bold w_1 \in \phi_{j_{i_1}}$
and  $\bold w_2
\in \phi_{j_{i_2}}$ for indices $i_1,i_2 \in \{ 1,\ldots, s \}$
satisfying 
$i_1 \not= i_2$. 
By (\ref{rthr}), one
of $\bold w_1 \in \gamma_{j_{i_1}}$ and  $\bold w_2 \in
\gamma_{j_{i_2}}$ fails. We may assume that $\bold w_1 \not\in
 \gamma_{j_{i_1}}$. The sets $\phi_{j_i}$ being disjoint, $\bold w_2 \in
\phi_{j_{i_2}}$ implies that there exists a path $\sigma$ 
from $\bold w_2$ to $D$
that is disjoint from $\gamma_{j_{i_1}}$. The fact that $\bold w_1
\not\in \gamma_{j_{i_1}}$ implies that the path formed by
prefixing $\bold w_1$ to $\sigma$ reaches $D$ from $\bold w_1$ and is disjoint
from  $\gamma_{j_{i_1}}$. We have reached the contradiction
that $\bold w_1 \not\in  \phi_{j_{i_1}}$ and have proved (\ref{rsix}). From
(\ref{rsix}) follows the first statement of the lemma, because $\ga_{j_i}$ is
$\LL$-connected. 

The second assertion has the same proof, with the first part of Lemma
\ref{lemtop} being applied, instead of the second part. The notational
changes consist of omitting each reference to `in $\Ba_n$', including in
the term $\partial_{vis(\bold y)}(\cdot)$ (which is independent of
$\bold y \in D$). $\Box$
\begin{lemma}\label{isomp}
Let $B \subseteq \mathbb{Z}^d$ denote the collection of sites that a
finite set $A  \subseteq \mathbb{Z}^d$ separates from infinity. Then
\[
 \big\vert A \big\vert \geq \big\vert B \big\vert^{1 - \frac{1}{d}}.
\]
Let $C,D,E \subseteq \Ba_n$, with $C$ and $D$ connected, and $C \cap D = \emptyset$. Suppose
that $E$ 
properly separates $C$ and
$D$ in $\Ba_n$. Then
\[
  \big\vert E \big\vert \geq   \frac{1}{2d} \min \Big\{ \big\vert C \big\vert^{1 -
  \frac{1}{d}},  \big\vert D \big\vert^{1 - \frac{1}{d}} \Big\}.
\]
\end{lemma}
\noindent{\bf Proof:}
To prove the first part of the lemma, note that, for each 
$i \in \{1,\ldots,d\}$ and $\bold x \in B$, $\big\{ \bold x + n \bold
e_i : n \in \mathbb{Z} \big\} \cap A \not= \emptyset$, where $\big\{ \bold e_i : i \in
\{1,\ldots,d \} \big\}$ denote unit vectors in the directions of the
co-ordinate axes of $\mathbb{R}^d$. It follows that
\begin{equation}\label{ott}
 \vert A \vert \geq \max_{i \in \{1,\ldots,d \} }{\vert {\rm proj}_i (B) \vert},
\end{equation}
where 
\[
 {\rm proj}_i (B) = \big\{ (\bold a_1,\ldots,\bold a_{d-1}) \in
 \mathbb{Z}^{d-1}: \, \exists \, \bold x \in \mathbb{Z} \, \,
 \textrm{such that} \, \, (\bold
 a_1,\ldots,\bold a_{i-1},\bold x,\bold a_i,\ldots,\bold a_{d-1}) \in
 B \big\}.
\]
By the Loomis-Whitney inequality \cite{loomwhit}, 
\begin{equation}\label{rtbc}
  \max_{i \in \{1,\ldots,d \} }{\vert {\rm proj}_i (B) \vert} \geq \vert
  B \vert^{1 - 1/d},
\end{equation}
so that we obtain the first part of the lemma from \req{ott} and
\req{rtbc}.

To begin the proof of the second part of the lemma, note that the
first part of Lemma \ref{lemmbul} permits us to assume that $E$ is
$\LL$-connected. We may also assume that $E \cap (C \cup D) =
\emptyset$,
because   $E \setminus (C \cup D)$ properly separates $C$ and $D$.
Recalling the definition \req{ronem}, note that
\begin{equation}\label{hdf}
 \tilde{C} \cap \tilde{D} = \emptyset.
\end{equation}
This is because, in the other case, we might construct a path in $E^c$
from $C$ to $D$. By \req{hdf}, at least one of the inequalities 
$\vert \tilde{C} \vert \leq \frac{1}{2}n^d$ and 
$\vert \tilde{D} \vert \leq \frac{1}{2}n^d$ holds. We suppose the
former for the time being. 
It follows from \cite[Theorem 19]{Leader} that,
for any $A \subseteq \Ba_n$ for which $\vert A \vert \leq n^d/2$, we
have that
\begin{equation}\label{dxc}
 \big\vert \partial_{\Ba_n} A \big\vert \geq \frac{1}{2d} \vert A\vert^{1 - \frac{1}{d}}.
\end{equation}
Note that $\partial_{\Ba_n}{\tilde{C}} \subseteq E$. Indeed, it was
noted after \req{rone} that $\hat{C} \subseteq E$, and
$\partial_{\Ba_n}{\tilde{C}} = \hat{C}$ in the current case, because
$C \cap E = \emptyset$.
We find that 
\[
 \vert E \vert \geq \vert \partial_{\Ba_n}{\tilde{C}} \vert
\geq \frac{1}{2d} \vert \tilde{C} \vert^{1 - \frac{1}{d}} 
\geq \frac{1}{2d} \vert C \vert^{1 - \frac{1}{d}},
\]
where \req{dxc} with the choice $A = \tilde{C}$ was applied in the second
inequality, this choice being valid because $\vert \tilde{C} \vert
\leq \frac{1}{2}n^d$.  The third inequality follows from $C \subseteq
\tilde{C}$, which is implied by $C \cap E = \emptyset$.

In the case where 
$\vert \tilde{D} \vert \leq \frac{1}{2}n^d$, we deduce that 
$ \vert E \vert 
\geq \frac{1}{2d} \vert D \vert^{1 - \frac{1}{d}}$. This completes the
proof of the second part of the lemma. $\Box$

Having assembled these preliminaries, we now state and prove a lemma
that is central to the proof of Theorem \ref{thmjkl}.
\begin{lemma}\label{lemmainr}
For given constants
$c>0$ and $\la,\rho<\infty$, we say that an $m$-box $\Ba=\Ba_{\bold x,m}$
satisfies condition ${\bf A}_{\la,c}^{\rho}$ if 
there exists $\gas \in \LA_{\Ba}$ such that
$S(\gas)=G_{\Ba} \geq c m^d$, $\vert \gas \vert \geq (\log m)^{\rho} +
1$, and  $D(\ga,\gas) \leq \rho m$ for
all $\ga \in \LA_{\Ba[2]}$ such that $|\ga| \geq (\log m)^\rho$.
We set
$
q_{m,\la,c,\rho} := \PPP(\Ba \ \textrm{satisfies condition}
\ {\bf A}_{\la,c}^{\rho} )$,
which is independent of the choice of the $m$-box $\Ba$.
Suppose that the distribution $F$
satisfies \req{martins} and is such that $N > 0$.
Then, for $\la, \rho<\infty$ sufficiently large
and $c \in (0,\liminf n^{-d}G_n)$, we have that
\[
\lim_{m \to \infty} q_{m,\la,c,\rho}   = 1 .
\]
\end{lemma}
\noindent 
{\bf Remark} Theorem 3.2 of \cite{GLAthree} states in the case that $N
> 0$, there exists a constant $c > 0$ such that $\liminf
n^{-d}G_n > c > 0$ almost surely. \\ 
\noindent
{\bf Proof:} Fixing $\eps>0$,
\req{martins} implies that $\int_0^\infty x^d dF(x)<\infty$
(this can be derived directly or by contrasting \cite[Theorem 2.2]{GLAthree}
with the note in page 207 of that paper). Following
\cite[proof of (3.17)]{GLAone}, this implies that $X_\bv \geq \|\bv\|$
for at most finitely many $\bv \in \Zd$, which implies, 
in view of $N>0$ and \cite[Theorem 3.2]{GLAthree}, that, for some $c_1>0$,
\begin{equation}\label{lmineq}
 L_m \geq c_1 m^{d-1} \, \textrm{for all sufficiently high $m$.}
\end{equation}
Thus, for all $m$ large enough,
\begin{equation}\label{lbdLn}
\PPP( L_\Ba \geq c_1 m^{d-1} ) > 1 - \eps
\end{equation}
for any $m$-box $\Ba$.

For $\la$ such that
\begin{equation}\label{pbound}
p:= \PPP (X_\orig \geq - \la) >  \max \{ p_c,  1 - p_c(\Zd,\LL) \},
\end{equation}
let $\WW$ denote the unique infinite cluster of $\la$-white sites in $\Zd$
(which exists by \cite[Theorem 8.1]{Gr}, since the $\la$-white sites form a supercritical
percolation). Increasing $\la$ as needed, we next show
that for any $\rho>1$, all $m$ large enough and any $m$-box $\Ba$,
\begin{equation}\label{ubdnew}
\PPP \Big( \exists \ga \in \LA_{\Ba[2]}
 \ \textrm{such that} \ |\ga| \geq (\log m)^\rho
 \ \textrm{and} \  \ga \cap \WW = \emptyset \Big) \leq \eps \,.
\end{equation}
For such a $\gamma$ as in \req{ubdnew},
there exists, by choosing $C = \gamma$, $D= \WW$ and $E = \{ \,
\textrm{black sites\}}$ in the second part of Lemma \ref{lemmbul}, 
an $\LL$-connected set
$\hat{\ga}$ of black sites that separates $\ga$ and $\WW$.  
Applying  
the first part of Lemma \ref{isomp}, we find that
\begin{equation}\label{bkj}
  \big\vert \hat{\gamma} \big\vert
  \geq \big( \log m \big)^{\rho (1 - \frac{1}{d})},
\end{equation}
because  $\hat{\ga}$ separates $\ga$ from $\infty$.
We set $A =
A_{\ga}$ according to $A = \inf \big\{ q \geq 2 : \hat{\ga} \cap B[q]
\not= \emptyset \big\}$.
Note that if $A > 2$, then $B[A-1]$ is separated
from $\WW$, and thus from $\infty$, by $\hat{\ga}$: indeed, a path from $B[A-1]$ to $\WW$
disjoint from $\hat{\ga}$ could be extended in $B[A-1]$ to such a path
from $\gamma$ because $B[A-1]\cap \hat{\gamma} = \emptyset$.
By   
Lemma \ref{isomp} again,
\begin{equation}\label{phjk}
 \big\{ A = q \big\} \subseteq \Big\{ \big\vert \hat{\gamma} \big\vert
  \geq \big( (2q-1) m \big)^{d-1} \Big\}, \, \textrm{
for $q \geq 3$, }
\end{equation}
since $\big\vert B[q-1] \big\vert = \big[ (2q- 1)m\big]^d$. Thus,
\begin{eqnarray}
& & \PPP \Big( \exists \ga \in \LA_{\Ba[2]}
 \ \textrm{such that} \ |\ga| \geq (\log m)^\rho
 \ \textrm{and} \  \ga \cap \WW = \emptyset \Big)  \nonumber \\
&\leq &
(5 m)^d \PPP \Big(\big| \BB_{\LL} (\orig) \big| \geq  (\log m)^{\rho(1 - 1/d)} \Big)
\nonumber \\
& & + \sum_{q=3}^{\infty}{\big[ (2 q + 1)m \big]^d  \PPP \Big( \big| \BB_{\LL}
(\orig) \big| \geq  (2 q - 1)^{d-1} m^{d-1} \Big) } \nonumber \\
& \leq &  (5m)^d \exp \big\{ - c_2 \big( \log m \big)^{\rho} \big\}
\, + \,  \sum_{q=3}^{\infty}{\big[ (2 q + 1)m \big]^d  \exp \big\{ -
c_2 (2 q - 1)^{d-1} m^{d-1}  \big\}  }, \label{rtu}
\end{eqnarray}
for all $m$ and any $m$-box $\Ba$. The first term after the first
inequality in \req{rtu} corresponds to the case where $A=2$, in which
case, $\ga \cap \Ba[2] \not= \emptyset$, and thus, $|\BB_{\LL}
(\bold x)| \geq (\log m)^{\rho(1 - 1/d)}$ for some $\bold x \in
\Ba[2]$
by \req{bkj}. The term indexed by $q$ in the sum after the first inequality
corresponds to the case where $A=q$, with \req{phjk} being used in
place of \req{bkj}. That the constant $c_2$ is positive follows from the
Aizenman-Newman Theorem in \cite[Section 2.4.2]{Hughes},
which proves an exponential rate of decay 
for the probability of a large subcritical cluster containing a given site
in any homogeneous lattice where each vertex has finite degree, 
including $\LL$. From \req{rtu}, 
we obtain (\ref{ubdnew}) for all but finitely many $m$.

Taking $\gas$ to be a greedy lattice animal of minimal size in $\Ba$, 
we find from \req{lbdLn} and \req{ubdnew} that
\begin{equation}\label{tosoh}
\PPP \Big( \big\{ \ga^* \cap \WW \not= \emptyset \big\} \cap \big\{ \,
\textrm{if} \, \, \ga \in
\LA_{\Ba[2]} \, \, \textrm{satisfies} \, \,
  |\ga| \geq (\log m)^\rho , \, \, \textrm{then} \, \,
 \ga \cap \WW \not= \emptyset \big\} \Big) \geq  1 - 2 \eps \, ,
\end{equation}
for all $m$ large enough and any $m$-box $\Ba$.

We apply \cite[Lemma 2.14]{GLAthree}, which is a variant of
a result in
\cite{antpis},
to the supercritical percolation of $\la$-white
sites. Using a union bound, we find that
there exists $\rho = \rho(\la,d)>1$ and $c_3 > 0$
such that, for all $m$ and any $m$-box $\Ba$,
\begin{equation}\label{lemtwod}
\PPP (\exists \bv,\bw \in \Ba[2] : \bv \leftrightarrow
\bw \ \textrm{and} \ D(\bv,\bw) \geq \rho m ) \leq e^{-c_3 m}.
\end{equation}
By \req{tosoh} and \req{lemtwod},
\begin{equation}\label{feqnthree}
\PPP \Big( \forall \ga \in \LA_{\Ba[2]}
 \ \textrm{such that} \ |\ga| \geq (\log m)^\rho : \ \,
D(\ga,\gas) < \rho m \Big) \geq  1 - 3 \eps \,,
\end{equation}
for all $m$ large enough and any $m$-box $\Ba$.

By the definition of $c$, we have that 
$\PPP\big( S(\gas)=G_\Ba \geq c m^d \big) \geq 1-\eps$ for
all $m$ large enough and any $m$-box $\Ba$. By \req{lbdLn},
we have that 
$\PPP\big( \vert \gas \vert = L_\Ba \geq (\log m)^{\rho} + 1 \big) 
\geq 1-\eps$ for such boxes $\Ba_m$.
 So, in view of
the assertion \req{feqnthree}, the proof of the
lemma is complete. $\Box$ \\
\noindent{\bf Proof of Theorem \ref{thmjkl}:}
By Lemma \ref{lemmainr}, we may and shall
set $\la < \infty$, $2<\rho<\infty$ and $c > 0$ such that
$\lim_{m \to \infty} q_{m,\la,c,\rho} = 1$.
\begin{defin}\label{deflact}
For $\ell \in \mathbb{N}$, we say that $\ba \in \Zd$ is {\it $\ell$-active} if the
$\ell$-box $\Ba_{\ell \ba,\ell}$ satisfies
the condition ${\bf A}_{\la,c}^{\rho}$, defined in Lemma \ref{lemmainr}.
\end{defin}
\begin{defin}
A random process $\tau$ taking values in subsets of $\mathbb{Z}^d$ is
said to be a $\rho$-near percolation of parameter $p \in (0,1)$
provided that for any $\bold x \in \mathbb{Z}^d$, $\mathbb{P}(\bold x
\in \tau) =
p$, and for any $\bold x,\bold y \in \mathbb{Z}^d$ for which $\| \bold
x - \bold y \| > \rho$, the events $\big\{ \bold x \in \tau \big\}$ and
$\big\{ \bold y \in \tau \big\}$ are independent.
\end{defin}
\begin{lemma}\label{lemact}
For any $\ell \in \mathbb{N}$, the collection of $\ell$-active sites forms a $(2\rho + 1)$-near percolation.
\end{lemma}
\noindent{\bf Proof:}
Note that, for given $\ba \in \mathbb{Z}^d$, the event 
\begin{equation}\label{evb}
 \Big\{ \Ba_{\ell \ba,\ell} \, \textrm{satisfies condition
 ${\bf A}_{\la,c}^{\rho}$} \Big\}
\end{equation}
is measurable with respect to $\sigma \big\{ X_{\bv}: \ell_{\infty}(\bv,\Ba_{\ell
\ba,\ell}) \leq \rho \ell \big\}$. Indeed, the event that there exists a lattice
animal $\gas \subseteq \Ba_{\ell \ba,\ell}$ such that
$S(\gas) = G_{\Ba_{\ell \ba,\ell}} \geq c \ell^d$ is measurable with respect to
$\sigma \big\{ X_{\bv}:  \bv \in \Ba_{\ell \ba,\ell} \big\}$. The
event $E$ that  $D(\gamma,\gamma^{*}) \leq \rho \ell$ whenever
$\gamma \in  \LA_{\Ba_{\ell \ba,\ell}[2]}$ satisfies $\vert \ga \vert
\geq (\log m)^{\rho}$ occurs if and
only if for each such $\gamma$, there exists a path $\tau_{\gamma,\gamma^{*}}$
of white sites of length at most $\rho \ell$ that has one site
$\bold w \in
\gas$ and another in $\gamma$. Each site $\bold\xi \in \tau_{\ga,\gas}$
satisfies $\ell_{\infty}(\bold\xi, \Ba_{\ell \ba,\ell}) \leq
\ell_{\infty}(\bold\xi,\bold w) \leq \rho \ell$
because $\bold w \in \gas \subseteq \Ba_{\ell \ba,\ell}$, whereas each
site $\bv$ of $\ga$ satisfies $\ell_{\infty}(\bv,\Ba_{\ell \bold
a,\ell}) \leq 2 \ell$. Thus, the event $E$
is measurable with respect to  $\sigma \big\{ X_{\bv}: d(\bv,\Ba_{\ell
\ba,\ell}) \leq \max\{ \rho , 2 \} \ell\big\}$. As $\rho > 2$, this establishes (\ref{evb}).  

For any pair ${\bold a_1},{\bold a_2} \in \mathbb{Z}^d$ that satisfy
$\| {\bold a_1 } - {\bold a_2} \| > 2 \rho + 1$,
\begin{equation}\label{evc}
 \Big\{ \bold v \in \mathbb{Z}^d : \ell_{\infty} \big( \bv, \Ba_{\ell \bold a_1,\ell}
 \big) \leq \rho \ell \Big\} \cap  \Big\{ \bv \in \mathbb{Z}^d : \ell_{\infty} \big( \bv, \Ba_{\ell \bold a_2,\ell}
  \big) \leq \rho \ell \Big\} = \emptyset.
\end{equation}
From (\ref{evb}) and (\ref{evc}), it follows that the events $\big\{
\bold a_1 \, \textrm{is $\ell$-active} \big\}$ and $\big\{
\bold a_2 \, \textrm{is $\ell$-active} \big\}$ are independent. Noting that
$\mathbb{P} \big( 
\bold a \, \, \textrm{is active} \big)$ is independent of $\ba$ completes
the proof of the lemma. $\Box$
\begin{lemma}\label{lemlss}
Let $\sigma > 0$ be given. For any $q \in (0,1)$, there exists $\epsilon
> 0$ such that if $\tau$ is any $\sigma$-near percolation in
$\mathbb{Z}^d$ of 
parameter exceeding $1 - \epsilon$, then there exists an independent 
percolation
$\tau'$ of parameter exceeding $q$ such that $\tau' \subseteq \tau$
almost surely.
\end{lemma}
\noindent{\bf Proof:} The statement of the lemma is implied by
\cite[Theorem $0.0$(i)]{lss}. $\Box$ \\
Recall from \cite[page 23]{Gr} that the percolation probability
$\theta(p)$ of a percolation of parameter $p$ is given by $\PPP_p
(\vert C(\orig) \vert = \infty)$, where $C(\orig)$ denotes the cluster
of open sites containing the origin.
We make use of the fact that
\begin{equation}\label{thetfact}
 \theta(p)
\to 1 \, \textrm{as $p \to 1$},
\end{equation}
which is implied by \cite[Theorem 8.8]{Gr} (there, the result is
being asserted for bond percolation, but the arguments used are valid
for the current case of site percolation).
 By Lemma \ref{lemlss} and \req{thetfact}, 
we may choose $\delta > 0$ such that any $(2
\rho + 1)$-near percolation of parameter exceeding $1 - \delta$
contains a percolation $P$ whose parameter $p$ is such that 
\begin{equation}\label{nthpineq}
\theta(p) > 1 - \eps.
\end{equation}
We now fix $\ell \in \mathbb{N}$ such that
$q_{\ell,\la,c,\rho} \geq 1 - \delta/2$. (We will also be requiring
that $\ell$ is sufficiently high relative to $\la,\rho$ and $c$, but
we prefer to state the particular bounds that are needed as they arise.) By Lemma \ref{lemact}, we may
find such a percolation $P$ satisfying $P \subseteq \big\{ \bold a \in
\mathbb{Z}^d: \bold a \, \,
\textrm{is active} \big\}$, where, now that $\ell$ is fixed, we write
active in place of $\ell$-active for the rest of the proof. For any $n
\in \mathbb{N}$, we write throughout the proof
$n = F \ell + r$ with $F \in \mathbb{N}$ and
$r \in \{ 0, \ldots , \ell- 1 \}$ so that $F$ implicitly depends on $n$.
For each greedy lattice animal $\xi$ in $\Ba_n$, let
$W_{\xi}$ denote the collection of $\ell$-boxes of the form
$\Ba_{\ell \ba,\ell}$ that are contained in $\Ba_n$
(i.e. $\ba \in \{0,\ldots,F-1\}^d$) and which $\xi$ intersects.
On several later occasions, we will use the following definition and lemma.
\begin{defin}\label{bigconndef}
For $P$ a percolation on $\Zd$,
we write $P_{F,C}$ for the largest connected component of $P \cap
\big\{ C, \ldots, F - 1 - C \big\}^d$.
\end{defin}
\begin{lemma}\label{lempva}
For any $j \in \mathbb{N}$, and $P$ a percolation on $\Zd$ of
supercritical parameter $p>p_c$, we have that
\begin{equation}\label{thrptwo}
\liminf_{n \to \infty}{\frac{\vert P_{n,j} \vert}{n^d}} \geq \theta(p)
\, \, \, \textrm{almost surely.}
\end{equation}
\end{lemma}
\noindent{\bf Proof:} 
We prove the lemma in the case where $j=0$, the other cases being no different.
Let $P_{\infty}$ denote the unique infinite component of $P$.
Given $\alpha \in (0,1)$, let the event $Q_n(\alpha)$ be given by 
\begin{eqnarray}
 Q_n(\al) & = & \Big\{ \exists \, \, \textrm{a connected component $C_n = C_n(\al)$ of $P
 \cap \Ba_n$ such that} \nonumber\\
   & & \qquad \qquad \textrm{if $D \subseteq P \cap \Ba_n$ has radius
   exceeding $\alpha n$, then $D \subseteq C_n$.} \Big\} \label{qnevent}
\end{eqnarray}
(Recall that the radius of a connected set $C \subseteq \Zd$ is the maximum
over pairs of sites $\bold x,\bold y \in C$ of the
minimal number of edges in a  path in $C$ from $\bold x$ to $\bold y$.)
By (2.24) of \cite{antpis}, there exists, for each $\al \in (0,1)$, a constant $c(\al) > 0$ such that
\begin{equation}\label{qnprob}
 \PPP \big( Q_n(\al) \big) \leq \exp \{ - cn \},
\end{equation}
for all sufficiently high values of $n$. (Note that the arguments in
\cite{antpis} are performed for bond percolation. The authors note
however that they may be applied to site percolation. Note also that
their argument is applied for the choice $\al = 1/{25}$, but is valid
for each $\al \in (0,1)$.) The Borel-Cantelli lemma
applied to \req{qnprob} implies that $Q_n(\al)$ occurs for all but
finitely many $n$.
We claim that, for each $\al \in (0,1)$,
\begin{equation}\label{polik}
P_{\infty} \cap \big\{ \lfloor 2 \al n \rfloor , \ldots, n - 1 - \lfloor 2 \al
n \rfloor \big\}^d \subseteq C_n(\al),
\end{equation}
for all $n$ sufficiently high. To derive \req{polik}, consider $\bold
x  \in 
P_{\infty} \cap \big\{ \lfloor 2 \al n , \ldots, n - 1 - \lfloor 2 \al
n \rfloor \big\}^d$.  Note
that the connected component of $P \cap \Ba_n$ in which the site
$\bold x$ lies has a radius of at least $\lfloor 2 \al n \rfloor > \al
n$. 
Thus,
if $Q_n(\al)$ occurs, then $\bold x \in C_n$. Hence, we obtain
\req{polik} for high values of $n$. 

We bound
\[
 \vert C_n \vert \geq \big\vert P_{\infty} \cap  \big\{ \lfloor 2 \al n \rfloor, \ldots, n - 1 - \lfloor 2 \al
n \rfloor \big\}^d \big\vert \geq  \big\vert P_{\infty} \cap \Ba_n
\big\vert \, - \, 4 d \al n^d,
\]
so that
\begin{equation}\label{ploki}
\liminf_{n \to \infty}{n^{-d}\vert C_n(\al) \vert} \geq
\liminf_{n \to \infty}{n^{-d}\vert P_{\infty} \cap \Ba_n \vert} \, -
\, 4d\al = \theta(p) - 4d\al,  
\end{equation}
the latter an almost sure equality that is due an application of the
ergodic theorem to the process $P$.
By the definition of the event $C_n(\al)$ appearing in \req{qnevent},
we have that $C_n(\al)=P_{n,0}$, provided that $P_{n,0}$ has radius at
least 
$\al n$. If $\al \in (0,1)$ is chosen to be small enough that
$\theta(p) > 4d \al + 2 \al^d$, then \req{ploki} implies that $\vert
C_n(\al) \vert > \al^d n^d$ for high $n$, whence $\vert P_{n,0} \vert
> \al^d n^d$ for such $n$. Noting that, for any finite connected set $B
\subseteq \Zd$,
\[
 B \subseteq \Big[ \min_{\bold x \in B}{\bold x_1},  \max_{\bold x \in
 B}{\bold x_1} \Big] \, \times \ldots \times \,
 \Big[ \min_{\bold x \in B}{\bold x_d},  \max_{\bold x \in
 B}{\bold x_d} \Big],
\]
and that $\max_{\bold x \in
 B}{\bold x_i} - \min_{\bold x \in
 B}{\bold x_i} \leq {\rm rad}(B)$ for each $i \in \{1,\ldots,d\}$, we find that $\vert B
 \vert \leq {{\rm rad} (B)}^d$.
We deduce that the radius of $P_{n,0}$ exceeds $\al n$, for high $n$,
so that indeed $C_n(\al) = P_{n,0}$ for such $n$.
Given that \req{ploki} holds
for each $\al \in (0,1)$, we deduce that 
\[
\liminf_{n \to \infty}{n^{-d}\vert P_{n,0} \vert} \geq \theta(p),
\]
as we sought.  $\Box$ \\
We define
\begin{eqnarray}
 E_1 & = & \Big\{  \textrm{for some greedy lattice animal $\xi$ in
 $B_n$}, \nonumber \\
 & & \qquad   \textrm{there exist $\bold a_1, \bold a_2 \in P_{F,\lfloor \rho \rfloor +
   1}, \, \, \textrm{such that} \, \, \Ba_{\ell \bold a_1,\ell} \in
   W_{\xi} \, \, \textrm{and} \, \, \Ba_{\ell \bold
   a_2,\ell} \not\in W_{\xi}$} \Big\} \nonumber 
\end{eqnarray}
and
\begin{equation}
 E_2  =  \Big\{  \textrm{for some greedy lattice animal $\xi$ in
 $B_n$,} \, \, \textrm{we have that} \, \, W_{\xi} \cap  \big\{  \Ba_{\ell \bold
   a,\ell} : \bold a \in P_{F, \lfloor \rho \rfloor + 1} \big\} = \emptyset \Big\}.
\end{equation}
We will now prove that the events $E_1$ and $E_2$ occur for finitely
values of $n$ almost surely. In the case of $E_1$, we firstly show
that, for all $n$ sufficiently large,
\begin{equation}\label{claone}
  \big\{ L_n > \big( \log \ell \big)^{\rho} \big\} \cap E_1 = \emptyset.
\end{equation}
To derive (\ref{claone}), suppose that the event on the left-hand-side
occurs. 
The set $P_{F,\lfloor \rho \rfloor + 1}$ being connected, we may
suppose that $\bold a_1$ and $\bold a_2$ are adjacent. Let $\ga_{\bold a_2}$ denote a lattice animal playing the role
of $\gas$ in the condition ${\bf A}_{\la,c}^{\rho}$ that $  \Ba_{\ell \bold
a_2,\ell}$ satisfies. We may locate a connected set $\phi \subseteq
 \xi $ satisfying $\vert \phi \vert = \lfloor (\log \ell)^{\rho}
 \rfloor + 1$ and $\phi \cap  \Ba_{\ell \bold
   a_1,\ell} \not= \emptyset$, because $L_n >  (\log \ell)^{\rho}$. 
 Requiring that $\ell \geq  \lfloor (\log \ell)^{\rho}
 \rfloor + 1$, we see that $\phi \in   \Ba_{\ell \bold
   a_2,\ell}[2]$.  The fact that $  \Ba_{\ell \bold
   a_2,\ell}$ satisfies condition ${\bf A}_{\la,c}^{\rho}$ implies
   that there exists a white path $\tau_{\ga_{\bold a_2},\phi}$ of length at most
   $\rho \ell$ from some site of $\ga_{\bold a_2}$ to some
   site of $\phi \subseteq \xi$. Consider the lattice animal
   $\xi^{+} = \xi \cup \tau_{\ga_{\bold a_2},\phi} \cup \ga_{\bold
a_2}$. Note firstly that $\xi^{+}$ is connected, because $\xi$ and
$\ga_{\bold a_2}$ are, and $\tau_{\ga_{\bold a_2},\phi}$ is a path
between them. Secondly, note that $\xi^{+} \subseteq \Ba_n$: indeed, $\xi \subseteq
\Ba_n$, $\gamma_{\bold a_2} \subseteq \Ba_{\ell \bold a_2,\ell}
\subseteq B_n$, while $\tau_{\gamma_{\bold a_2},\phi}$ is a path of
length at most $\rho \ell$, with  $\tau_{\gamma_{\bold a_2},\phi} \cap
\Ba_{\ell \bold a_2,\ell} \not= \emptyset$. Since $\bold a_2 \in P_{F,
\lfloor \rho \rfloor + 1}$,
$$
 \ell_{\infty} \big(  \Ba_{\ell \bold a_2,\ell} , \Ba_n^c \big) \geq 
 \big( \lfloor \rho \rfloor + 1 \big) \ell  \geq \rho \ell,
$$ 
and thus, $\tau_{\gamma_{\bold a_2},\phi} \subseteq B_n$.

We bound from below the weight of $\xi^{+}$:
\begin{eqnarray}
 S \big( \xi^{+} \big) & = & S \big( \xi \big) + S \big( \gamma_{\bold
 a_2} \big) + S \Big( \tau_{\gamma_{\bold a_2},\phi} \setminus \big( \xi \cup
 \gamma_a \big) \Big) \label{evd} \\
 & \geq &  S \big( \xi \big) + S \big( \gamma_{\bold
 a_2} \big) - \lambda \big\vert  \tau_{\gamma_{\bold a_2},\phi} \big\vert
 \geq   S \big( \xi \big) + S \big( \gamma_{\bold
 a_2} \big) - \lambda \rho \ell, \nonumber
\end{eqnarray}
where, in the equality, we used $\gamma_{\bold a_2} \cap \xi =
\emptyset$, which follows from $\gamma_{\bold a_2} \subseteq  \Ba_{\ell
\bold a_2,\ell} \not\in W_{\xi}$. In the first inequality, the fact that the
path $\tau_{\gamma_{\bold a_2},\phi}$ is white was used. From
$S(\gamma_{\bold a_2}) \geq c \ell^d$ and (\ref{evd}), we deduce that
$S(\xi^{+}) > S(\xi)$, provided that $\ell$ was chosen so that 
$\ell \geq (\la \rho/c)^{1/(d-1)}$. We
have however shown that $\xi^{+}$ is a lattice animal in $\Ba_n$, so
this contradicts the fact that $\xi$ is a greedy lattice animal in
$B_n$. We have proved (\ref{claone}). That $E_1$ may occur for only
finitely many values of $n$ almost surely follows from \req{lmineq}.

It remains to rule out the implausible scenario that a greedy animal
in $\Ba_n$ might fail to meet any $\ell$-box $\Ba_{\ell \bold a,\ell}$
with $\bold a \in P_{F,\lfloor \rho \rfloor + 1}$. 
To prove that $E_2$ may occur for only finitely many values
of $n$ almost surely, we
firstly show, for any $\epsilon > 0$, for $\chi$ satisfying
\begin{equation}\label{nin}
\chi > (1 + \epsilon) d/(d-1),
\end{equation}
and 
for all $n$ sufficiently high, that
\begin{equation}\label{derit}
  \big\{ L_n \geq \big( \log n \big)^{\chi} \big\}
 \, \cap \, \big\{ \bold x \in \Ba_n \implies \big\vert B_{\LL} (\bold
 x)
 \big\vert <  \big( \log n
 \big)^{1 + \epsilon} \big\} \, \cap \, E_2 = \emptyset. 
\end{equation}

To derive (\ref{derit}), suppose now that the event on its
left-hand-side occurs.  We will connect the greedy lattice animal $\xi$
in $B_n$ that exists because event $E_2$ occurs to a weighty lattice
animal $\Psi$ in $\Ba_n$ formed from animals in boxes corresponding to
active sites. To construct $\Psi$, note that for each pair of
adjacent sites $\bold a_1,\bold a_2 \in P_{F,\lfloor \rho \rfloor +
1}$, we may find a white path $\phi_{\bold a_1,\bold a_2}$ from
$\ga_{\bold a_1}$ to $\ga_{\bold a_2}$ of length
at most $\rho \ell$. This is because the condition ${\bf
A}^{\rho}_{\la,c}$ satisfied by $\Ba_{\ell \bold a_2,\ell}$ ensures
that $\vert \gamma_{\bold a_2} \vert \geq (\log \ell)^{\rho} + 1$, so
that the path  $\phi_{\bold a_1,\bold a_2}$ may be found by putting
$\ga = \gamma_{\bold a_2}$ in the  condition ${\bf
A}^{\rho}_{\la,c}$ satisfied by $\Ba_{\ell \bold a_1,\ell}$. We set 
\begin{equation}\label{eqnpsi}
 \Psi = \bigcup_{\ba \in P_{F, \lfloor \rho \rfloor + 1}}{\ga_{\ba}}
 \, \,
 \cup \bigcup_{\bold a_1,\bold a_2 \in P_{F,\lfloor \rho \rfloor + 1:
 \vert\bold a_1 - \bold a_2 \vert  =1 }}{\phi_{\bold a_1,\bold a_2}}.
\end{equation} 
We now check that $\Psi \in \LA_{\Ba_n}$. It is 
connected, because each $\ga_{\ba}$ is,
and to each adjacent pair  $(\bold a_1,\bold a_2)$ of sites in
the connected set $P_{F,\lfloor \rho \rfloor + 1}$, there corresponds a path
$\phi_{\bold a_1,\bold a_2}$ that joins $\ga_{\bold a_1}$ and
$\ga_{\bold a_2}$. 
Note that for $\bold a_1,\bold a_2 \in P_{F,\lfloor \rho \rfloor +
1}$, we have that $\ell_{\infty}(\phi_{\bold a_1,\bold a_2},\Ba_n^c) \geq
\ell_{\infty}(\Ba_{\ell \bold a_1,\ell},\Ba_n^c) - \vert \phi_{\bold
a_1,\bold a_2} \vert \geq (\lfloor \rho \rfloor + 1) \ell - \rho \ell
> 0$, so that $\phi_{\bold a_1,\bold a_2} \subseteq \Ba_n$, and thus,
$\Psi \subseteq \Ba_n$.

Note also that
\begin{equation}\label{evf}
 \big\vert \Psi \big\vert \geq \big\vert P_{F,\lfloor \rho \rfloor +1} \big\vert
 \geq F^d \big( 1 -\epsilon \big) \geq \frac{(n - \ell)^d}{\ell^d}(1-\epsilon),
 \end{equation}
the second inequality following for high choices of $n$ from Lemma \ref{lempva} and \req{nthpineq}.

Let $\xi$ be a greedy lattice animal in $\Ba_n$. We are aiming to find
a path from $\Psi$ to $\xi$ that is white except perhaps for its
endpoints. We may thus assume that $\Psi \cap \xi = \emptyset$, the
other case being trivial. Any set $F$
properly separating $\Psi$ and $\xi$ in $\Ba_n$ satisfies
\begin{eqnarray}
\big\vert F \big\vert & \geq &  \frac{1}{2d} \min \Big\{ \big\vert \xi \big\vert^{1 -
  \frac{1}{d}},  \big\vert \Psi \big\vert^{1 - \frac{1}{d}}
  \Big\} \nonumber \\
 & \geq &   \frac{1}{2d} \min \Big\{ \big( \log n\big)^{\chi(1 - 1/d)},
(1-\epsilon)^{1-1/d}  \frac{(n-\ell)^{d-1}}{\ell^{d-1}} \Big\}
\geq  
\big( \log
n\big)^{1 + \epsilon}, \label{crs}
\end{eqnarray}
for high values of $n$. In the first inequality here, we made use of
the second assertion in Lemma \ref{isomp} (which requires that $\Psi
\cap \xi = \emptyset$), 
while the second is valid
by the occurrence of the first event on the left-hand-side of
\req{derit}, and by \req{evf}. The third is due to \req{nin}. 
By the occurence of the second event in
\req{derit}, and \req{crs}, no black $\LL$-cluster properly 
separates $\xi$ and
$\Psi$ in $\Ba_n$. Applying the first part of Lemma \ref{lemmbul} with
the choices $C = \xi$, $D = \Psi$ and $E = \{ \, \textrm{black sites}
\} \setminus (\xi \cup \Psi)$, we deduce that the black sites do not properly separate
$\xi$ and $\Psi$ in $\Ba_n$, and thus, we may locate 
the desired path $T_{\bold x,\bold y}$ in $\Ba_n$ from $\bold x
\in \xi$
to $\bold y \in \Psi$ that is white with the possible exception of its endpoints.

Consider the case where $\bold x$ and $\bold y$ are white.
Let $\tau_{\bold x,\bold y} = \big( \bold x = \bold x_0,\bold
x_1,\ldots,\bold x_r = \bold y \big)$ denote some
path from $\bold x$ to $\bold y$ in $\Ba_n$ for which $r \leq dn$ and
whose sites may or may not be white. We now modify
$\tau_{\bold x,\bold y}$ to form a white path $\sigma$ from $\bold x$ to
$\bold y$ in $\Ba_n$ such that 
\begin{equation}\label{saq}
 \vert \sigma \vert \leq  d n  (3^d - 1) (\log n)^{1+ \eps}. 
\end{equation}
If $\tau_{\bold x,\bold y}$ is white, we are
done. Otherwise, let $r_1 \in \{ 0,\ldots,r-1 \}$ denote the smallest
value for which $\bold x_{r_1}$ is white and $\bold x_{r_1 + 1}$ is black. Note that
$\bold x_{r_1} \in \partial_{vis(\bold x),n}\big( \BB_{n,\LL}(\bold
x_{r_1 + 1}) \big)$: indeed,
$(\bold x_{r_1},\bold x_{r_1 - 1},\ldots,\bold x_0 = \bold x)$ is disjoint from $\BB_{n,\LL}(\bold x_{r_1 + 1})$. 
We claim that there exists $r_2 \in \{ r_1 + 2,\ldots,r \}$ for
which 
\begin{equation}\label{tepm}
\bold x_{r_2} \in \partial_{vis(\bold x),n}\big( \BB_{n,\LL}(\bold
x_{r_1 + 1}) \big) .
\end{equation}
For, taking $r_2 = 1 + \sup
 \big\{ r' \in \{ r_1 +1,\ldots,r-1 \} : \bold x_{r'} \in  \BB_{n,\LL}(\bold x_{r_1 + 1})
 \big\}$, the path $(\bold x_{r_2},\ldots,\bold x_r)$ is disjoint
 from $ \BB_{n,\LL}(\bold x_{r_1 + 1})$ and may be prefixed to the reversal of
 $T_{\bold x,\bold y}$ (which is white and hence  disjoint
 from $ \BB_{n,\LL}(\bold x_{r_1 + 1})$). The result is a path from
 $\bold x_{r_2}$ to
 $\bold x$ in $\Ba_n$ that does not intersect  $ \BB_{n,\LL}(\bold x_{r_1 + 1})$. Thus,
 \req{tepm}.

By the second assertion of Lemma \ref{lemtop}, we may
find a path $(\bold x_{r_1} = \bold y_0,\bold
y_1,\ldots,\bold y_{r_3} = \bold x_{r_2})
\subseteq \partial_{vis(\bold x),n}  \BB_{n,\LL}(\bold x_{r_1 + 1})$. 
The same path-altering procedure may be applied to 
\[
 \big( \bold x = \bold x_0,\bold x_1,\ldots, \bold x_{r_1} = \bold
 y_0,\bold y_1 , \ldots, \bold y_{r_3} =
 \bold x_{r_2},\ldots, \bold x_r = \bold y \big),
\]
and then iterated. The effect of each iteration is to replace the
passage of the original path through a black $\LL$-cluster by one in
its visible boundary with respect to $\bold x$. After at most $r \leq dn$ applications, the
iteration produces a white path from $\bold x$ to $\bold y$ in
$\Ba_n$. 
The length
of the path increases by at most $\max\{ \vert \partial_{vis(\bold
x),n}(B) \vert - 1 : B \, \,
\textrm{a black $\LL$-cluster in $\Ba_n$} \}$ at each step. Note that $\vert
\partial_{vis(\bold x),n} (B) \vert \leq (3^d - 1) \vert B \vert \leq
(3^d - 1)  (\log
n)^{1 + \epsilon}$, the latter inequality by the occurence of the
second event on the left-hand-side of \req{derit}. Thus, setting
$\sigma$ to be equal to the white path that the iteration
produces, we have obtained \req{saq}.
In the case where at least one of $\bold x$ and $\bold y$ is black, we may
recolour them white for the course of the argument to produce a path
$\sigma$ which is white except for its endpoints $\bold x$ and $\bold y$.

We form the lattice animal $\Phi =
\Psi \cup \sigma \cup \xi$. We showed after \req{eqnpsi} that $\Psi \subseteq \Ba_n$,
from which $\Phi \subseteq
\Ba_n$ immediately follows. We compute
\begin{eqnarray}
 S \big( \Phi \big) & = & S \big( \xi \big) + \sum_{\ba \in
 P_{F,\lfloor \rho \rfloor + 1}}{S \big( \ga_{\ba} \big)} + S \bigg(
 \Big( \sigma \cup  \bigcup_{\bold a_1,\bold a_2 \in P_{F,\lfloor \rho \rfloor + 1:
 \vert\bold a_1 - \bold a_2 \vert  =1 }}{\phi_{\bold a_1,\bold a_2}}
 \Big) \setminus \Big(  \bigcup_{\ba \in P_{F,\lfloor \rho
 \rfloor + 1}}{\ga_{\ba}} \cup \xi  \Big) \bigg) \nonumber \\
 & \geq & S \big( \xi \big) + \big\vert P_{F,\lfloor \rho \rfloor + 1}
 \big\vert c \ell^d  - \la \Big(  \big\vert \sigma \big\vert +
  \rho \ell \Big\vert \Big\{ \{ \bold a_1,\bold a_2\}: \bold a_1,\bold
  a_2 \in P_{F,\lfloor \rho \rfloor + 1}, \vert\bold a_1 - \bold a_2 \vert  =1 \Big\} \Big\vert \Big)
 \nonumber\\
 & \geq & S\big( \xi \big) + (n-\ell)^d \big( 1 - \epsilon\big)c - \la \Big( 
 d n (3^d - 1) (\log n)^{1 + \eps} + d \big\vert  P_{F,\lfloor \rho \rfloor + 1} \big\vert
 \rho \ell \Big) \nonumber \\
 & \geq &  S\big( \xi \big) + n^d \Big[ \big( 1 - \epsilon\big)c -
 \frac{d \la \rho}{\ell^{d-1}} \Big] - C n^{d-1}\ell - \la d n (3^d -
 1) (\log n)^{1 + \eps}, \label{thisineq}
 \end{eqnarray}
where, in the equality, we used the fact that $\ga_{\bold a} \cap \xi =
\emptyset$ for each $\bold a \in P_{F,\lfloor \rho \rfloor + 1}$,
which follows from $W_{\xi} \cap \{ B_{\ell \bold a,\ell}: \bold a \in
P_{F,\lfloor \rho \rfloor + 1 } \} = \emptyset$.
Regarding the first inequality, note that, if either of the endpoints
$\bold x$ and $\bold y$ of $\sigma$ is black, then that endpoint lies in  
$\bigcup_{\ba \in P_{F,\lfloor \rho \rfloor + 1}}{\ga_{\ba}} \, \cup \xi$. 
In the second inequality, we made use of \req{evf} and
\req{saq}. The third follows from the fact that $\vert P_{F,\lfloor
\rho \rfloor + 1} \vert \leq F^d$. Provided that $\ell$ was chosen so that
$\ell > \big( d \la \rho c^{-1}(1 - \eps)^{-1}  \big)^{1/(d-1)}$,
the inequality \req{thisineq} is admissible for at most finitely 
many values of $n
\in \mathbb{N}$, because $S\big( \Phi \big) > S \big( \xi \big)$ would
contradict $\xi$ being a greedy lattice animal in $\Ba_n$. Thus, \req{derit}.

We now check that the two events other than $E_2$ that appear on the
left-hand-side of \req{derit} may occur for only finitely many values
of $n$.
That  $L_n \geq (\log n)^\chi$ for all high $n$ is implied by \req{lmineq}.
Note also that
\[
\PPP \Big(  \exists \, \bold x \in B_n : \big\vert \BB_{\LL} (\bold
x) \big\vert > \big( \log n
 \big)^{1 + \epsilon}  \Big) \leq n^d \PPP \Big( \big\vert
 \BB_{\LL} (\orig) \big\vert > \big( \log n
 \big)^{1 + \epsilon} \Big) \leq n^d \exp \Big\{ - c \big( \log n
 \big)^{1 + \epsilon} \Big\},
\]
with the constant $c = c(\la) > 0$ by \cite[Theorem 2.4.2]{Hughes} and
\req{pbound}. 
Since $\eps > 0$, we see from the Borel-Cantelli lemma that,
for all high $n$, every black $\LL$-cluster intersecting $\Ba_n$ has at
most $(\log n)^{1 + \eps}$ sites.

We conclude that the event $E_2$ may occur for only finitely many
values of $n$ almost surely. Note that  $E_1^c \cap E_2^c$
occurs if and only if for each greedy lattice animal $\xi \in \Ba_n$,
$W_{\xi} \supset \big\{ \Ba_{\ell \ba , \ell} : \ba \in P_{F,\lfloor \rho
\rfloor + 1} \big\}$. This completes the proof of Theorem
\ref{thmjkl}. $\Box$

\section{Existence of $G$; proof of Theorem \ref{pthmthree}}\label{secthree}

In this section, we strengthen a result of \cite{GLAthree}, with hypotheses as
general as those of that paper.

\noindent
{\bf Proof of Theorem \ref{pthmthree}:}
We set $y = \liminf n^{-d} G_n$ and $y + E = \limsup n^{-d} G_n$. Note
that $y,E \in \cap_{n\geq 1} \sigma \{ Y_m : m \geq n \}$, where $Y_m$ is a
vector whose components comprise $\big\{ X_{\bv} : \bv \in \Ba_n - \Ba_{n-1}
\big\}$ in an arbitrary order. The family of random variables $\big\{
Y_n : n \in \mathbb{N} \big\}$ being independent, Kolmogorov's
zero-one law \cite[Theorem 1.8.1]{Durrett} implies that $y$ and $E$ are almost sure
constants. It thus suffices for proving the theorem to show that $E >
\epsilon > 0$ results in a contradiction. 
To this end, we aim to produce a box
percolation whose members contain lattice animals whose weight per
unit volume is close to the value $y$ and which may be connected by
paths of negligible weight. 
If the sidelength of the boxes is chosen to be large, then the
percolation has a high density. As
such, the animals lying in members of the largest connected component
of the box percolation inside  
any very large box may be joined to form a well-spread lattice animal
of weight per unit volume close to $y$. The assumption that $E > \epsilon$
should allow us to identify lattice animals whose weight per unit volume
exceeds $y + \epsilon$, which may replace parts of the constructed
animal, thereby increasing its weight per unit volume. These heavier animals
must be found in a uniform fraction of space and be capable of being
joined to nearby structure at negligible cost if a sufficient increase
in weight is to result from the proposed modification.

By Lemma \ref{lemmainr}, we may and shall set $\la < \infty$, $2 < \rho < \infty$ and $c \in \big(
y - \eps/{5^{6d}},y \big)$ such that $\lim_{m \to
\infty}{q_{m,\la,c,\rho}}=1$. 
By the argument presented
after \req{thetfact} in the proof of Theorem \ref{thmjkl}, 
we may find $\ell \in \mathbb{N}$ so that there exists
a percolation $P$ whose parameter $p$ satisfies 
\begin{equation}\label{thpineq}
\theta(p) > 1 - \eps/(2^{8d}y)
\end{equation}
and for which $P \subseteq \big\{ \bold a \in 
\mathbb{Z}^d: \bold a \, \,
\textrm{is active} \big\}$, where, now that $\ell$ is fixed, we write
active in place of $\ell$-active for the rest of the proof. 
(We also require that $\ell$ be
chosen high relative to $\la,\rho$ and also to $\eps$. We state the
precise bounds as each one arises.)
We will
use $n \in \mathbb{N}$ to denote the large scale in the proof, that of
the patchwork of joined animals, and, similarly 
to the
proof of Theorem \ref{thmjkl}, will write $n = F \ell + r$ with $F
\in \mathbb{N}$ and $r \in \{0,\ldots,\ell - 1\}$. 

By Lemma \ref{lempva} and \req{thpineq}, we find that, for any given
$C \in \mathbb{N}$, and all $n$ (and thus $F$) sufficiently high,
\begin{equation}\label{tkc}
\big\vert P_{F,C} \big\vert \geq \Big( 1 - \frac{\eps}{2^{8d} y} \Big) F^d.
\end{equation}
Each member of the $\ell$-box percolation $\{ \Ba_{\ell \bold a,\ell} :
\bold a \in P \}$ satisfies condition ${\bf A}_{\la,c}^{\rho}$, and
thus contains a lattice animal of
weight exceeding $(y - \epsilon/5^{6d}) \ell^d$ that may be connected to
another such in an adjacent box. We will obtain a backdrop lattice
animal by joining together such animals that lie
in $\ell$-boxes $\Ba_{\ell \bold a,\ell}$ for sites $\bold a$
belonging to 
a large connected component of $P \cap \Ba_F$.
We now make precise the notion of the heavier animal instances of
which we seek to stitch into this patchwork.
The following definition is convenient.
\begin{defin}\label{recdef}
For any $m \in \mathbb{N}$, $m \geq \ell$ and an $m$-box $\Gamma =
B_{\bold x,m}$, we set, for $q \in \mathbb{N}$,
\[
 \Gamma(\ell,q) = \bigcup \big\{ B_{\ell \bv,\ell}[q]: \bv \, \textrm{such
 that $B_{\ell \bv,\ell} \cap \Gamma \not= \emptyset$} \big\}
\]
equal to the union of those $\ell$-boxes whose sup-norm distance from
some $\ell$-box intersecting $\Gamma$ is at most $q$. We also write
\[
 w_{\Gamma} = \big\{ \bold a \in \Zd : \Ba_{\ell \bold a,\ell} \cap
 \Gamma \not= \emptyset \big\}, \, \, \textrm{and} \, \,
 D_{\infty}(w_{\Gamma},q) = \big\{ \bold v \in \Zd : l_{\infty}(\bold v,w_{\Gamma})
 \leq q \big\},
\]
so that $\Gamma(\ell,q) = \bigcup \big\{ \Ba_{\ell \bold v, \ell} : \bold v \in
D_{\infty}(w_{\Gamma},q) \big\}$.
\end{defin}
\begin{defin}\label{pdefnto}
For $m \geq \ell$, an $m$-box  $\Gamma =
B_{\bold x,m}$ is said to be $(c_1,\la,\rho)$-high provided that
\begin{itemize}
\item  
there exists $\gas \in \LA_{\Gamma}$ such that
$S(\gas)=G_{\Gamma} \geq (y + c_1) m^d$ and  $D(\ga,\gas) \leq \rho m$ for
all $\ga \in \LA_{\Gamma[2]}$ such that $|\ga| \geq (\log m)^\rho$,
\item
any two sites $\bold v_1,\bold v_2 \in D_{\infty} \big( w_{\Gamma} , \lfloor
\log(m/\ell) \rfloor \big) \setminus  D_{\infty} \big( w_{\Gamma} , \lfloor
\log(m/\ell) \rfloor - 1 \big)$ that are connected by a path in $P
\cap  D_{\infty} \big( w_{\Gamma} , \lfloor
\log(m/\ell) \rfloor \big)$, are connected by a path in  
 $P
\cap  \big( D_{\infty} \big( w_{\Gamma} , \lfloor
\log(m/\ell) \rfloor \big) \setminus w_{\Gamma} \big)$.  
\end{itemize}
We write `high' for $(\eps,\la,\rho)$-high.  
\end{defin}
We now construct a disjoint collection of high boxes that fill out a
uniform fraction of a large box in $\Zd$.
\begin{lemma}\label{lemkap}
For any $m_1 \geq \ell$, there exists $m_2 \geq m_1$, and a collection
$\kappa_{m_1,m_2}$ of high boxes in $\Zd$ such that if \, $\Gamma =
 \Ba_{\bold x,m} \in \kappa_{m_1,m_2}$, then $m \in \big\{ m_1,\ldots,m_2
\big\}$, if $\Gamma_1,\Gamma_2 \in \kappa_{m_1,m_2}$, then
$\Gamma_1 [1] \cap \Gamma_2 [1] = \emptyset$, while 
\begin{equation}\label{prfour}
\liminf_{m}{
m^{-d} \Big\vert \big( \bigcup_{\Gamma \in
\kappa}{\Gamma} \big) \ \cap B_m \Big\vert} >  \frac{1}{2. 7^d}
\end{equation}
and
\begin{equation}\label{prfive}
\limsup_{m}{
m^{-d} \Big\vert \big( \bigcup_{\Gamma \in
\kappa}{\Gamma} \big) \ \cap B_m \Big\vert} \leq  \frac{1}{3^d}.
\end{equation}
\end{lemma}
\noindent{\bf Proof:}
We claim firstly that
\begin{equation}\label{namit}
\mathbb{P} \big( \Ba_m \ \textrm{is high for infinitely
many $m \in \mathbb{N}$} \big) = 1.
\end{equation}
To derive \req{namit}, note that the box $\Ba_m$ satisfies the first
requirement in the definition of high provided that it satisfies
condition ${\bf A}_{\la,y + \eps}^{\rho}$. The fact that $\limsup_m
m^{-d} G_m = y + E > y + \eps$ implies that there are almost surely
infinitely many values of $m \in \mathbb{N}$ for which there exists
$\xi \in \LA_{\Ba_m}$ such that $S(\xi) = G_m \geq (y + \eps)
m^d$. The proof of Lemma \ref{lemmainr} may be applied to show that 
condition ${\bf A}_{\la,y + \eps}^{\rho}$ holds for all but finitely
many of those $m$ for which such $\xi$ exists.

To handle the second requirement,
suppose that for a pair of sites $\bold v_1,\bold v_2 \in 
D_{\infty} \big( w_{\Ba_m} , \lfloor
\log(m/\ell) \rfloor \big) \setminus  D_{\infty} \big( w_{\Ba_m} , \lfloor
\log(m/\ell) \rfloor - 1 \big)$, we may find a path $\big( \bv_1 = \bx_1,
\bx_2,\ldots,\bx_r = \bv_2 \big)$ with 
\[
\bx_i \in P \cap  
D_{\infty} \big( w_{\Ba_m} , \lfloor \log(m/\ell) \rfloor \big)
\]
 for $i \in
\{1,\ldots,r \}$. If $\bv_1$ and $\bv_2$ are not connected by a path
in  $P \cap \big(  
D_{\infty} \big( w_{\Ba_m} , \lfloor \log(m/\ell) \rfloor \big) \setminus w_{\Ba_m} \big)$,
then the path $(\bx_1,\ldots,\bx_r)$
visits $w_{\Ba_m}$, and we can set $r_1 = \inf \big\{ i \in \{1,\ldots,r\} : \bold x_i
\in w_{\Ba_m} \big\}$ and $r_2 = \sup \big\{ i \in \{1,\ldots,r\} : \bold x_i
\in w_{\Ba_m} \big\}$. The absence of this latter type of path implies that
the sets $C =\{ \bold x_1,\ldots,\bold x_{r_1 - 1}\}$ and  $D = \{ \bold
x_{r_2 + 1},\ldots,\bold x_r \}$ are
separated by $P^c$ in $D_{\infty}(w_{\Ba_m},\lfloor \log(m/\ell)
\rfloor) \setminus w_{\Ba_m}$. Put differently, the sets $C$ and $D$ are separated by
$E = P^c \cup w_{\Ba_m}$ in the box $B(w_{\Ba_m},\lfloor \log(m/\ell)
\rfloor)$. Noting the pairwise disjointness of $C$,$D$ and $E$, we may
apply the first part of Lemma \ref{lemmbul}, to deduce that there is
an $\LL$-connected set $\chi \subseteq P^c \cup w_{\Ba_m}$ that
separates $\{ \bold x_1,\ldots,\bold x_{r_1 - 1} \}$ and  $\{ \bold
x_{r_2 + 1},\ldots,\bold x_r \}$ in
$D_{\infty}(w_{\Ba_m},\lfloor \log(m/\ell) \rfloor)$. Note that 
\begin{equation}\label{mimat}
 \chi \cap \big( D_{\infty}(w_{\Ba_m},i) \setminus  D_{\infty}(w_{\Ba_m},i-1) \big) \not= \emptyset
\end{equation} 
for each $i \in \big\{ 1,\ldots,\lfloor \log(m/\ell) \rfloor \big\}$. 
To see this, note that $l_{\infty}(\bold x_1,w_{\Ba_m}) = \lfloor
\log(m/\ell)\rfloor$ and $\bold x_{r_1} \in w_{\Ba_m}$
imply that for any $i \in \big\{ 1,\ldots,\lfloor \log(m/\ell) \rfloor
\big\}$, there exists $j_1 \in \{1,\ldots,r_1 - 1\}$ for which
$\ell_{\infty}(\bold x_{j_1},w_{\Ba_m}) = i$. Similarly,  there exists
$j_2 \in \{r_2 + 1,\ldots,r \}$ for which
$\ell_{\infty}(\bold x_{j_2},w_{\Ba_m}) = i$. Let $(\bold x_{j_1} =
\bold y_1,\ldots,\bold y_{r_3}
= \bold x_{j_2})$ denote a path from $\bold x_{j_1}$ to $\bold
x_{j_2}$ 
in the connected
set $D_{\infty}(w_{\Ba_m},i) \setminus D_{\infty}(w_{\Ba_m},i-1)$. This path must intersect
$\chi$ by the separation property that $\chi$ satisfies. We have proved
\req{mimat}.

We have seen that if the second requirement for the box $\Ba_m$ to be
high fails, then there exists $\bold x \in D_{\infty}(w_{\Ba_m}, \lfloor
\log(m/\ell) \rfloor)$ such that $\vert C_{P^c,\LL}(\bold x)\vert \geq
\lfloor \log(m/\ell) \rfloor$, where $C_{P^c,\LL}(\bold x)$ denotes the
$\LL$-cluster of $P^c$ that contains $\bold x$. By applying a union bound
and \cite[Theorem 2.4.2]{Hughes},
we see that, for high values of $m$, this latter event has probability at most
\begin{eqnarray}
  & & \Big\vert D_{\infty}(w_{\Ba_m}, \lfloor \log(m/\ell) \rfloor)
  \Big\vert \, \, 
 \PPP \Big( \vert C_{P^c,\LL}(\orig)\vert \geq
\lfloor \log(m/\ell) \rfloor \Big) \label{glj} \\
 & \leq &
 \big( \lfloor m/\ell \rfloor + 2 + 2 \lfloor \log(m/\ell) \rfloor
 \big)^d \exp \Big\{ -c \log(m/\ell)\Big\}, \nonumber
\end{eqnarray}
where $c$ is any constant satisfying $c \in (0,(1-p)^2/(2\chi(1-p)^2))$,
$\chi(t)$ denoting the mean cluster size for a percolation of
parameter $t$ on $(\Zd,\LL)$. Given that $\chi(t) \to 0$ as $t \to 0$,
which follows from \cite[Theorem 6.108]{Gr} in the case of bond
percolation on $\Zd$ (the case of site percolation having the same
proof),
we may suppose that $p$ has been chosen so close to one that 
the sequence in $m$ of terms \req{glj} is
summable. (This lower bound on $p$ is ensured by fixing $\ell$ at a
high enough value, in the same way that \req{thpineq} was obtained.) We deduce from the Borel-Cantelli lemma that the second
condition in the definition of `high' applies to all but finitely many
of the boxes $\Ba_m$ almost surely. We have shown \req{namit}.   

Given $m_1 \in \mathbb{N}$, there exists by \req{namit} $m_2 \geq m_1$
for which
\begin{equation}\label{rlat}
\mathbb{P}\big( \Ba_m \ \textrm{is high for some} \
m \in \{m_1,\ldots, m_2 \} \big) > 1/2. 
\end{equation}
Let $\Delta = \big\{ \bold v \in \Zd: \Ba_{\bold v,m} \, \textrm{is
high for some $m \in \{m_1,\ldots,m_2\}$} \big\}$. Let $\{ \bold z_j:j \in
\mathbb{N} \}$ denote an ordering of $\mathbb{N}^d$ that enumerates
each shell $\Ba_j - \Ba_{j-1}$ in turn. 
Note that the event
$\big\{ \bold x \in \Delta \big\}$ is measurable with respect to $\sigma
 \{ X_{\bold v} : \ell_{\infty}(\bold v,\Ba_{\bold x,m_2}) \leq \rho
 m_2  \}$, provided that $m_1$ (and thus $m_2$) is high enough 
 relative to $\ell$: indeed, the first condition that defines a high box is
 determined by the values of these random variables by \req{evb},
 whereas the occurence of the second is measurable with respect to  
 $\sigma
 \{ X_{\bold v} : \ell_{\infty}(\bold v,\Ba_{\bold x,m_2}) \leq \ell
 \lfloor \log (m_2/\ell) \rfloor + 1\}$. Setting $\hat{m}= \lfloor \rho
 m_2 \rfloor + 1$ and using 
the fact that $\Ba_{\bold x,\hat{m}} \cap
\Ba_{\bold x',\hat{m}} = \emptyset$ if $\| \bold x - \bold x' \| \geq
\hat{m}$, it follows that the sequence of events $\big\{ \bold a + \hat{m}
\bold z_j
\in \Delta: j \in \mathbb{N}\big\}$ is independent, for any given
$\bold a \in \Ba_{\hat{m}}$. This sequence is identically
distributed because the law of $\{ X_{\bv} : \bv \in \mathbb{N}^d \}$ is
translation invariant. As such, the strong law of large numbers
\cite[Theorem 1.7.1]{Durrett}
may be applied to the sequence of random variables $1 \! \! 1 \big\{  \bold a + \hat{m} z_{j}
\in \Delta: j \in \mathbb{N}
\big\}$. By considering the sequence of partial sums corresponding to those $j$ at which
the enumeration of the set $\Ba_j$ is completed by $\bold z_j$, we
deduce from \req{rlat} that
\begin{equation}\label{nand}
 \textrm{the limit} \, \, \lim_{m} \frac{\big\vert
 \Ba_m \cap \Delta \cap \{ \bold a + \hat{m} \mathbb{N}^d
 \}\big\vert}{\big\vert  \Ba_m \cap \{ \bold a + \hat{m}
 \mathbb{N}^d \} \big\vert} \, \, \textrm{exists and exceeds $1/2$}.
\end{equation}
From \req{nand} and the fact that $\big\vert \Ba_m \cap \{ \bold a +
\hat{m} \mathbb{N}^d \} \big\vert = m^d \hat{m}^{-d} + \hat{m}^{1-d} O(m^{d-1})$, it
follows that
\begin{eqnarray}
 \lim_{m} m^{-d} \big\vert \Ba_m \cap \Delta \big\vert & = & \lim_{m}
 m^{-d} \sum_{\bold a \in \Ba_{\hat{m}}}{\big\vert  \Ba_m \cap \Delta  \cap \{ \bold a + \hat{m}
 \mathbb{N}^d \} \big\vert} \nonumber \\
 & = & \lim_{m} \big[ \hat{m}^{-d} + \hat{m}^{1-d}
 O(m^{-1}) \big] \sum_{\bold a \in \Ba_{\hat{m}}}{ \frac{\big\vert
 \Ba_m \cap \Delta \cap \{ \bold a + \hat{m} \mathbb{N}^d
 \}\big\vert}{\big\vert  \Ba_m \cap \{ \bold a + \hat{m}
 \mathbb{N}^d \} \big\vert}} > 1/2. \label{zxe}
\end{eqnarray}
For $\bold x \in \Delta$, let $\Gamma_{\bold x} = \Ba_{\bold
x,m}$ with $m
\in \{ m_1,\ldots,m_2 \}$ maximal such that $\Ba_{\bold x,m}$ is
high. Let $\Delta' = \{ \Gamma_{\bold x} : \bold x \in \Delta \}$. 
It remains to disjointify the collection of boxes $\Delta'$ while
retaining enough of its members so that their union has
positive density. 
To do this, enumerate 
\[
\Delta' = \{ \Ba_{\bold x_{m,j},m} : m
\in \{ m_1,\ldots,m_2 \} , j \in \mathbb{N} \},
\]
so that $\{ \bold x_{m,j}: j \in \mathbb{N} \}$ is an ordering of
those $\bold x
\in \Delta$ for which $\Gamma_{\bold x}$ has sidelength $m$. We will
iteratively examine the indices $(m,j)$ that label members of
$\Delta'$,
admitting one at each
step into a
set of {\it accepted} indices $\bold A$ while at the same time placing
others in a set of {\it rejected} indices $\bold R$. We will allow these
symbols to denote those indices currently accepted or rejected at each
step, without using further labelling. At the start, $\bold A = \bold
R = \emptyset$. We begin by examining the indices $\{ (m_2,j): j \in
\mathbb{N} \}$. At the first step, we put $(m_2,1)$ in $\bold A$, and reject (put
in $\bold R$) those $(m,i)$ (except for $(m_2,1)$) for which
$\Ba_{\bold x_{m,i},m}[1] \cap
\Ba_{\bold x_{m_2,1},m_2}[1] \not= \emptyset$. At the generic step for
boxes of sidelength $m_2$, we put $(m_2,i)$ in $\bold A$, where $i
\in \mathbb{N}$ is minimal for which $(m_2,i)$ is not currently in
$\bold A \cup \bold R$, and put in $\bold R$,
\[
 \big\{ (m,j): m \in \{ m_1,\ldots,m_2 \}, j \in \mathbb{N},
 (m,j)\not=(m_2,i):  \Ba_{\bold x_{m,j},m}[1] \cap  \Ba_{\bold x_{m_2,i},m_2}[1] \not= \emptyset  \big\}.
\]
After at most countable many iterations, $(m_2,i) \in \bold A \cup
\bold R$ for each $i \in \mathbb{N}$. We proceed to deal with those
$\{ (m,i) : i \in \mathbb{N} \}$ not yet in $\bold A \cup
\bold R$, for $m = m_2 - 1$, then for each $m$ in descending order
until we finish with $m = m_1$. At the generic step when $m \in \{
m_1,\ldots,m_2 \}$ is some fixed value, $(m,i)$ is admitted to $\bold
A$ for the least $i$ for which it is not already in $\bold A \cup
\bold R$, while all those other $(m',j)$ for which $\Ba_{\bold x_{m,i},m}[1]$ intersects  $\Ba_{\bold x_{m',j},m'}[1]$ enter $\bold R$. At the end of the
 procedure, each $(m,i)$ lies in $\bold A \cup \bold R$. We set
 $\kappa = \big\{ \Ba_{\bold x_{m,i},m} : (m,i) \in \bold A \big\}$, with
 $\bold A$ now denoting the collection of accepted indices at the end.

The first two properties asserted for the collection $\kappa$ follow
directly by its construction. We claim that 
\begin{equation}\label{zxi}
\bigcup_{\Gamma \in \Delta'}{\Gamma} \subseteq 
\bigcup_{\Gamma \in \kappa}{\Gamma[3]}.
\end{equation}
To show \req{zxi}, note that each index $(m,i)$ of some box in
$\Delta'$ is eventually either
accepted and so lies in $\kappa$ (so that the box certainly lies in the set on the
right-hand-side of \req{zxi}), or is rejected by the
algorithm. If it is rejected, consider the index $(m',j)$ whose
admission to $\bold A$ resulted in $(m,i)$ joining $\bold R$. The key
point is that this can only happen if $m' \geq m$, because all boxes
in $\Delta'$ whose sidelength exceeds $m'$ have been dealt with by the
time $(m',j)$ is admitted to $\bold A$. This fact along with the
criterion for the rejection of $(m,i)$, namely $\Ba_{\bold
x_{m',j},m'}[1] \cap 
\Ba_{\bold x_{m,i},m}[1] \not= \emptyset$ imply that $\Ba_{\bold x_{m',j},m'}[3]
\supseteq \bold \Ba_{\bold x_{m,i},m}$. 
Thus, \req{zxi}. We may now estimate
\begin{eqnarray}
 & & \Big\vert \Big( \bigcup_{\Gamma \in \kappa}{\Gamma} \Big) \cap \Ba_k
 \Big\vert  \geq  \Big\vert  \bigcup \{ \Gamma \in \kappa,\Gamma
 \subseteq \Ba_k \big\} \Big\vert = \sum_{\Gamma \in \kappa,\Gamma
 \subseteq \Ba_k}{\vert \Gamma \vert} = 7^{-d}  \sum_{\Gamma \in
 \kappa,\Gamma \subseteq \Ba_k}{\vert \Gamma[3] \vert} \nonumber \\ 
 & \geq &  7^{-d} \Big\vert  \bigcup \big\{ \Gamma[3]:\Gamma \in \kappa,\Gamma
 \subseteq \Ba_k \big\} \Big\vert \geq 7^{-d} \Big\vert  \Big( \bigcup \big\{
 \Gamma[3]:\Gamma \in \kappa \big\} \Big) \cap \big\{ \bold x \in \Ba_k :
 \ell_{\infty}(x,\Ba_k^c) \geq 4 m_2 \big\} \Big\vert \nonumber \\
 & \geq & 7^{-d} \Big\vert  \Big( \bigcup \big\{
 \Gamma: \Gamma \in \Delta' \big\} \Big) \cap \big\{ \bold x \in \Ba_k :
 \ell_{\infty}(\bold x,\Ba_k^c) \geq 4 m_2 \big\} \Big\vert \nonumber \\
 & \geq & 7^{-d} \Big\vert
 \Delta \cap \big\{ \bold x \in \Ba_k :
 \ell_{\infty}(\bold x,\Ba_k^c) \geq 4 m_2 \big\} \Big\vert \nonumber \\
 & \geq & 7^{-d} \big\vert \Delta \cap \Ba_k \big\vert - 8dm_2 7^{-d}
 k^{d-1} \geq \big( 1/2 + o(1) \big) 7^{-d} k^d  \, \, - \, \, 
 O(k^{d-1}), \nonumber
\end{eqnarray}
which implies \req{prfour}. In the first equality, we used the fact
that $\kappa$ is a disjoint collection of sets.
In the second equality, we used the fact that
$\vert \ga[3] \vert = 7^d \vert \ga \vert$, which is valid for any box
$\ga$. In
the second inequality, we used that if for any $\Ba_{\bold y,m} \in \kappa$,
there exists $\bold x \in \Ba_k \cap \Ba_{\bold y,m}$ such that
$l_{\infty}(\bold x,B_k^c) \geq 4m_2$, 
then $\Ba_{\bold y,m}[3] \subseteq \Ba_k$, this following directly
from $m \leq m_2$. In the third inequality, \req{zxi} was used, and, in
the fourth, that $\Ba_{\bold x,m} \in \Delta' \implies \bold x \in \Delta$. In the
final inequality, \req{zxe} was used. 
The property \req{prfive} is derived by a similar estimate, that makes
use of  the disjointness of the collection $\big\{
\Gamma[1]:\Gamma \in \kappa \big\}$ and the fact that $\vert \Gamma[1]
\vert = 3^d \vert \Gamma \vert$ for any box $\Gamma$. $\Box$\\    
In the application, we insist that the value of $m_1 \in \mathbb{N}$
satisfy the inequality
\begin{equation}\label{monebd}
 \frac{\rho \lambda}{m_1^{d-1}}  < \frac{\eps}{7^{5d}}.
\end{equation}
We will mention conditions stipulating that $m_1$ must be high relative to $\ell$,
$\rho$, $\la$ and $c$ as they arise.
We will be joining animals lying in the boxes of $\kappa$ to a
structure of joined lattice animals that lies in $\ell$-boxes.
To do so, our first step is
to make space in the fabric of joined lattice animals for the high
boxes lying in $\kappa$. We now claim that, for $n$ sufficiently high
almost surely,
\begin{equation}\label{bal}
 P_{F,C_1} \setminus \bigcup_{\Ba_{\bold x,m} \in \kappa}{D_{\infty} \big( w_{\Ba_{\bold x,m}},
 \lfloor \log (m/\ell) \rfloor \big)}
\end{equation}
lies in a connected component of  $P_{F,\lfloor \rho \rfloor + 1} \setminus \cup_{\Gamma \in
\kappa}{w_{\Gamma}}$.
where $C_1 = \lfloor m_2/\ell + 2 \log(m_2/\ell) \rfloor + \lfloor
\rho \rfloor + 2$. To show this, consider $\bold a_1,\bold a_2$ that lie in the
set on the left-hand-side of \req{bal}. Let $\tau = \big( \bold a_1 =
\bold y_1,\ldots,\bold y_r = \bold a_2 \big) \subseteq P \cap \{ C_1,\ldots,F-1-C_1 \}^d$ be a path from $\bold a_1$ to $\bold
a_2$. We aim to modify the path $\tau$ to find a new one from 
$\bold a_1$ to $\bold a_2$ that avoids each of the sets $w_{\Gamma}$
 for $\Gamma \in \kappa$ while
staying in $P \cap \{ \lfloor \rho \rfloor + 1,\ldots, F-2 - \lfloor
\rho \rfloor \}^d$. We may assume then that $\tau$ does not itself
satisfy these requirements, so that 
$\tau \cap w_{\Ba_{\bold x,m}} \not= \emptyset$ 
for some $\Ba_{\bold x,m} \in \kappa$. Set 
\[
 r_1 = \inf \big\{ i \in \{ 1,\ldots,r \} : \bold y_i \in
 D_{\infty}\big(w_{\Ba_{\bold x,m}}, \lfloor \log(m/\ell) \rfloor \big) \big\}
\]
and
\[
 r_2 = \inf \big\{ i \in \{ r_1+1,\ldots,r \} : \bold y_i \not\in
 D_{\infty}\big(w_{
 \Ba_{\bold x,m}}, \lfloor \log(m/\ell) \rfloor \big) \big\} - 1.
\]
Note that $r_1 > 1$,  because $\bold a_1 \not\in D_{\infty}\big(w_{
 \Ba_{\bold x,m}}, \lfloor \log(m/\ell) \rfloor \big)$, and that
 $\ell_{\infty}(\bold y_{r_i},w_{\Ba_{\bold x,m}}) = \lfloor
\log(m/\ell) \rfloor$ for $i \in \{1,2\}$. The subpath $\big(
 \bold y_{r_1},\ldots,\bold y_{r_2}\big)$ is an excursion of $\tau$
 inside the set $D_{\infty}\big(w_{
 \Ba_{\bold x,m}}, \lfloor \log(m/\ell) \rfloor \big)$. The $m$-box
 $\Ba_{\bold x,m} \in
 \kappa$ being high, we may, by the second requirement in the
 definition of `high', find a path $\big(
 \bold y_{r_1}= \bold z_1,\ldots,\bold z_{r_3}=\bold y_{r_2}\big)$ such that 
\begin{equation}\label{cna}
 \bold z_i \in P \cap \Big(   D_{\infty}\big(w_{
 \Ba_{\bold x,m}}, \lfloor \log(m/\ell) \rfloor \big) \setminus w_{\Ba_{\bold x,m}} \Big)
\end{equation} 
for $i \in \{1,\ldots r_3 \}$.
Note that 
\begin{eqnarray}
 \ell_{\infty}(\bold z_i,B_F^c) & \geq & \inf \big\{ \ell_\infty
 (\bold z,B_F^c) :
 \, \textrm{$\bold z$ such that $\ell_\infty (\bold z,w_{\Ba_{\bold x,m}}) = \lfloor
 \log (m/\ell) \rfloor $} \big\} \label{cnc} \\
 & \geq & \ell_\infty (\bold y_{r_1},B_F^c) - \big( \lfloor m/\ell \rfloor + 1
 + 2 \lfloor \log(m/\ell) \rfloor \big) \nonumber \\
 & \geq &  C_1 -  \big( \lfloor m/\ell \rfloor + 1
 + \lfloor \log(m/\ell) \rfloor \big) \geq  \lfloor \rho \rfloor + 1, \nonumber
\end{eqnarray}
the third inequality valid by $\bold y_{r_1} \in \{C_1,\ldots,F-1-C_1
\}^d$ and the fourth due to $m \leq m_2$.
For any $\Ba_{\bold x',m'} \in \kappa$ for which
$(\bold x',m')\not=(\bold x,m)$, we have that
\begin{equation}\label{cnd}
 \bold z_i \in   D_{\infty}\big(w_{\Ba_{\bold x,m}}, \lfloor \log(m/\ell) \rfloor \big) \subseteq
 D_{\infty}\big(w_{\Ba_{\bold x',m'}}, \lfloor \log(m'/\ell) \rfloor \big)^c, 
\end{equation}
where the containment follows from 
$\Ba_{\bold x,m}[1]\cap \Ba_{\bold x',m'}[1] = \emptyset$
and a choice for $m_1$ that satisfies $m_1 \geq C \ell$.
By \req{cna},\req{cnc} and \req{cnd}, we see that the path 
\[
 \big( \bold a_1 =
 \bold y_1,\ldots,\bold y_{r_1}=\bold z_1,\ldots,\bold z_{r_3}=\bold
 y_{r_2},\ldots,\bold y_r = \bold a_2
 \big) \subseteq P \cap \big\{ \lfloor \rho \rfloor + 1,\ldots,F -
 \lfloor \rho \rfloor - 2 \big\}^d 
\]
has removed any instance of a visit to $w_{\Ba_{\bold x,m}}$ during the
excursion in $  D_{\infty}\big(w_{\Ba_{\bold x,m}}, \lfloor \log(m/\ell) \rfloor \big)$ from $\bold y_{r_1}$ to
 $\bold y_{r_2}$
 without introducing any new visits to 
\[
\bigcup_{\Ba_{\bold x',m'} \in \kappa}   D_{\infty}\big(w_{\Ba_{\bold x',m'}}, \lfloor \log(m'/\ell) \rfloor \big).
\]
We modify the path in such a way, for each example of an excursion into
the set $D_{\infty}\big(w_{
 \Ba_{\bold x',m'}}, \lfloor \log(m'/\ell) \rfloor \big)$ for any $\Ba_{\bold x',m'}
 \in \kappa$. After a finite number of such alterations, we obtain a
 path $\phi$ from $\bold a_1$ to $\bold a_2$ in $P \cap  \{ \lfloor
 \rho \rfloor + 1,\ldots,F - \lfloor \rho \rfloor - 2 \}^d$ that is disjoint
 from 
$\bigcup \big\{ w_{\Ba_{\bold x',m'}} : \Ba_{\bold x',m'} \in \kappa \big\}$: indeed, 
 any excursion of $\phi$ in a set  $D_{\infty}\big(w_{\Ba_{\bold
 x,m}}, \lfloor \log(m/\ell) \rfloor \big)$ will not
 visit $w_{\Ba_{\bold x,m}}$, nor $  D_{\infty}\big(w_{
 \Ba_{\bold x',m'}}, \lfloor \log(m'/\ell) \rfloor \big) \supseteq w_{\Ba_{\bold x',m'}}$ by construction. We have proved that the set in
 \req{bal} lies in a connected component of $P_{F,\lfloor \rho \rfloor
 + 1} \setminus \cup_{\Gamma \in \kappa}{w_\Gamma}$, as we sought to do. We
 call this connected component the backdrop and denote it by $BD = BD(n,\ell)$.

We will require the following lower bound on $\vert BD \vert$:
\begin{equation}\label{tbe}
\big\vert BD \big\vert 
\geq 
\big\vert P_{F,C_1} \big\vert 
- \Big( 1 + \frac{\eps}{10^{10d}y} \Big) \ell^{-d}
\Big\vert \Big( \bigcup_{\Gamma \in
\kappa}\Gamma \Big) \cap \Ba_n \Big\vert.
\end{equation}
To obtain this, note that
\begin{eqnarray}
 \vert BD \vert & \geq & \vert P_{F,C_1} \vert - \Big\vert \Big(
 \bigcup \Big\{ D_{\infty}\big( w_{\Ba_{\bold x,m}} , \lfloor \log (m/\ell)
 \rfloor \big) : \Ba_{\bold x,m} \in \kappa \Big\} \Big)
 \cap P_{F,C_1}\Big\vert \label{cva} \\
  & \geq & \vert P_{F,C_1} \vert - \Big\vert 
 \bigcup \Big\{ D_{\infty}\big( w_{\Ba_{\bold x,m}} , \lfloor \log (m/\ell) \rfloor
 \big) : \Ba_{\bold x,m} \in \kappa, \Ba_{\bold x,m} \subseteq \Ba_n \Big\} 
 \Big\vert, \nonumber
\end{eqnarray}
the first inequality being due to the definition of $BD$. The second is valid
because $m \leq m_2$ implies that $C_1 \geq \lfloor m/\ell \rfloor + 2 + 2 \lfloor \log(m/\ell)
\rfloor$, from which it follows that if $D_{\infty}\big( w_{\Ba_{\bold x,m}} , \lfloor \log (m/\ell) \rfloor
 \big) \cap \big\{C_1,\ldots,F-1-C_1 \big\}^d \not= \emptyset$, then
 $D_{\infty}\big( w_{\Ba_{\bold x,m}} , \lfloor \log (m/\ell) \rfloor
 \big) \subseteq B_F$ and thus $\Ba_{\bold x,m} \subseteq \Ba_n$ (for,
 $\Ba_{\bold x,m}
 \not\subseteq \Ba_n \implies w_{\Ba_{\bold x,m}} \cap B_F^c \not=
 \emptyset$). Note that
\begin{equation}\label{cvb}
 \big\vert w_{\Ba_{\bold x,m}} \big\vert \geq \lfloor m/\ell \rfloor^d 
\end{equation} 
and
\begin{equation}\label{cvc}
 \big\vert w_{\Ba_{\bold x,m}} \big\vert \leq \big( \lfloor m/\ell \rfloor
 + 2 \big)^d \leq \big( 1 + \frac{\eps}{3 \cdot 10^{10d}y} \big) \frac{m^d}{\ell^d} = \big( 1 +
  \frac{\eps}{3 \cdot 10^{10d}y} \big) \frac{\vert \Ba_{\bold x,m} \vert}{\ell^d}, 
\end{equation} 
given that $m_1 \geq C(\eps,y) \ell$. We find that
\begin{eqnarray}
 \Big\vert D_{\infty}\big( w_{\Ba_{\bold x,m}} , \lfloor \log (m/\ell) \rfloor
 \big) \Big\vert & \leq & \Big( \lfloor m/\ell \rfloor + 2 + 2 \lfloor \log(m/\ell)
\rfloor \Big)^d \nonumber \\
  & \leq & \Big( 1 + \frac{ 2(1 +
  \log(m/\ell) )}{\lfloor m/\ell \rfloor} \Big)^d  \big\vert w_{\Ba_{\bold
  x,m}} \big\vert \nonumber \\
 & \leq & \Big( 1 +  \frac{\eps}{3\cdot 10^{10d}y} \Big)  \big\vert w_{\Ba_{\bold x,m}} \big\vert \nonumber \\
 & \leq & \Big(
 1 +  \frac{\eps}{3\cdot 10^{10d}y} \Big)^2 \frac{\vert \Ba_{\bold x,m}
 \vert}{\ell^d} \leq  \Big(
 1 +  \frac{\eps}{10^{10d}y} \Big) \frac{\vert \Ba_{\bold x,m}
 \vert}{\ell^d} , \nonumber 
\end{eqnarray}
the second inequality by \req{cvb}, the third by  $m_1 \geq C(\eps,y)
\ell$, the fourth due to \req{cvc} and
the fifth from $\eps < 3\cdot 10^{10d}y$. Thus,
\begin{eqnarray}
 & & \sum_{\Ba_{\bold x,m} \in \kappa,\Ba_{\bold x,m} \subseteq \Ba_n}\big\vert D_{\infty}\big( w_{\Ba_{\bold x,m}} , \lfloor \log (m/\ell) \rfloor
 \big) \big\vert  \nonumber \\
 & \leq &  \Big( 1 + \frac{\eps}{10^{10d}y} \Big) \ell^{-d} \sum_{\Gamma
 \in \kappa,\Gamma \subseteq \Ba_n}{\big\vert \Gamma \big\vert} \leq 
 \Big( 1 + \frac{\eps}{10^{10d}y} \Big)  \ell^{-d}  \Big\vert \Big( \bigcup_{\Gamma \in
 \kappa} \Gamma \Big) \cap \Ba_n \Big\vert. \label{cvd}
\end{eqnarray}
By \req{cva},
 \req{cvd} and the triangle inequality, follows \req{tbe}.

We now define the lattice animal $\hat{\Psi}$ formed from animals in
$\ell$-boxes corresponding to active
sites of the backdrop $BD$ and into which we will stitch animals from
high boxes in $\kappa$: let
\begin{equation}\label{gbnm}
 \hat{\Psi} = \bigcup_{\ba \in BD}{\ga_{\ba}}
 \, \cup \bigcup_{\bold a_1,\bold a_2 \in BD:
 \vert\bold a_1 - \bold a_2 \vert  =1 }{\phi_{\bold a_1,\bold a_2}}.
\end{equation}
(Recall from after \req{claone} that $\ga_{\ba}$ is the animal $\ga^{*}$ in the condition
${\bf A}_{\la,c}^{\rho}$ satisfied by $\Ba_{\ell \bold a,\ell}$, and,
from before \req{eqnpsi},
that each $\la$-white path $\phi_{\bold a_1,\bold a_2}$ 
satisfies $\vert \phi_{\bold a_1,\bold a_2} \vert \leq \rho \ell$, and
intersects each of $\ga_{\bold a_1}$ and $\ga_{\bold a_2}$.)

There may be a few high boxes of $\kappa$ and contained in  $\Ba_n$ whose greedy lattice animal
(or animals) cannot be connected to $\hat{\Psi}$ in the intended way, if
it so happens that $\hat{\Psi}$ does not reach into a neighbourhood
of these high boxes that would allow the greedy animals therein to
attach to $\hat{\Psi}$. In addition, when a path can be formed from
inside the high box to $\hat{\Psi}$, we must ensure that the path stays in
$\Ba_n$, which amounts to insisting that the box be at a certain
distance from the complement of $\Ba_n$. We now define the set of high
boxes in $\kappa$ whose greedy lattice animal we will connect to
$\hat{\Psi}$, bearing in mind these two requirements.
\begin{defin} 
Let the set of {\it useful high} boxes $UH$ be given by
\begin{eqnarray}
UH & = & \Big\{ \Ba_{\bold x,m} \in \kappa : \Ba_{\bold x,m} \subseteq \big\{
\lfloor \rho m_2 \rfloor + 1,\ldots, n - \lfloor \rho m_2 \rfloor - 2 \big\}^d, \nonumber
\\
  & & \qquad BD \cap \big( D_{\infty}(w_{\Ba_{\bold x,m}}, \lfloor m/(2\ell)
  \rfloor ) \setminus  D_{\infty}(w_{\Ba_{\bold x,m}}, \lfloor \log(m/\ell) \rfloor) \big) \not=
  \emptyset \Big\}.
\end{eqnarray}   
\end{defin}
We now show that it is only a few boxes in $\kappa$ contained in
$\Ba_n$ that do not make it into $UH$. Specifically, we prove that,
for all sufficiently high $n$,
\begin{equation}\label{reone}
 \Big\vert  \bigcup_{\Gamma \in UH}{\Gamma} 
 \Big\vert \geq  \Big\vert  \Big( \bigcup_{\Gamma \in \kappa}{\Gamma}
 \Big) \cap \Ba_n \Big\vert - \frac{\epsilon n^d}{2^{8d-1}y}. 
\end{equation}
To this end, note that
\begin{eqnarray}
  \Big( \bigcup_{\Gamma \in \kappa}{\Gamma}
 \Big) \cap \Ba_n & \subseteq & \Big( \bigcup_{\Gamma \in UH}{\Gamma}
   \Big) \, \cup \, \bigcup \Big\{ \Gamma \in \kappa \setminus UH :
 \Gamma \subseteq \big\{ \lfloor \rho m_2 \rfloor + 1,\ldots, n -
 \lfloor \rho m_2 \rfloor - 2 \big\}^d \Big\} \nonumber \\
 & & \qquad  \cup \, \Big\{ \bold x \in \Ba_n :
 \ell_\infty (\bold x,\Ba_n^c) \leq \lfloor \rho m_2 \rfloor + m_2
 \Big\}. \label{era}
\end{eqnarray}
To show that the second set on this right-hand-side is small, we adopt
a temporary notation, saying that the box $\Ba_{\bold x,m}$ is `far from
$BD$' if $\Ba_{\bold x,m} \in \kappa \setminus UH$ and $\Ba_{\bold x,m} 
\subseteq \big\{
\lfloor \rho m_2 \rfloor + 1,\ldots, n - \lfloor \rho m_2 \rfloor -
2 \big\}^d$. For such a box,
\begin{equation}\label{erb}
 BD \cap \Big(  D_{\infty}(w_{\Ba_{\bold x,m}}, \lfloor m/(2 \ell) \rfloor) \setminus
 D_{\infty}(w_{\Ba_{\bold x,m}}, \lfloor  \log(m/\ell) \rfloor) \Big) = \emptyset.
\end{equation}
Note that, for $\Ba_{\bold x',m'} \in \kappa$, 
$(\bold x',m')\not=(\bold x,m)$,
\begin{eqnarray}
 & &  D_{\infty}\big(w_{\Ba_{\bold x,m}}, \lfloor m/(2\ell) \rfloor \big) \subseteq 
\big\{ \bold a \in \Zd: \Ba_{\ell \bold a,\ell} \subseteq 
\Ba_{\bold x,m}[1] \big\} \nonumber \\
 & \subseteq &  
\big\{ \bold a \in \Zd: \Ba_{\ell \bold a,\ell} \subseteq 
\Ba_{\bold x',m'}[1] \big\}^c \subseteq   D_{\infty}\big(w_{\Ba_{\bold
x',m'}}, \lfloor m'/(2\ell) \rfloor \big)^c, \label{erc}
\end{eqnarray}
the first and third containments requiring that $m_1 \geq 4 \ell$,
since this ensures that any $m \geq m_1$ satisfies $\ell \lfloor
m/(2\ell) \rfloor +1 \leq m/\ell - 1$.
The second containment follows from $\Ba_{\bold x,m}[1] \cap
\Ba_{\bold x',m'}[1] = \emptyset$.
We have that 
\begin{equation}\label{erd}
 D_{\infty} \big( w_{\Ba_{\bold x,m}} , \lfloor m/(2\ell) \rfloor \big) \subseteq B_F.
\end{equation}
Indeed,
\[
 \rho m_1 \leq \rho m_2 \leq \lfloor \rho m_2 \rfloor +1 \leq
 \ell_{\infty}(\Ba_{\bold x,m},\Ba_n^c) \leq \ell_{\infty}( w_{\Ba_{\bold x,m}}, B_F^c) \ell + 2 \ell,
\] 
so that \req{erd} follows, given that $m_1$ may be chosen so that $m_1
> 2 \ell/\rho$.
We now show that  the
collection of sets
\begin{equation}\label{ere}
 \Big\{ P_{F,C_1} , \Big(  D_{\infty}(w_{\Ba_{\bold x,m}}, \lfloor m/(2\ell) \rfloor 
  ) \setminus  D_{\infty}(w_{\Ba_{\bold x,m}}, \lfloor \log (m/\ell) \rfloor) \Big) : \,
  \textrm{$\Ba_{\bold x,m}$ far from $BD$} \Big\}
\end{equation}
are disjoint, with union contained in $B_F$. 
By \req{erb},\req{erc},\req{erd} and $BD \subseteq \Ba_F$, we find
that \req{ere} is true with $P_{F,C_1}$ replaced by $BD$. However, if
$\bold y \in P_{F,C_1} \cap \big( D_{\infty}(w_{\Ba_{\bold x,m}}, \lfloor m/(2\ell) \rfloor 
  ) \setminus  D_{\infty}(w_{\Ba_{\bold x,m}}, \lfloor \log (m/\ell) \rfloor)
  \big)$ for some $\Ba_{\bold x,m} \in \kappa$, then, using \req{erc}
  along with 
\[
  D_{\infty}(w_{\Ba_{\bold x',m'}}, \lfloor \log (m'/\ell) \rfloor)
\subseteq  D_{\infty}(w_{\Ba_{\bold x',m'}}, \lfloor m'/(2\ell) \rfloor), 
\]
it follows that $\bold y$ belongs to the set in \req{bal}, so that
  $\bold y \in BD$. Thus, \req{ere}.
We estimate
\begin{eqnarray}
 & & \Big\vert \bigcup \big\{ \Gamma \,  \textrm{far from $BD$}
 \big\} \Big\vert = \sum_{ \Gamma \,  \textrm{far from
 $BD$}}{\big\vert \Gamma \big\vert} \leq \ell^d  \sum_{ \Gamma \,
 \textrm{far from $BD$}}{\big\vert w_{\Gamma} \big\vert} \nonumber \\
 & \leq &   \ell^d  \sum_{ \Ba_{\bold x,m} \,
 \textrm{far from $BD$}}{  \big\vert D_{\infty}(w_{\Ba_{\bold x,m}}, \lfloor
 m/(2\ell) \rfloor) \setminus  D_{\infty}(w_{\Ba_{\bold x,m}}, \lfloor \log(m/\ell) \rfloor) \big\vert} \nonumber \\
 & \leq & \ell^d
 \big( F^d - \vert P_{F,C_1} \vert \big) \leq \frac{\eps \ell^d
 F^d}{2^{8d}y} \leq  \frac{\eps n^d}{2^{8d}y}, \label{erf}
\end{eqnarray}
where the second inequality follows from the fact that for any $\Ba_{\bold x,m} \in \kappa$, 
\[
 \big\vert D_{\infty}(w_{\Ba_{\bold x,m}, \lfloor
 m/(2\ell) \rfloor}) -  D_{\infty}(w_{\Ba_{\bold x,m}, \lfloor \log(m/\ell)
 \rfloor}) \big\vert \geq \big(2m/\ell - 3 \big)^d - \big( m/\ell +
 2\log(m/\ell) + 2\big)^d,
\]
allied with \req{cvc}, the inequality $\eps < 3 \cdot 10^{10d}y (2^d -
 2)$ and the lower bound $m \geq m_1 \geq C \ell$.
The third inequality in \req{erf} follows from the claimed property of the
collection in \req{ere}, and the fourth from \req{tkc}. We may now
bound
\begin{eqnarray}
 \Big\vert \bigcup_{\Gamma \in UH}{\Gamma}
  \Big\vert & \geq &
 \Big\vert
 \Big( \bigcup_{\Gamma \in \kappa}{\Gamma}
 \Big) \cap \Ba_n \Big\vert -
 \frac{\eps n^d}{2^{8d}y} -
\Big\vert
 \Big\{ \bold x \in \Ba_n :
 \ell_\infty (\bold x,\Ba_n^c) \leq \lfloor \rho m_2 \rfloor + m_2 \Big\}
\Big\vert \nonumber \\
 & \geq  & \Big\vert 
\Big( \bigcup_{\Gamma \in \kappa}{\Gamma}
 \Big) \cap \Ba_n  \big\vert - \frac{\eps n^d}{2^{8d}y}  - 2d \big(
 \lfloor \rho m_2 \rfloor + m_2 \big) n^{d-1} \nonumber \\
 & \geq &
\Big\vert \Big( \bigcup_{\Gamma \in \kappa}{\Gamma}
 \Big) \cap \Ba_n \Big\vert  -  \frac{\eps n^d}{2^{8d-1}y}, \nonumber
\end{eqnarray}
the first inequality valid by \req{era} and \req{erf}, the final one
valid for all high $n$. We have obtained \req{reone}.

We need to connect greedy lattice animals in the boxes of $UH$ to the
backdrop animal $\hat{\Psi}$. For $\Gamma \in UH$, consider then $\gamma^{*} =
\ga^*_{\Gamma} \in
\LA_{\Gamma}$,
provided by the first of the two conditions that the high box $\Gamma$
satisfies. 
 As $\Gamma
\in UH$, there exists $\ba \in BD$ for which $\ell_\infty (\bold a, w_{\Gamma}) \leq \lfloor
m/(2\ell)\rfloor$ if $\Gamma$ is of the form $\Ba_{\bold x,m}$. Any $\bold y
\in \gamma_{\ba} \subseteq \Ba_{\ell \bold a,\ell}$ satisfies 
\begin{equation}\label{crso}
\ell_{\infty}(\bold y,\Gamma) \leq \ell \big(
\lfloor m/(2 \ell) \rfloor  + 1 \big) \leq m/2 + 2 \ell \leq 3m/4, 
\end{equation}
given that
$m \geq m_1 \geq 8 \ell$. 
Note that, by the disjointness of $\{ \gamma_{\ba}: \ba \in BD \}$,
\req{tbe}, \req{gbnm}, \req{tkc} and \req{prfive}, 
\begin{equation}\label{psd}
 \big\vert \hat{\Psi} \big\vert \geq \vert BD \vert \geq \Big( 1 -
 \frac{\eps}{2^{8d}y}\Big) F^d - 2\cdot 3^{-d} \Big( 1 + \frac{\eps}{10^{10d}y}\Big)
 (n/\ell)^d \geq \frac{n^d}{10 \ell^d}, 
\end{equation}
for high values of $n$, given that $\eps$ may be chosen so that 
$\eps < y \big(2^{-8d} + 2
\cdot 3^{-d} \cdot 10^{-10d} \big)^{-1} \big( 9/{10} - 2\cdot 3^{-d} \big)$.
Given that $m$ is at most the fixed constant $m_2$, we may by
\req{psd}
find a
lattice animal $\chi \subseteq \hat{\Psi}$ with $\vert \chi \vert =
\lfloor (\log m)^{\rho} \rfloor + 1$ and $\bold y \in \chi$. If $\bold
z \in \chi$, then
$\ell_\infty (\bold z,\Gamma) \leq \lfloor (\log m)^{\rho} \rfloor +
\ell_\infty (\bold y,\Gamma) \leq m$, by \req{crso} and the fact that
we may choose $m_1$ high enough that for each $m \geq m_1$, $(\log m)^{\rho}
\leq m/4$. Thus $\chi \in \Gamma[1]$. By the first condition that the
high box $\Gamma$ satisfies, we may locate a $\lambda$-white path
$\hat{\phi}_{\Gamma}$ from a site of $\ga^{*}$ to one of $\chi$, with
$\vert \hat{\phi}_{\Gamma} \vert \leq \rho m$. We can now define the
lattice animal, modified from $\hat{\Psi}$ in the way
that we sought:
\begin{equation}\label{ghnl}
\Phi = \hat{\Psi} \cup \bigcup_{\Gamma \in UH}{ \big( \ga^{*}_{\Gamma} \cup
\hat{\phi}_{\Gamma} \big)}
\end{equation}  
It remains to verify that $\Phi$ has the required properties. 
It is indeed a lattice animal, for each $\ga^*_{\Gamma}$ is connected to the
animal $\hat{\Psi}$ by a path $\hat{\phi}_{\Gamma}$. We claim that
$\Phi \subseteq \Ba_n$. We may show that $\hat{\Psi}
\subseteq \Ba_n$, in the same way that we showed that $\Psi
\subseteq \Ba_n$ after \req{eqnpsi}. Note also that 
$\ga^{*}_{\Gamma} \subseteq \Gamma \subseteq
\Ba_n$, because $\Gamma \in UH$. If $\bold y \in \hat{\phi}_{\Ba_{\bold x,m}}$, then
\begin{equation}\label{eqnonet}
\ell_\infty (\bold y,\Ba_{\bold x,m}) \leq \ell_\infty (\bold
y,\ga^*_{\Ba_{\bold x,m}})
\leq \vert \hat{\phi}_{\Ba_{\bold x,m}} \vert \leq \rho m \leq \rho m_2,
\end{equation}
the
second inequality due to $\ga^*_{\Ba_{\bold x,m}} \cap \hat{\phi}_{
\Ba_{\bold x,m}} \not= \emptyset$. However, 
\begin{equation}\label{eqntwot}
  \Ba_{\bold x,m}  \subseteq UH \implies
\ell_\infty(\Ba_n^c, \Ba_{\bold x,m}) \geq \lfloor \rho m_2 \rfloor + 1.
\end{equation}  
From \req{eqnonet} and \req{eqntwot}, we deduce that $\bold y \in \Ba_n$. We
have shown that $\Phi \subseteq \Ba_n$.
Note that
\begin{eqnarray}
 S(\Phi) & = & \sum_{\ba \in BD}{S(\ga_{\ba})} + \sum_{\Gamma \in
 UH}{S(\ga^*_{\Gamma})} \nonumber \\
 & & \quad + \, \, S \bigg( \Big( \bigcup_{\Gamma \in
 UH}{\hat{\phi}_{\Gamma}} \, \cup \, \bigcup_{\bold a_1,\bold a_2
 \in BD : \vert \bold a_1 - \bold a_2 \vert  = 1}{\phi_{\bold
 a_1,\bold a_2}} \Big) \setminus \Big( \bigcup_{a \in BD}{\ga_{\ba}} \, \cup \,
 \bigcup_{\Gamma \in UH}{\ga^*_{\Gamma}} \Big) \bigg), \label{vivb}
\end{eqnarray}
since, for $\bold a \in BD$, $\ga_{\ba} \subseteq B_{\ell \ba,\ell}$
is disjoint from any $\ell$-box intersecting any $\Gamma \in \kappa$,
and thus from each $\chi^*_{\Gamma} \subseteq \Gamma$, by the
definition of $BD$.
We bound 
\begin{equation}\label{cto}
 \sum_{\ba \in BD}{S \big( \ga_{\ba} \big)} \geq \vert BD \vert c
 \ell^d \geq \Big( 1 - \frac{\eps}{2^{8d}y} \Big) \Big( y -
 \frac{\eps}{5^{6d}}\Big) F^d \ell^d  \, - \, \Big( 1 +
 \frac{\eps}{10^{10d}y}\Big) y \Big\vert \Big( \bigcup_{\Gamma \in \kappa}
 \Gamma \Big) \cap \Ba_n  \Big\vert, 
\end{equation}
where $S(\ga_{\ba}) \geq c \ell^d$ was used in the first inequality,
the second due to \req{tbe}, \req{tkc} and $y \geq c
\geq y - \eps/(5^{6d})$. Note also that
\begin{eqnarray}\label{ctt}
 \sum_{\Gamma \in UH}{S \big( \ga^*_{\Gamma} \big)} \geq ( y +
 \eps ) \Big\vert \bigcup_{\Gamma \in UH}{\Gamma} \Big\vert \geq
 ( y + \eps ) \Big\vert \Big( \bigcup_{\Gamma \in
 \kappa}{\Gamma} \Big) \cap \Ba_n \Big\vert - 
( y + \eps ) \frac{\eps n^d}{2^{8d - 1}y},
\end{eqnarray}
the first inequality following from the disjointness of the boxes
$\Gamma \in UH$ and the definition of the animals $ \ga^*_{\Gamma}$,
the second due to \req{reone}.
Note that
\begin{equation}\label{dft}
 \sum_{\Ba_{\bold x,m} \in UH}{m} \leq
 \frac{1}{m_1^{d-1}} \sum_{\Gamma \in UH}{\vert \Gamma \vert}
 \leq \frac{n^d}{m_1^{d-1}},
\end{equation}
because $m \geq m_1$ for each $\Ba_{\bold x,m} \in UH$, and the
collection $UH$ is disjoint, with its union contained in $\Ba_n$. We find that
\begin{eqnarray}
 & &  S \bigg(  \Big( \bigcup_{\Gamma \in
 UH}{\hat{\phi}_{\Gamma}} \, \cup \, \bigcup_{\bold a_1,\bold a_2
 \in BD : \vert \bold a_1 - \bold a_2 \vert  = 1}{\phi_{\bold
 a_1,\bold a_2}} \Big) \setminus \Big( \bigcup_{a \in BD}{\ga_{\ba}} \, \cup \,
 \bigcup_{\Gamma \in UH}{\ga^*_{\Gamma}} \Big) \bigg) \nonumber\\
 & \geq & - \la \Big( \sum_{\Gamma \in UH}{\vert \hat{\phi}_{\Gamma}\vert}
 +  \sum_{\bold a_1,\bold a_2
 \in BD : \vert \bold a_1 - \bold a_2 \vert  = 1}{\vert \phi_{\bold
 a_1,\bold a_2} \vert} \Big) \nonumber \\
 & \geq & - \la \rho \Big( \sum_{\Ba_{\bold x,m} \in UH}{m} \, + \, \ell
 \big\vert \big\{ \{ \bold a_1,\bold a_2 \}: \bold a_1,\bold a_2 \in
 BD , \, \vert \bold a_1 -
 \bold a_2 \vert = 1  \big\} \big\vert \, \Big) \nonumber \\
 & \geq & - \la \rho \big( n^d/m_1^{d-1} + d \ell F^d \big) \geq -
 \frac{\eps n^d}{7^{5d}}  -  d\la \rho \ell F^d, \label{ctv}
\end{eqnarray}
the final inequality by \req{monebd}.
Substituting the bounds \req{cto}, \req{ctt} and
\req{ctv} into \req{vivb} yields
\begin{eqnarray}
 S(\Phi) & \geq & \bigg[ y + \eps - \Big(  1 + \frac{\eps}{10^{10d}y}
 \Big) y   \bigg] \Big\vert \Big( \bigcup_{\Gamma \in
 \kappa}{\Gamma}  \Big) \cap \Ba_n \Big\vert \label{vvc} \\
 & & + \,  \Big(  1 - \frac{\eps}{2^{8d}y}
 \Big) \Big( y - \frac{\eps}{5^{6d}} \Big) F^d \ell^d  -
 \frac{\eps}{2^{8d-1} y}  \big( y + \eps \big) n^d
 -\frac{\eps}{7^{5d}} n^d -   d\la\rho \ell F^d
 \nonumber \\
 & \geq & \bigg[ y + \eps - \Big(  1 + \frac{\eps}{10^{10d}y}
 \Big) y   \bigg] \frac{n^d}{2\cdot 7^d} \nonumber \\
 & & + \,  \Big(  1 - \frac{\eps}{2^{8d}y}
 \Big) \Big( y - \frac{\eps}{5^{6d}} \Big)\big( n - \ell \big)^d  -
 \frac{\eps}{2^{8d-1} y}  \big( y + \eps \big) n^d
 -\frac{\eps}{7^{5d}} n^d -   \frac{\eps}{3^{10d}}n^d,
 \nonumber
\end{eqnarray}
the second inequality using \req{prfour} and the inequality $\ell >
\big( 3^{10d} d \la \rho \eps^{-1} \big)^{1/(d-1)}$ that we may
require that $\ell$ satisfies.
For large values of $n \in
\mathbb{N}$, the dominant term in the last expression is the one in
$n^d$, whose coefficient is bounded below by
\begin{equation}\label{aftt}
 y \, + \, \eps \Big( \frac{1}{2 \cdot 7^d } - \frac{1}{2\cdot 7^d \cdot 10^{10d}} -
 \frac{1}{2^{8d}}  - \frac{1}{2^{8d-1}}  - \frac{1}{5^{6d}} -
 \frac{1}{7^{5d}}  -   \frac{1}{3^{10d}} - \frac{\eps}{2^{8d-1} y} \Big),
\end{equation}
which strictly exceeds $y$, provided that $\eps < y$.

The lattice animal $\Phi$ may be formed for all sufficiently large
$n$. We conclude that
\begin{displaymath}
\liminf_{n \to \infty}{\frac{S(\Phi_n)}{n^d}} > y ,
\end{displaymath}
which implies that
\begin{displaymath}
\liminf_{n \to \infty}{\frac{G_n}{n^d}} > y ,
\end{displaymath}
an inconsistency which completes the proof. $\Box$

\section{Proof of Theorem \ref{pthmsix}}\label{secold}
We require a lemma.
\begin{lemma}\label{lemthm}
Let $P$ denote a percolation of parameter 
$p \in (p_c,1]$. 
For any $C \in \mathbb{N}$, there exists $n_0 \in \mathbb{N}$ such
that
\[
 \bigcap_{n \geq n_0}{P_{n,C}} \not= \emptyset,
\]
where the sets $P_{n,C}$ were specified in Definition \ref{bigconndef}.
\end{lemma}
\noindent{\bf Proof.} 
It follows from \cite[Theorem 7.2]{Gr} and the assumption that $p
>p_c$ that there exists almost surely an infinite 
cluster $P^{+}_{\infty}$ of the process $P \cap \Zd_{+}$, where $\Zd_{+} = \big\{ \bv \in \Zd: \bold v_i \geq 0, \, i
\in \{1,\ldots,d\} \big\}$. Let $\bold x \in P^{+}_{\infty}$. Note
that the connected component of $P \cap \Ba_n$ in which the site
$\bold x$ lies has radius at least $n - \| \bold x \|$.
Note that, for any $\al \in (0,1)$, if
the event $Q_n(\al)$ defined in \req{qnevent} occurs, and $n > (1-\al)^{-1}
\| \bold x \|$, then $\bold x$ lies in the connected set 
$C_n(\al)$, also defined in \req{qnevent}. 
Recall from the proof of Lemma \ref{lempva} that, if $\al \in (0,1)$
is small enough that $\theta(p) > 2\al^d + 4d\al$, then $C_n(\al) =
P_{n,0}$ for high values of $n$. Recalling also that $Q_n(\al)$ occurs
for all but finitely many $n$ almost surely, we deduce that $\bold x \in P_{n,0}$ for all
high choices of $n$.
The
statement of the lemma for a positive value of $C$ is obtained
by 
translating the process $P$ by the vector $(-C,\ldots,-C)$, and
applying the result for $C=0$. $\Box$ \\
\noindent{\bf Proof of Theorem \ref{pthmsix}:}
Given $\eps > 0$, let $C,\ell \in \mathbb{N}$ and the percolation $P$  
be those to which the statement of Theorem \ref{thmjkl} refers.
By Lemma \ref{lemthm}, there exists $F_0
\in \mathbb{N}$ such that we may 
choose $\bold v \in \bigcap_{F \geq F_0}{P_{F,C}}$. 
Let $\xi_n \in \LA_{\Ba_n}$ satisfy $S(\xi_n)
= G_n$ and $\vert \xi_n
\vert = L_n$. Provided that $n \geq F_0 \ell$ is also chosen to be so
high that we may apply Theorem \ref{thmjkl} to $\xi_n \subseteq \Ba_n$, 
we find that 
$\xi_n \cap \Ba_{\ell \bv,\ell} \not= \emptyset$. 
Let $\tau = ({\bf \tau_0},\ldots,{\bf \tau_r})$ denote a path in $\Zd_{+}$ such
that $\orig \in \tau$ and $\Ba_{\ell \bold v,\ell} \subseteq
\tau$. Let $V = \min \big\{ S(\tau^s) : \tau^s =
({\bf \tau_0},\ldots,{\bf \tau_s}) , \, s \in \{0,\ldots,r \} \big\}$ be equal
to the minimal weight of any initial subpath of $\tau$.
For each sufficiently high $n$, we may choose $s(n) \in
\{0,\ldots r\}$ such that $\tau^{s(n)} \subseteq \Ba_n \setminus \xi_n$ and
$\tau^{s(n)} \cap \partial \xi_n \not= \emptyset$.
Note that
\begin{equation}\label{brf}
N_{\vert \xi_n \vert + \vert \tau \vert} \geq S \big( \xi_n \cup \tau^{s(n)}
\big)
 =  S \big( \xi_n \big) + S\big( \tau^{s(n)} \big) \geq G_n + V.
\end{equation}
Given that $\vert \xi_n \vert =L_n$, we find from \req{lmineq} and  \cite[Theorem 2.1]{GLAthree} 
that
\begin{equation}\label{brgminus}
 N_{\vert \xi_n \vert + \vert \tau \vert} \leq (N + \eps)(L_n + \vert
 \tau \vert),
\end{equation}
for all $n$ sufficiently high. From \req{brf} and \req{brgminus}, we
deduce that
\begin{equation}\label{brg}
G_n \leq \big( N + \eps \big) L_n +   \big( N + \eps \big) \vert \tau
\vert - V.
\end{equation}
Given that $\eps > 0$ is arbitrary, and that the path $\tau$ is fixed,
we obtain, by taking a liminf of the $n^{-d}$-th multiple of
\req{brg}, the inequality $G \leq N L$ that we sought. $\Box$

\section{Critical behaviour, proof of Theorem \ref{thm-crit}}\label{secfour}

We aim to prove that the quantity $N$ is positive under the assumption
that, for
some $\epsilon > 0$,
\begin{equation}\label{thehyp}
   \limsup_{n \to \infty}{n^{-1}(\log n)^{-\frac{d}{d-1} -
\epsilon}G_n } > 0
\end{equation}
with positive probability. 
This limsup is non-random, similarly to $\limsup n^{-d} G_n$, as
explained at the beginning of the proof of Theorem \ref{pthmthree}. 
Hence, the hypothesis \req{thehyp} allows us to fix $\delta > 0$ for which
\begin{equation}\label{rep}
 \limsup_{n \to \infty}{ n^{-1} (\log n)^{-\frac{d}{d-1} -
\epsilon} G_n} > \delta  \, \, \, \textrm{almost surely.}  
\end{equation}
Let $\la_0 = \inf \{ \la \in \mathbb{R}: \PPP(X_{\orig} \geq - \la) >
p_c \}$. Recall from after \req{pbound} that, for $\la > \la_0$, we denote
by $\WW$ the unique infinite component of $\la$-white sites in $\Zd$. 
\begin{defin}
For $\la>\la_0$, $\rho>0$, $n \in \mathbb{N}$ and $\bold x \in \Zd$, the event
$E(\bold x,n,\la,\rho)$ occurs if there exists $\ga \in \LA_{\Ba_{\bold x,n}}$
such that $S(\ga) = G_{\Ba_{\bold x,n}} > \delta n (\log n)^{d/(d-1) +
\eps}$, $\vert \ga \vert = L_{\Ba_{\bold x,n}} >  \delta (\log n)^{d/(d-1) +
\eps}$, and a site $\bold v \in \ga \cap \WW$ satisfying $D(\bold u,\bold v) \leq \rho
\ell_{\infty} (\bold u,\bold v)$ for each $\bold u \in \WW \cap 
 \Ba_{\bold x,n}[1]^c$.
\end{defin}
\begin{lemma}\label{justgone}
For $\la,\rho$ sufficiently high, 
\begin{eqnarray}
 & & \big\{ G_n > \delta n (\log n)^{\frac{d}{d-1} + \epsilon} \
\textrm{occurs for infinitely many $n$} \big\} \nonumber \\
 = & & \big\{ E(\orig,n,\lambda,\rho) \ \textrm{occurs for infinitely many
 $n$} \big\}, \nonumber
\end{eqnarray}
up to a set of measure zero.
\end{lemma}
\noindent{\bf Proof:} We must show that, for high enough values of $n$, $G_n >
\delta n (\log n)^{d/(d-1) + \eps}$ implies the occurrence of
$E(\orig,n,\la,\rho)$ for given choices of $\la$ and $\rho$. 
As noted before \req{lmineq},
$X_{\bold v} \geq \| \bold v \|$ for at
most finitely many $\bold v \in \Zd$ almost surely, so that $ G_n > \delta n
(\log n)^{\frac{d}{d-1} + \epsilon}$ implies that 
\begin{equation}\label{vitn}
L_n > \delta (\log
n)^{\frac{d}{d-1} + \epsilon}
\end{equation}
 for high $n$. It follows from \req{rtu}, written with $\rho'$ in
 place of $\rho$,
and the Borel-Cantelli lemma, that for any $\rho' > 1$, and for  all
$n$ sufficiently high, each $\ga \in \LA_{\Ba_n}$
satisfying $\vert \ga \vert \geq (\log n)^{\rho'}$
intersects $\WW$. 
By \req{vitn}, each greedy lattice animal in $\Ba_n$ intersects $\WW$
for all high $n$. 
Let $\ga$ be  a greedy lattice
animal in $\Ba_n$, with $n$ chosen to be high enough that we may
locate a site $\bold v \in \WW \cap \ga$. 
If a site  $\bold u \in \Ba_n[1]^c \cap \WW$ satisfies
$D(\bold v,\bold u)
\geq \rho \ell_{\infty} (\bold v,\bold u)$, then $D(\bold v,\bold u) >
(\rho/4) \| \bold u \|$, since, setting $\bold c = (\lfloor n/2
\rfloor,\ldots,\lfloor n/2 \rfloor)$, 
\begin{eqnarray}
\ell_{\infty} (\bold u,\bold v) & = &
\ell_{\infty} (\bold u - \bold c,\bold v - \bold c) \geq 
\| \bold u - \bold c
\| - \| \bold v - \bold c 
\| \geq \| \bold u - \bold c \| - \lfloor n/2 \rfloor -
1 \nonumber \\
 & \geq & \frac{2}{3} \| \bold u  - \bold c \| -1  \geq
\frac{2}{3} \big( \| \bold u \| -  \| \bold c
\| \big) -1   \geq  
\frac{1}{3} \| \bold u \| -1, \nonumber
\end{eqnarray}
the third inequality valid by $\| \bold u - \bold c \|
\geq 3n/2$, and the fifth by $\| \bold u \| \geq n \geq 2 \| \bold c \|$.
By \cite[Lemma 2.14]{GLAthree}, with $\rho(p,d)$ set equal to $4 \rho$
as it is here, and a union bound, the probability that such a site
$\bold u$
exists is at most $\exp{-cn}$, for some positive constant $c$. The
Borel-Cantelli lemma implies that each $\bold u \in \Ba_n[1]^c \cap \WW$
satisfies $D(\bold u,\bold v) \leq \rho \ell_{\infty}(\bold u,\bold
v)$, 
provided that $n$ is
high enough. We have shown that $G_n > \delta n (\log
n)^{\frac{d}{d-1} + \epsilon}$ implies the occurence of $E(\orig,n,\lambda,\rho)$
for high values of $n$, as required. $\Box$

Defining the event $D(\bold x,m,\la,\rho,C)$, for $\bold x \in \Zd$ and
$C > 0$, according to
\begin{eqnarray}
 D(\bold x,m,\la,\rho,C) & = & \big\{ \exists \ga \in \LA_{\Ba_{\bold x,m}},
 \bold v \in \ga : \bold v \, \, \textrm{is $\la$-white}, \nonumber \\
 & & \qquad S(\ga) \geq C m, D(\bold v,\bold u) \leq
 \rho m \, \textrm{for each corner $\bold u$ of $\Ba_{\bold x,m}$} \big\}, \nonumber
\end{eqnarray}
we will now show that for any $\eps_0 > 0$, we may choose $\la,\rho$ and $C$
sufficiently high that
\begin{equation}\label{zyk}
\PPP \big( D(\orig,m,\la,\rho,C) \big) \geq 1 - \eps_0 \, ,
\end{equation}
for high values of $m$.
Let $c_1 \in (0,\infty)$ and $N_0 \in \mathbb{N}$ be chosen so that,
for $n \geq N_0$, 
\begin{equation}\label{cbone}
\big( \log n \big)^{\frac{d}{d-1} + \epsilon} > \frac{3 \big( 2 c_1 +
(2d-1) \lambda \rho \big)}{\delta \big( 1 - {\eps_0}/2 \big)}.
\end{equation}
By Lemma \ref{justgone} and \req{rep}, we may fix $N_1 > N_0$ for which
\begin{equation}\label{shn}
\mathbb{P}\big(E(\orig,n,\lambda,\rho) \ \textrm{occurs for some $n \in \{N_0, \ldots, N_1 \}$}\big) > 1 - \eps_0^2/4 . 
\end{equation}
Declare any site $\bold x \in \mathbb{Z}^d$ to be full if the event
$E(\bold x,n,\lambda,\rho)$ occurs for some $n \in \{N_0, \ldots, N_1 \}$. 
To any full site $\bold x$, we may associate a lattice animal
$\gamma_{\bold x}$,
a site $\bold v_{\bold x} \in \gamma_{\bold x}$ and the box
$\Gamma_{\bold x} = \Ba_{\bold x,n_{\bold x}}$,
these objects arising from the definition of the event
$E(\bold x,n_{\bold x},\lambda,\rho)$, for the minimal $n_{\bold x} \in \{N_0, \ldots, N_1 \}$
for which this event occurs.
Allowing $\bold e_1$ to denote the unit vector $(1,0,\ldots,0)$, we
set $\bold x_j
= j \bold e_1$, for $j \in \mathbb{N}$. We now form a subsequence $\{
\bold y_j :
j \in \mathbb{N} \}$ of the sequence $\{ \bold x_j : j \in \mathbb{N}
\}$. The first element $\bold y_1$ is taken to be $\bold x_j$, where $j$ is the
lowest natural number for which $\bold x_j$ is a full site. Having
constructed an initial segment of the $\bold y$-sequence, $\{ \bold y_j : j \in
\{1,\ldots, K \} \}$, say, we set $\bold y_{K+1}$ equal to the
lowest-labelled site in the $\bold x$-sequence which is full and has
$\bold e_1$-coordinate exceeding that of any site lying in the box
$\Gamma_{\bold y_K}[1]$.

Noting that $\bold v_{\bold y_{i+1}} \not\in \Gamma_{\bold y_i}[1]$, it follows from
the definition of the event $E(\bold y_i,n_{\bold y_i},\lambda,\rho)$ that
we may join $\bold v_{\bold y_i}$ and $\bold v_{\bold y_{i+1}}$ by a
path $\tau_i$ in $\WW$ of length at most $\rho
 \ell_{\infty}(\bold v_{\bold y_i}, \bold v_{\bold y_{i+1}})$. For each $J \in \mathbb{N}$, form the animal
\begin{displaymath}
\kappa_J = \gamma_{\bold y_1} \cup \tau_1 \cup \gamma_{\bold y_2} \cup \tau_2 \cup
\ldots \cup \tau_{J-1} \cup \gamma_{\bold y_J}.
\end{displaymath} 
Let $\xi$ be the collection of sites $\bold x_i$ that are not full and
that lie between
$\Gamma_{\bold y_j}[1]$ and $\bold y_{j+1}$ for some $j \in \mathbb{N}$, or before
$\bold y_1$. 
Writing $H  = \big\{ \frac{\vert \xi \cap \{ \bold x_1, \ldots, \bold x_m
\} \vert}{m} > {\eps_0}/2 \big\}$, we claim that, for any 
$m \in \mathbb{N}$,
\begin{equation}\label{weqnfour}
\mathbb{P} \big( H^c  \big) > 1 - {\eps_0}/2.
\end{equation}
To see this, we perform an experiment in which we
sample $z \in \{1,\ldots, m\}$ uniformly at random, and ask whether
the site $\bold x_z$ is full. If $\mathbb{P}(H^c) \leq 1 - {\eps_0}/2$, then 
\begin{displaymath}
\mathbb{P} \big( \bold x_z \ \textrm{is not full} \big) \geq \mathbb{P}
\big( \bold x_z \ \textrm{is
not full}
\big\vert H \big) \mathbb{P}(H) \geq \eps_0^2/4 ;
\end{displaymath} 
however, $\mathbb{P}(\bold x_z \ \textrm{is not full})$ is the 
probability that a
given site is not full, contradicting \req{shn} and establishing
\req{weqnfour}.

For $m \in \mathbb{N}$, let $J (= J(m))$ be maximal such that
$\Gamma_{\bold y_J}[1]$ has maximum $\bold e_1$-co-ordinate at most $m-1$. Let us
estimate the weight $S(\kappa_J)$ of the animal $\kappa_J$, for fixed
$m$. The animals $\gamma_j$ are disjoint for distinct $j$, and the
paths $\tau$ lie in $\WW$. Thus,
\begin{equation}\label{cbtwo}
S(\kappa_J) \geq \sum_{j=1}^{J}{S(\gamma_{\bold y_j})} - \lambda \sum_{j=1}^{J - 1
}{\vert \tau_j \vert}.
\end{equation}
In bounding the first term on the right-hand-side of \req{cbtwo},
note that,
for $j \in \{ 1, \ldots, J \}$, we have that 
$$
S \big( \gamma_{\bold y_j} \big) \geq 
\delta n_{\bold y_j} (\log n_{\bold y_j})^{\frac{d}{d-1} + \epsilon}.
$$
Since $n_{\bold y_i} \geq N_0$ for any such $j$, from \req{cbone}, it
follows that
\begin{equation}\label{cbthree}
\sum_{j = 1}^{J}{S \big( \gamma_{\bold y_j} \big) } 
\geq \frac{3 \big(2 c_1 + (2d-1)
\rho \lambda\big) }{1 - {\eps_0}/2} \sum_{j=1}^{J}{n_{\bold y_j}}.
\end{equation}
To bound from below the quantity $\sum_{j=1}^{J}{n_{\bold y_J}}$, note the
following inclusion:
\begin{equation}\label{cbfour}
\{\bold x_0,\ldots,\bold x_{m-1} \} \subseteq \Big(
\bigcup_{j=1}^{J}{\Gamma_{\bold y_j}[1]} \Big) \cup \xi \cup R,
\end{equation}
where the set $R$ denotes the final $2N_1 - 1$ sites of the interval $\{
\bold x_1,\ldots, \bold x_{m-1} \}$, and appears because of the possibility that the
site $\bold y_{J+1}$ lies in this interval. 
From \req{cbfour}, it follows that, on the event $H^c$,
\begin{equation}\label{cbfive}
\sum_{j=1}^{J}{n_{\bold y_j}} \geq \frac{m(1 - {\eps_0}/2)}{3} -
\frac{2N_1 - 1}{3}.
\end{equation}
We must also bound from above the quantity 
$\sum_{j=1}^{J-1}{\big\vert \tau_j \big\vert}$.
Note that 
\begin{eqnarray}
 \sum_{i=1}^{J-1}{\ell_{\infty} \big( \bold v_{\bold y_i},\bold
 v_{\bold y_{i+1}} \big)} & \leq & \sum_{i=1}^{J-1}
 \sum_{l=1}^{d}{\big\vert \bold v_{\bold y_i}^l - \bold v_{\bold y_{i+1}}^l \big\vert} \nonumber \\
 & \leq & \sum_{i=1}^{J-1}{\big\vert \bold v_{\bold y_i}^1 - \bold
 v_{\bold y_{i+1}}^1 \big\vert} +
 \big( d-1 \big) \sum_{i=1}^{J-1}{ \max \big\{ 
  n_{\bold y_i} ,  n_{\bold y_{i+1}} \big\}} \nonumber
 \\
 & \leq & \sum_{i=1}^{J-1}{\big\vert \bold v_{\bold y_i}^1 - \bold
 v_{\bold y_{i+1}}^1 \big\vert} 
 +
 2 \big( d-1 \big) \sum_{i=1}^{J}{n_{\bold y_i}} 
 \leq  \big( 2d - 1 \big) m, \label{anl}
\end{eqnarray}
where, in the second inequality, we used the fact that, for each $i \in
\mathbb{N}$ and $l \in
\big\{ 2,\ldots, d \big\}$,
\begin{displaymath}
\big\vert \bold v_{\bold y_i}^{l} - \bold v_{\bold y_{i+1}}^{l} \big\vert \leq  \max\big\{ 
 n_{\bold y_i}, n_{\bold y_{i+1}} \big\},
\end{displaymath}
while in the fourth, we used the bounds
\[
 \sum_{i=1}^{J-1}{\big\vert \bold v_{\bold y_i}^1 - \bold v_{\bold
 y_{i+1}}^1 \big\vert}  \leq  m, \, \,\textrm{and} \ \ 
  \sum_{i=1}^{J}{n_{\bold y_i}}  \leq  m.
\]
(In the first of these, we used the fact that $\bold v_{\bold y_i}^1$ is
increasing, which is true because the boxes $\Gamma_{\bold y_i}$ are
disjoint. The second also uses this disjointness).
From \req{anl} and $\vert \tau_i \vert \leq \rho
\ell_{\infty}(\bold v_{\bold y_i},\bold v_{\bold y_{i+1}})$, it follows that
\begin{equation}\label{cbsix}
\sum_{j=1}^{J-1}{\big\vert \tau_j \big\vert} \leq \big( 2d - 1 \big)
\rho m.
\end{equation}
Substituting the bounds \req{cbthree} and \req{cbsix} into \req{cbtwo} yields
\begin{displaymath}
S(\kappa_J) \geq  \frac{3 \big(2 c_1 + (2d-1)
\rho \lambda\big) }{1 - {\eps_0}/2} \sum_{j=1}^{J}{n_{\bold y_j}} \, - \,
\lambda \rho  (2d-1) m.
\end{displaymath}  
Substituting \req{cbfive} into this inequality, we see that, for high
values of $m$,
\begin{equation}\label{anm}
 H^c \subseteq \Big\{ S(\kappa_J) \geq 2  c_1 m - \frac{\big(2N_1 -
 1\big)\big( 2 c_1 + (2d-1) \rho \lambda\big)}{1 - {\eps_0}/2}
 \Big\} \subseteq \big\{  S(\kappa_J) \geq c_1 m \big\}.
\end{equation} 
We now claim that
\begin{equation}\label{bcz}
 H^c \subseteq \Big\{  \kappa_J \subseteq  \big\{ \bold v \in \Zd :
 \ell_{\infty}(\bold v,\Ba_m) \leq 
\rho \big( {m\eps_0}/2  + (d+2) N_1 \big) \big\} \Big\}.
\end{equation}
To show \req{bcz}, note that $\ga_{\bold y_j} \subseteq \Gamma_{\bold y_j}
\subseteq \Ba_m$ for each $j \in \{1,\ldots,J \}$, the latter inclusion valid
by the definition of $J=J(m)$. For $j \in \{ 1,\ldots,J-1 \}$ and
$\bold v
\in \tau_j$,
\begin{equation}\label{vasd}
 \ell_{\infty} (\bold v, \Ba_m) \leq \vert \tau_j \vert \leq \rho
 \ell_{\infty}(\bold v_{y_j}, \bold v_{y_{j+1}}) \leq \rho \big( \vert
 \bold v_{\bold y_j}^{1} - \bold v_{\bold y_{j+1}}^{1}
 \vert + (d-1) \max \{ n_{\bold y_j} , n_{\bold y_{j+1}} \} \big),
\end{equation}
the third inequality following similarly to \req{anl}. Note that
\begin{equation}\label{vase}
 \vert \bold v_{\bold y_j}^{1} - \bold v_{\bold y_{j+1}}^{1} \vert \leq \vert \xi
 \vert + 2 n_{\bold y_j} + n_{\bold y_{j+1}},
\end{equation}
because $\bold v_{\bold y_j}^{1}$ is at most $2 n_{\bold y_j}$ less than the
maximum $\bold e_1$-coordinate of the box $\Gamma_{\bold y_j}[1]$, $\bold
v_{\bold y_j}^1$ is at most $n_{\bold y_{j+1}}$ more than $\bold y_j$, the minimum
$\bold e_1$-coordinate of the box $\Gamma_{\bold y_j}$, while each
site $\bold x_{j}= j \bold e_1$ for which $j$ lies strictly between this
maximum and this minimum belongs to $\xi$. Given that $\vert \xi \vert
\leq m {\eps_0}/2$ on the event $H^c$, and that $\max \{
n_{\bold y_j},n_{\bold y_{j+1}} \} \leq N_1$, we find that \req{vasd} and \req{vase} imply \req{bcz}. 

By \req{anm}, \req{bcz} and the bound $\PPP (H^c) > 1 -
{\eps_0}/2$, we find, provided that $\eps_0$ has been chosen
so that $\eps_0 < 2 \rho^{-1}$, that, for $m$ sufficiently high,
\begin{eqnarray}
& & 
\PPP \Big( \exists \ga \in \LA_{\Ba_m[1]}, \bold v \in \ga \cap \WW :
S(\ga) \geq c_1 m, \nonumber \\ 
& & \qquad  \bold u \in \WW \cap \big\{ \bold x \in \Zd:
\ell_{\infty}(\bold x,\Ba_m) > m/2 \big\}
\implies D(\bold v, \bold u) \leq \rho \ell_{\infty}(\bold
v,\bold u) \Big) > 1 - {\eps_0}/2, \label{zxn}
\end{eqnarray}
the role of $\bold v$ in \req{zxn} being played by any $\bold v_{\bold y_i}$ for
$i \in \{1,\ldots,J \}$. (We are using the fact that $\Ba_{\bold
y_i,n_{\bold y_i}}[1] \subseteq \big\{ \bold x \in \Zd:
\ell_{\infty}(\bold x,\Ba_m) \leq m/2 \big\}$, which is implied,
provided that $m \geq 2N_1$, by 
$\Ba_{\bold y_i,n_{\bold y_i}} \subseteq \Ba_m$ and $N_1 \geq n_{\bold
y_i}$.) Note also that
\begin{equation}\label{zxd}
 \PPP \Big( \, \textrm{each corner of $\Ba_m[1]$ lies in $\WW$} \Big)
 \geq 1 - 2^d \PPP ( \orig \not\in \WW ) \geq 1 - {\eps_0}/2,
\end{equation}
since \req{thetfact} permits us to
choose $\lambda \in \mathbb{R}$ so that the second inequality is
valid.
By \req{zxn}, \req{zxd} and the translation invariance of the process
$\{ X_{\bv}: \bv \in \Zd \}$, we find that, for $m$ large and
divisible by three,
\[
\PPP \Big( \exists \ga \in \LA_{\Ba_m}, \bold v \in \ga \cap \WW : \, \,
S(\ga) \geq (c_1/3) m, \, D(\bold v, \bold u) \leq \rho m \, \textrm{for
each corner $\bold u$ of $\Ba_m$} \Big) > 1 - \eps_0,
\]
The condition that $m$ is divisible by three occurs because the
sidelength of the box $\Ba_m [1]$ must satisfy this. It may be dropped
by replacing $\Ba_m[1]$ in \req{zxd} by a box that extends by one or
two sites further on one half of its faces.
Writing 
$c_1 = 3C$, we have shown \req{zyk}.

By the proof of Lemma \ref{lemact}, the process 
\begin{equation}\label{ahs}
 \Big\{ \bold a \in \Zd : D (m \bold a,m,\la,\rho,C) \, \,  \textrm{occurs} \Big\}
\end{equation}
is a $(2 \rho + 1)$-near percolation, for any given $m \in
\mathbb{N}$. By \req{zyk}, Lemma \ref{lemlss} and \req{thetfact}, 
we may fix $\la,\rho,C>0$ and a high value of $m$ so that there exists
a subset $P$ of \req{ahs} that is a percolation of supercritical
parameter. Let $\{ \bold a_i : i \in \mathbb{N} \}$ denote an infinite
self-avoiding path in
$P$. Note that, for each $i \in \mathbb{N}$, $D(\ga_{\bold
a_i},\ga_{\bold a_{i+1}}) \leq 2\rho m + 1$, where $\ga_{\bold a}
\subseteq \Ba_{m \bold a,m}$ and $\bold v_{\bold a} \in \ga_{\bold a}$
denote the
lattice animal and the site therein resulting from
the occurence of $D(m \bold a,m,\la,\rho,C)$: the inequality is due to  
$D(\ga_{\bold a_i},\ga_{\bold a_{i+1}}) \leq D(\bold v_{a_i},\bold w)
+  D(\bold w',\bold v_{a_{i+1}}) + 1$, where $\bold w$ $\bold w'$ 
are adjacent corners
of $\Ba_{m \bold a_i,m}$ and $\Ba_{m \bold a_{i+1},m}$. Let
$\phi_i$ denote a white path of length at most $2 \rho m + 1$ from
$\ga_{\bold a_i}$ to $\ga_{\bold a_{i+1}}$. Let $\phi_0$ denote an
arbitrary path from $\orig$ to $\ga_{\bold a_1}$. We form an increasing
sequence of lattice animals $\{R_i : i \in \mathbb{N} \}$, satisfying
$\vert R_i \vert = i$, with $R_0 = \{ \orig \}$, which successively
collect the sites of the path $\phi_0$, then the animal $\ga_{\bold
a_i}$ and the path $\phi_i$, for each $i \in \{1,2,\ldots\}$ in
turn.

Let $n_i = \inf \{ j \in \mathbb{N}: R_j \supset \ga_{\bold a_i}
\}$. Then
\[
R_{n_i} = \phi_0 \, \cup \, \bigcup_{j=1}^{i-1}{\big( \ga_{\bold a_j}
\cup \phi_j \big)} \, \cup \, \ga_{\bold a_i}.
\]
We find that
\begin{equation}\label{gne}
S(R_{n_i}) \geq \sum_{j=1}^{i}{S(\ga_{\bold a_j})} \, - \, S(\phi_0)
\, - \, \la \sum_{j=1}^{i-1}{\big( \vert \phi_j \vert - 2 \big)} \geq
C m i  - \la(2\rho m -1)(i - 1) - S(\phi_0),
\end{equation}
where the first inequality follows from the animals $\ga_{\bold a_j}
\subseteq \Ba_{m \bold a_j,m}$ being disjoint and the paths $\phi_j$
being white, their endpoints lying in $\ga_{\bold a_j}$ or 
$\ga_{\bold a_{j+1}}$. 
Note also that
\begin{equation}\label{gen}
n_i \leq \vert \phi_0 \vert + m^d i + (2\rho m -1)(i-1),
\end{equation}
since $\vert \ga_{\bold a_j} \vert \leq m^d$ and $\vert \phi_{\bold
a_j} \vert - 2 \leq 2 \rho m - 1$. By \req{gne}, \req{gen} and $\vert
R_{n_i}\vert = n_i$, we obtain
\[
 \liminf_{i \to \infty} \frac{S(R_{n_i})}{n_i} \geq \frac{Cm - \la
 (2\rho m -1)}{m^d + 2 \rho m - 1} > 0 \, ,
\]
provided that the constant $C$ is chosen so that $C > 2\la\rho$.
Thus, on our hypothesis \req{thehyp}, $N > 0$, by the definition of
$N$. This completes the proof. $\Box$ \\

\section{Appendix: the proofs of Proposition \ref{thm-conc} and
Corollary \ref{webcor}}\label{secfive}
We break Proposition \ref{thm-conc} into three parts, stated now as
lemmas. We will prove each of the three lemmas in turn.
\begin{lemma}\label{pthmtwo}
For any law $F$ satisfying condition \req{martins}, the function
$\tilde{G}(\al):= \liminf_n n^{-d} \tilde{G}_n(\al)$ is finite valued, non-random
and concave on $(0,1)$.
\end{lemma}
\begin{lemma}\label{pthmfour}
For any distribution $F$ satisfying condition \req{martins},
the limit $\lim_{n \to \infty}{n^{-d} \tilde{G}_n(\al)}$ exists almost surely,
for each $\al \in (0,1)$.
\end{lemma}
\begin{lemma}\label{pthmfive}
For a distribution $F$ satisfying condition \req{martins}
which is bounded below, we have that $\al^{-1} \tilde{G}(\al) \to N$ as
$\al \downarrow 0$.
\end{lemma}
\noindent{\bf Proof of Lemma \ref{pthmtwo}:}
For $\alpha \in (0,1)$, $\tilde{G}(\alpha)$ is a tail event for 
the collection $\{ X_{\bv} : \bv \in \Zd \}$ (see the start of the
proof of Theorem \ref{pthmthree} for a formal statement). By Kolmogorov's
zero-one law, $\tilde{G}(\alpha)$ is an almost sure constant. Note that $\tilde{G}(\alpha)
\leq \liminf_n G_n = G < \infty$, the last inequality by
\cite[Theorem 3.2]{GLAthree}.  

We adapt the notion that an $m$-box contains a greedy lattice animal
that is well connected to large animals in its surroundings to the
current context where that `greedy' animal must comprise a fixed
fraction of the sites of the box, merely adjusting the form of
condition  ${\bf A}_{\la,c}^{\rho}$ given in Lemma \ref{lemmainr} 
to include this requirement: 
\begin{defin}
For $\al \in (0,1)$, write $\tilde{G}_{\Ba}(\al) = \max \{ S(\xi): \xi \in
\LA_{\Ba}, \, \vert \xi \vert = \lfloor \al m^d \rfloor  \}$ 
for any $m$-box $\Ba = \Ba_{\bold x,m}$. 
Given constants $c>0$ and $\la,\rho < \infty$, we
say that the box $\Ba$ satisfies condition
${\bf A}_{\la,c}^{\rho}[\alpha]$ if there exists 
$\ga^{*} \in \LA_{\Ba}$ satisfying $S(\ga^{*}) = G_{\Ba}(\alpha)$,
$\vert \ga^{*} \vert = \lfloor m^d \alpha \rfloor$ and $D(\ga,\ga^{*})
\leq \rho m$ for all $\ga \in \LA_{\Ba[2]}$ such that $\vert \ga \vert
\geq (\log m)^{\rho}$.
\end{defin}
\begin{lemma}\label{lemff}
Set
$
q_{m,\la,c,\rho}(\alpha) := \PPP(\Ba \ \textrm{satisfies condition}
\ {\bf A}_{\la,c}^{\rho}[\alpha] )$,
which is independent of the choice of the $m$-box $\Ba$.
Then, for $\alpha \in (0,1)$, $c \in (-\infty,\tilde{G}(\alpha) )$, and $C
\leq \la,\rho < \infty$, with $C>2$ a universal constant, we have that
\[
\lim_{m \to \infty} q_{m,\la,c,\rho}(\alpha)   = 1 .
\]
\end{lemma}
\noindent{\bf Proof:} A few alterations to the proof of Lemma \ref{lemmainr}
are required. Taking $\ga^{*} \in \LA_{\Ba}$ that satisfies  $\lfloor \ga^* \rfloor =
\lfloor \alpha m^d \rfloor$ and $S(\ga^*) = \tilde{G}_{\Ba}(\alpha)$, 
we find from \req{ubdnew} that if $m$ is high enough that $\al m^d
\geq (\log m)^{\rho}$ + 1, then
\begin{equation}\label{qhj}
 \PPP \big( \ga^* \cap \WW  = \emptyset \big) \leq \eps. 
\end{equation}
By \req{qhj} and another use of \req{ubdnew}, we find that \req{tosoh} holds for the current
choice of $\ga^*$. We deduce \req{feqnthree}
from \req{tosoh} and  \req{lemtwod} as in the earlier proof. 
 That $\PPP \big( S(\ga^*) = \tilde{G}_{\Ba}(\alpha) \geq c m^d \big) >
1 - \eps$ follows from $\tilde{G}(\alpha)$ being non-random and $c <
\tilde{G}(\alpha)$. From this and \req{feqnthree}, the statement of the lemma 
follows. $\Box$ \\
We now fix $\alpha_1,\alpha_2 \in (0,1)$ satisfying $\alpha_1 <
\alpha_2$, and $p \in (0,1)$.
We aim to show that 
\begin{equation}\label{geeeqn}
 \tilde{G} \big( p \alpha_1 + (1-p)\alpha_2 \big) \geq p \tilde{G}(\alpha_1) +
 (1-p)\tilde{G}(\alpha_2),
\end{equation}
and begin by choosing $\eps > 0$ that satisfies 
\begin{equation}\label{epscond}
(1- \eps)\alpha_2
> p \alpha_1 + (1-p)\alpha_2.
\end{equation}

By Lemma \ref{lemff}, we choose $\la<\infty$ and $2<\rho<\infty$ such that,
for each $\alpha \in (0,1)$ and $\eps_1 > 0$, 
\[
\lim_{m \to \infty}{q_{m,\la,\tilde{G}(\alpha) -
\eps_1,\rho}(\alpha)} = 1.
\]
Adapting Definition \ref{deflact}, we say that $\ba \in \Zd$ is
$(\ell,\alpha)$-active if the $\ell$-box $\Ba_{\ell \ba,\ell}$
satisfies condition ${\bf A}_{\la,c}^{\rho}[\alpha]$. 
Recall from the paragraph after Lemma \ref{lemlss} that 
we may choose $\delta > 0$ such that any  $(2 \rho + 1)$-near
percolation of parameter exceeding $1 - \delta$ contains a percolation
of parameter $p$ for which 
\begin{equation}\label{hatit}
\theta(p) > 1 - \eps.
\end{equation}
 We fix $\ell \in
\mathbb{N}$ such that 
\begin{equation}\label{stit}
q_{\ell,\la,\tilde{G}(\alpha_i) -\eps_1,\rho} (\alpha_i) > 1 - \delta/2 \, \,
\textrm{for $i \in \{ 1,2 \} $,}
\end{equation}
stating other required lower bounds on $\ell$ as they arise.
Having fixed $\ell$, we write $\alpha$-active in place of
$(\ell,\alpha)$-active from now on.  
Note that the argument in the proof of Lemma \ref{lemact} shows that
for any $\ba \in \Zd$ and for any $\al \in (0,1)$, the event $\big\{
B_{\ell \ba,\ell} \, \, \textrm{satisfies condition ${\bf
A}_{\la,c}^{\rho}[\alpha]$} \big\}$ has the same property of
measurability as does that in \req{evb} (to use this argument, we have
required that
$\rho > 2$). Thus, the collection of sites
that are both $\al_1$- and $\al_2$-active is 
a $(2 \rho + 1)$-near
percolation whose parameter exceeds $1 - \delta$ by \req{stit}. 
We may find a percolation $P$ of parameter $p$ such that 
$P \subseteq \big\{ \ba \in \Zd : \ba \, \, \textrm{is $\alpha_1$- and
$\alpha_2$-active}\big\}$. By Lemma
\ref{lempva} and \req{hatit}, there exists $F_0 \in \mathbb{N}$ such that, for $F \geq F_0$,
\begin{equation}\label{adf}
  \big\vert P_{F,\lfloor \rho \rfloor + 1} \big\vert \geq F^d \big( 1
  - \eps \big).
\end{equation}
Let $n \geq \ell F_0$ be written in the form $n = F \ell + r, r \in
\{0,\ldots,F-1\}$, with $F \geq F_0$. For $\ba \in \Zd$
$\alpha$-active, let $\ga_{\ba}(\alpha)$ denote a lattice animal
playing the role of $\ga^*$ in the  condition  ${\bf
A}_{\la,c}^{\rho}[\alpha]$ that $\Ba_{\ell \ba,\ell}$ satisfies. If
$\bold a_1,\bold a_2 \in P$ are adjacent, then for any choice of the
pair $i,j \in \{1,2\}$, we may find a $\la$-white path
$\phi^{i,j}_{\bold a_1,\bold a_2}$ satisfying $\vert \phi^{i,j}_{\bold
a_1,\bold a_2} \vert \leq \rho \ell$ that connects $\ga_{\bold
a_1}(\alpha_i)$ and $\ga_{\bold
a_2}(\alpha_j)$, since  $\vert \ga_{\bold
a_1}(\alpha_i) \vert = \lfloor \alpha_1 \ell^d \rfloor \geq (\log \ell
)^\rho$, and similarly for  $\vert \ga_{\bold
a_1}(\alpha_j) \vert$, provided that $\ell$ has been fixed at a high enough value. Consider the two lattice animals
\[
 \Psi_1 = \bigcup_{\ba \in P_{F,\lfloor \rho \rfloor + 1}}{\ga_{\ba}(\alpha_1)}
 \ \cup \bigcup_{\bold a_1,\bold a_2 \in  P_{F,\lfloor \rho \rfloor + 1}:
 \vert\bold a_1 - \bold a_2 \vert  =1 }{\phi^{1,1}_{\bold a_1,\bold a_2}},
\]
and
\[
 \Psi_2 = \bigcup_{\ba \in P_{F,\lfloor \rho \rfloor + 1}}{\ga_{\ba}(\alpha_2)}
 \ \cup \bigcup_{\bold a_1,\bold a_2 \in  P_{F,\lfloor \rho \rfloor + 1}:
 \vert\bold a_1 - \bold a_2 \vert  =1 }{\phi^{2,2}_{\bold a_1,\bold a_2}}.
\]
That $\Psi_1$ and $\Psi_2$ are lattice animals in $\Ba_n$ follows by
the argument,  given after \req{eqnpsi},  that shows this for $\Psi$. 
Note that 
\begin{eqnarray}
 \vert \Psi_1 \vert & \leq & \vert P_{F,\lfloor \rho \rfloor + 1} \vert
 \lfloor \alpha_1 \ell^d \rfloor + \big\vert \big\{ \{\bold a_1,\bold
 a_2\}: \bold a_1,\bold a_2 \in P_{F,\lfloor \rho \rfloor + 1}, \vert \bold a_1 - \bold a_2
 \vert = 1 \big\} \big\vert \rho \ell \nonumber \\ 
 & \leq & \alpha_1 F^d \ell^d +
 d\rho F^d \ell \leq \alpha_1 n^d + d \rho \frac{n^d}{\ell^{d-1}}, \nonumber
\end{eqnarray}
where $\vert P_{F,\lfloor \rho \rfloor + 1} \vert \leq F^d$ was used
in the second inequality. Note also that
\[
 \vert \Psi_2 \vert \geq \vert P_{F,\lfloor \rho \rfloor + 1} \vert
 \lfloor \alpha_2 \ell^d \rfloor \geq
 F^d \big( \alpha_2 \ell^d - 1 \big) (1 - \eps) \geq \alpha_2
 (1-\eps)(n-\ell)^d - \frac{n^d}{\ell^d} (1 - \eps),
\]
the second inequality in this case valid by \req{adf}.
Thus, provided that $\ell$ is chosen so that $\ell \geq
C(\eps,\al_1,\al_2,p,\rho,d)$, then, for $n$ satisfying $n
\geq n(\ell,\eps,\al_1,\al_2,p,\rho,d)$, \req{epscond} implies that 
\begin{equation}\label{xhj}
 \big\vert \Psi_1 \big\vert \leq n^d \big( p \alpha_1 + (1-p)\alpha_2
 \big) \leq  \big\vert \Psi_2 \big\vert.
\end{equation}
For each subset $R \subseteq P_{F,\lfloor \rho \rfloor +1}$, we
define
\begin{eqnarray}
 \Psi [R] & = & \bigcup_{\ba \in R}{\ga_{\ba}(\alpha_1)} \, \, \cup
  \bigcup_{\ba \in  P_{F,\lfloor \rho \rfloor + 1} \cap
  R^c}{\ga_{\ba}(\alpha_2)}  \, \,  \cup 
 \bigcup_{\bold a_1,\bold a_2 \in  R:
 \vert\bold a_1 - \bold a_2 \vert  =1 }{\phi^{1,1}_{\bold a_1,\bold
 a_2}}  \nonumber \\
 & & \qquad  \cup  \bigcup_{\bold a_1 \in  R,\bold a_2 \in  P_{F,\lfloor \rho \rfloor + 1} \cap
  R^c :
 \vert\bold a_1 - \bold a_2 \vert  =1 }{\phi^{1,2}_{\bold a_1,\bold
 a_2}} \, \, \cup  \bigcup_{\bold a_1,\bold a_2 \in  P_{F,\lfloor \rho \rfloor + 1} \cap
  R^c :
 \vert\bold a_1 - \bold a_2 \vert  =1 }{\phi^{2,2}_{\bold a_1,\bold
 a_2}}. \label{powp}
\end{eqnarray}
For each such $R$, $\Psi[R]$ is a lattice animal in $\Ba_n$, by the
argument that applied to $\Psi_1$ and $\Psi_2$. We write $\big(
R_0,R_1,\ldots,R_M\big)$, $M = \vert P_{F,\lfloor \rho \rfloor + 1}\vert$,
for an arbitrary sequence satisfying $R_i \subseteq  P_{F,\lfloor \rho
\rfloor + 1}$, $\vert R_i \vert = i$, $R_i \subseteq R_{i+1}$ for $i
\in \{ 0,\ldots,M-1 \}$. 
Note that if $R_{i+1} \setminus R_i = \{ \bold a \}$, 
$\Psi[R_{i+1}]$ may be formed from
$\Psi[R_i]$ 
in two steps. In the first one, we remove from 
$\Psi[R_i]$  that part of 
\[
 \ga_{\bold a}(\al_1) \, \cup  \bigcup_{\bold a' \in  R_i:
 \vert\bold a - \bold a' \vert  =1 }{\phi^{1,1}_{\bold a,\bold
 a'}} \, \cup  \bigcup_{\bold a' \in  P_{F,\lfloor \rho \rfloor + 1}
 \cap R_i^c :
 \vert \bold a - \bold a' \vert  = 1}{\phi^{1,2}_{\bold a,\bold
 a'}}
\]
that is disjoint from the union of the other sets in the expression for
$\Psi[R_i]$ given by
\req{powp}. In the second step, we add to this altered set the
sites of
\[
 \ga_{\bold a}(\al_2) \, \cup  \bigcup_{\bold a' \in  R_{i+1}:
 \vert\bold a - \bold a' \vert  =1 }{\phi^{2,1}_{\bold a,\bold
 a'}} \, \cup  \bigcup_{\bold a' \in  P_{F,\lfloor \rho \rfloor + 1}
 \cap R_{i+1}^c :
 \vert \bold a - \bold a' \vert  = 1}{\phi^{2,2}_{\bold a,\bold a'}}
\] 
to obtain 
$\Psi[R_{i+1}]$. From this, we see that
\begin{equation}\label{halfp}
 \big\vert \Psi[R_{i+1}] \big\vert -  \big\vert \Psi[R_i] \big\vert
 \leq  \vert \ga_{\bold a}(\alpha_2) \vert 
 + 2d\rho \ell \leq \alpha_2 \ell^d + 2d \rho \ell , 
\end{equation}
for this second term is an upper bound on the number of sites
gained at the second step.
Thus, provided that $\ell \geq \big( 2 d \rho/\alpha_2
\big)^{1/(d-1)}$, by \req{xhj}, the facts that $\Psi_1 = \Psi[R_0]$,
$\Psi_2 = \Psi[R_M]$ and \req{halfp}, 
we may find $k \in \{0,\ldots, M \}$ for which
\begin{equation}\label{rfg}
 \vert  \Psi[R_k] \vert  \in  \big\{ \lfloor n^d \big( p \alpha_1 +
 (1-p)\alpha_2 \big) \rfloor - \lfloor 2 \alpha_2  \ell^d \rfloor - 1, \ldots,
   \lfloor n^d \big( p \alpha_1 +
 (1-p)\alpha_2 \big) \rfloor  \big\}.
\end{equation}
We now bound the value of $k$. Note that
\begin{eqnarray}
 \vert \Psi[R_k] \vert & \geq & \sum_{\bold a \in R_k} \vert
 \ga_{\ba}(\alpha_1) \vert \, + \, \sum_{\ba \in  P_{F,\lfloor \rho
\rfloor + 1} \cap R_k^c}{ \vert
 \ga_{\ba}(\alpha_2) \vert} \nonumber \\
 & \geq & k \lfloor \alpha_1 \ell^d \rfloor + \big( \big\vert
  P_{F,\lfloor \rho
\rfloor + 1} \big\vert - k \big) \lfloor \alpha_2 \ell^d \rfloor \geq
k \big(  \lfloor \alpha_1 \ell^d \rfloor - \lfloor \alpha_2 \ell^d
\rfloor \big)\, + \, (1-\eps)F^d  \lfloor \alpha_2 \ell^d \rfloor, \label{rfh}
\end{eqnarray}
where in the first inequality, we used the disjointness of
$\ga_{\bold a}(\alpha_i)$ for distinct $\bold a$, with
$k = \vert R_k \vert$ being applied in the second inequality and \req{adf} in the
third. Note also that
\begin{eqnarray}
 \vert \Psi[R_k] \vert & \leq & \sum_{\bold a \in R_k} \vert
 \ga_{\ba}(\alpha_1) \vert \, + \, \sum_{\ba \in  P_{F,\lfloor \rho
\rfloor + 1} \cap R_k^c}{ \vert
 \ga_{\ba}(\alpha_2) \vert}  \, + \, d\rho\ell F^d \label{rfi} \\
 & \leq & k \lfloor \alpha_1 \ell^d \rfloor + \big( F^d - k \big)
 \lfloor \alpha_2 \ell^d \rfloor \, + \, d\rho \ell F^d, \nonumber
\end{eqnarray}
where in the first inequality, we used the fact the $\phi$-paths in
$\Psi[R_k]$ are each of length at most $\rho \ell$ and are indexed by
pairs of adjacent sites in $P_{F,\lfloor \rho \rfloor +1}$.
By \req{rfg} and \req{rfh},
\begin{equation}\label{klb}
 k \geq \Big( p - \frac{\eps \alpha_2}{\alpha_2 - \alpha_1} \Big) \big( n/\ell
 \big)^d \, - \, o(n^d).
\end{equation}
By \req{rfg} and \req{rfi},
\begin{equation}\label{kub}
 k \leq p \big( n/\ell \big)^d \, +
 \, o(n^d).
\end{equation}
By \req{rfg}, a set $E$ of at most $2 \alpha_2 \ell^d + 1$
sites of $\Ba_n$ may be needed to be added to $\Psi[R_k]$ in such a
way that
\begin{equation}\label{sgtyk}
\Psi[R_k] \cup E \in \LA_{\Ba_n} \, \, \textrm{and} \, \, \big\vert \Psi[R_k] \cup E \big\vert = \lfloor n^d \big( p \alpha_1 +
 (1-p) \alpha_2 \big) \rfloor.
\end{equation}
We now show that we may choose $E$ from a fixed set of sites, so that
its weight is bounded below, uniformly in $n$. Note that, for any $\ba
\in \Zd$,
\begin{eqnarray}
 \big\vert \Ba_{\ell \bold a,\ell} \cap  \Psi[R_k]  \big\vert & \leq &
 \max\big\{ \lfloor \ga_{\bold a}(\al_1) \rfloor, \lfloor \ga_{\bold
 a}(\al_2) \rfloor \big\}  \nonumber \\
 & & + \, \sum\Big\{ \max_{i,j \in
 \{1,2\}}{\vert \phi_{\bold a_1,\bold a_2}^{i,j} \vert}: \big\{ \bold
 a_1,\bold a_2 \big\} \, \, \textrm{such that $\| \bold a_1 - \bold
 a_2 \| \leq 1$,} \nonumber \\
 & & \qquad \qquad \qquad \qquad \qquad   \qquad \qquad \qquad  \textrm{and $\| \bold a_l - \bold a
 \| \leq \lfloor \rho \rfloor + 1$ for $l \in \{1,2\}$} \Big\} \nonumber
 \\
 & \leq & \lfloor \al_2 \ell^d \rfloor \, + \, \big( 2 \lfloor \rho
 \rfloor + 3 \big)^d d \rho \ell. \label{crct}
\end{eqnarray} 
That is, the part of $\Psi[R_k]$ that may intersect a given $\ell$-box
$\Ba_{\ell \bold a,\ell}$ consists of the animal $\ga_{\bold
a}(\al_1)$ or $\ga_{\bold
a}(\al_2)$ in that box, and some part of the $\phi$-paths that connect
neighbouring boxes 
$\Ba_{\ell \bold a_1,\ell}$ and 
$\Ba_{\ell \bold a_2,\ell}$ for which $\| \bold a_i - \bold a \| \leq
\lfloor \rho \rfloor + 1$ for $i \in \{1,2\}$. (It suffices to consider
only these $\phi$-paths, because each such path has length at most
$\rho \ell$.) By choosing $\ell$ to satisfy 
\[
 \ell \geq \Big( \frac{2d\rho}{1 - \al_2 }\Big)^{\frac{1}{d-1}} (2\rho
 + 3)^{\frac{d}{d-1}},
\]
we obtain from \req{crct} that, for each $\bold a \in \Zd$,
\begin{equation}\label{scrto}
 \big\vert \Ba_{\ell \bold a,\ell} \cap \Psi[R_k] \big\vert \leq \Big(
 \frac{1}{2} + \frac{\al_2}{2}\Big) \ell^d.
\end{equation}
By Lemma \ref{lemthm}, we may select a collection $\big\{ \bold
b_1,\ldots,\bold b_D \big\}$, with $D \geq \frac{6\al_2}{1+\al_2}$,
such that $\bold b_i \in P_{F,\lfloor \rho \rfloor + 1}$ for all $n$
(and thus, $F$) sufficiently high.

By \req{scrto}, the lower bound on $D$, and a choice of $\ell$ for
which $\ell \geq \al_2^{-d}$, there are at least $2\al_2
\ell^d +1$ sites in the set 
\[
 \Big( \bigcup_{i=1}^D{\Ba_{\ell \bold b_i,\ell} }\Big) \setminus \Psi[R_k].
\]
Choosing a set $E$ that lies in this set and satisfies \req{sgtyk}, we
note that
\[
 S(E) \geq W,
\] 
where $W = \min \big\{ S(\ga): \ga \subseteq \bigcup_{i=1}^D{\Ba_{\ell
\bold a_i,\ell}} \big\}$.
 We bound
\begin{equation}\label{poli}
 S \big( \Psi[R_k] \cup E \big) =   S \big( \Psi[R_k] \big) + S \big(
 E \big) \geq   S \big( \Psi[R_k] \big)  + W. 
\end{equation}
In this way, we obtain the necessary small adjustment in the size of
$\Psi[R_k]$ by means of the set $E$, whose weight is uniformly bounded
below in $n$.

Note that 
\begin{eqnarray}
 S \big( \Psi[R_k] \big) & \geq & \sum_{\bold a \in R_k} S
 \big( \ga_{\ba}(\alpha_1) \big) \, + \, \sum_{\ba \in  P_{F,\lfloor \rho
\rfloor + 1} \cap R_k^c}{ S \big( \ga_{\ba}(\alpha_2) \big) } \, - \,
d\la \rho \ell F^d \nonumber \\
 & \geq &  k \big(  \tilde{G}(\alpha_1) - \epsilon_1 \big) \ell^d + \big(  \vert
 P_{F,\lfloor \rho
\rfloor + 1} \vert - k \big)  \big(  \tilde{G}(\alpha_2) - \epsilon_1 \big)
\ell^d \, - \, 
d\la \rho \ell F^d \nonumber \\
 & \geq & \big( p - \eps \alpha_2/(\alpha_2 - \alpha_1) \big) \big(  \tilde{G}(\alpha_1)
 - \epsilon_1 \big) n^d \, + \, \big( (1-\eps)(F \ell)^d - p n^d \big) 
 \big(  \tilde{G}(\alpha_2) - \epsilon_1 \big) - o(n^d) \nonumber \\
 & \geq & \big( p \tilde{G}(\alpha_1)
 + (1-p)\tilde{G}(\alpha_2) \big)n^d + (\eps +
 \eps_1) O(n^d) - o(n^d) \label{alr},
\end{eqnarray}
by the fact that the $\phi$-paths in $\Psi[R_k]$ are $\la$-white in
the first inequality, and
by \req{adf}, \req{klb} and \req{kub} in the third. 
From \req{sgtyk}, \req{poli}, \req{alr}, the definition of the
function $\tilde{G}:(0,1) \to \mathbb{R}$, and the fact that $\eps,\eps_1>0$
may be chosen to be arbitrarily small, we obtain \req{geeeqn}. $\Box$

\noindent{\bf Proof of Lemma \ref{pthmfour}:} 
We mimic the proof of Theorem \ref{pthmthree}. Similarly to that which was noted
at the start of the proof of Lemma \ref{pthmtwo}, if $y = \tilde{G}(\alpha)$ and $y +
E = \limsup n^{-d} \tilde{G}_n (\alpha)$, then $y$ and $E$ are
non-random. Recalling from the proof of Lemma \ref{pthmtwo} that the
collection of $\al$-active sites forms a $(2 \rho + 1)$-near
percolation, we may find, as we did in the second paragraph of the
proof of Theorem \ref{pthmthree}, 
a percolation $P$ whose parameter $p$ satisfies \req{thpineq},
and for which 
$$
P \subseteq \big\{ \bold a \in \Zd : \bold a \, \, \textrm{is} \, \, \al-\textrm{active} \big\}.
$$
In the current case, the animal $\hat{\Psi}$ is defined by
\req{gbnm}, with instances of animals $\ga_{\ba}$ being replaced
by $\ga_{\ba}(\al)$, and the $\phi$-paths existing by use of condition
${\bf A}^{\rho}_{\la,y + \eps}[\alpha]$ in place of  ${\bf A}_{\la,y +
\eps}^{\rho}$.
We alter the definition of an $m$-box being high so that its first
requirement takes the form of condition $A_{\la,y + \eps}^{\rho}
[\al]$, the second requirement being unchanged, so that now $\ga^{*}_{\Ba_{\bold x,m}}
= \lfloor \al m^d \rfloor$ if the box $\Ba_{\bold x,m}$ is high. We obtain as in the proof of
Theorem \ref{pthmthree}, a sequence of lattice animals $\Phi_n \subseteq \Ba_n$
satisfying $\liminf n^{-d} S(\Phi_n) > y$. A contradiction has not
been reached however, because $\vert \Phi_n \vert$ may not be equal to
$\vert \al n^d \vert$.
By \req{ghnl} and the analogue of \req{gbnm}, we find that
\begin{equation}\label{vaxo}
 \big\vert \Phi \big\vert \geq \sum_{\bold a \in BD}{\vert \ga_{\bold
 a} \vert} \, + \, \sum_{\Gamma \in UH}{\vert \ga_{\Gamma}^* \vert} \, ,
\end{equation}
because the sets $\bigcup \big\{ \Ba_{\ell \bold a,\ell}:
\bold a \in BD \big\}$ and $\bigcup \big\{ \Gamma : \Gamma \in UH \big\}$ are
disjoint (which follows from $BD$ being a connected component of
$P_{F,\lfloor \rho \rfloor + 1} - \cup_{\Gamma \in
\kappa}{w_{\Gamma}}$). 
We also have that
\begin{equation}\label{vaxr}
 \sum_{\bold a \in BD}{\vert \ga_{\bold a}(\al) \vert} = \lfloor \al \ell^d
 \rfloor \vert BD \vert \geq  \lfloor \al \ell^d
 \rfloor \Big( \vert P_{F,C_1} \vert -  \big( 1 + \frac{\eps}{10^{10d}y}\big) \ell^{-d}
 \big\vert \big( \bigcup_{\Gamma \in \kappa} \Gamma  \big) \cap \Ba_n \big\vert \Big),
\end{equation}
the inequality in \req{vaxr} valid by \req{tbe}. Note that
\begin{eqnarray}
 \sum_{\Gamma \in UH}{\vert \ga_{\Gamma}^* \vert} & = & \sum_{\Gamma \in
 UH}{\lfloor \al \vert \Gamma \vert \rfloor} \geq \big( \al - 
 m_1^{-d} \big)  \sum_{\Gamma \in UH}{\vert \Gamma \vert} \nonumber
 \\
 & \geq & \big( \al - 
 m_1^{-d} \big)  \Big(  \big\vert \big( \bigcup_{\Gamma \in \kappa}
 \Gamma \big) \cap \Ba_n \big\vert \, - \, \frac{\eps n^d}{2^{8d-1} y} \Big),
 \label{vaxt} 
\end{eqnarray}
where \req{reone} was used in the latter inequality. From \req{vaxo},
\req{vaxr}, \req{vaxt} and \req{tkc}, it
follows that
\begin{eqnarray}
\vert \Phi \vert & \geq & ( \al \ell^d - 1 ) \Big( \big( 1 - \eps
2^{-8d} y^{-1} \big) F^d -  \big( 1 + \frac{\eps}{10^{10d}y} \big) \ell^{-d}
 \big\vert \big( \bigcup_{\Gamma \in \kappa} \Gamma  \big) \cap \Ba_n
 \big\vert  \Big) \nonumber \\
 & & \qquad + \, \, \big( \al - m_1^{-d} \big) \Big( \big\vert \big( \bigcup_{\Gamma \in \kappa} \Gamma \big) \cap \Ba_n
 \big\vert  \, - \, \frac{\eps n^d}{2^{8d-1}y} \Big) \nonumber \\
 & \geq & \al  \big( 1 - \eps
2^{-8d} y^{-1} \big) (\ell F)^d  \, - \, \frac{\al \eps
n^d}{2^{8d-1}y} \, - \, \Big( \frac{\eps \al}{10^{10d}y} +  m_1^{-d} \Big) \big\vert
\big( \bigcup_{\Gamma \in \kappa} \Gamma  \big) \cap \Ba_n
 \big\vert  \, - \, F^d \nonumber \\
 & \geq & \al n^d - \big( \eps + m_1^{-d} + \ell^{-d} \big) O(n^d). \label{vaxy} 
\end{eqnarray} 
Note that, by \req{ghnl} and the current definition of $\hat{\Psi}$,
\begin{equation}\label{vaxv}
 \vert \Phi \vert \leq \sum_{\bold a \in BD}{\vert \ga_{\bold a}(\al) \vert}
 \, + \, \sum_{\Gamma \in UH}{\vert \ga_{\Gamma}^* \vert} \, + \,
 \sum_{\bold a_1,\bold a_2 \in BD : \vert \bold a_2 - \bold a_1 \vert
 = 1}{\vert \phi_{\bold a_1,\bold a_2} \vert} \, + \, \sum_{\Gamma \in
 UH}{\vert \hat{\phi}_{\Gamma}}.
\end{equation}
We have that
\begin{equation}\label{vaxw}
\sum_{\bold a \in BD}{\vert \ga_{\bold a}(\al) \vert} + \sum_{\Gamma \in
UH}{\vert \ga_{\Gamma}^* \vert} \leq \al \Big( 
\sum_{\bold a \in BD}{\vert \Ba_{\ell a,\ell} \vert} + \sum_{\Gamma \in
UH}{\vert \Gamma \vert} \Big) \leq \al n^d,
\end{equation}
the second inequality following from the disjointness of  $\cup \big\{ \Ba_{\ell \bold a,\ell}:
\bold a \in BD \big\}$ and $\cup \big\{ \Gamma : \Gamma \in UH
\big\}$, and from $\Ba_{\ell \bold a,\ell}
\subseteq \Ba_n$ for $\bold a \in BD$, and $\Gamma \subseteq \Ba_n$ for
$\Gamma \in UH$. Note that
\begin{equation}\label{vaqb}
 \sum_{\Gamma \in UH}{\vert \hat{\phi}_{\Gamma} \vert} \leq \rho
 \sum_{\Ba_{\bold x,m} \in UH}{m} \leq \frac{\rho n^d}{m_1^{d-1}},
\end{equation}
the latter inequality valid by \req{dft}. Note also that 
\begin{equation}\label{vaqa}
\sum_{\bold a_1,\bold a_2 \in BD : \vert \bold a_1 - \bold a_2 \vert =
1}{\vert \phi_{\bold a_1,\bold a_2} \vert} \leq d \rho \ell F^d \, .
\end{equation}
From \req{vaxv},
\req{vaxw}, \req{vaqa} and \req{vaqb}, it follows that
\begin{equation}\label{vacc}
 \vert \Phi \vert \leq \al n^d + 
 \frac{\rho n^d}{m_1^{d-1}}  + d \rho \ell F^d = \big( \al + \rho m_1^{1-d} + d\rho \ell^{1-d}
 \big) n^d. 
\end{equation}
We now adjust $\Phi$ so that it contains precisely the right number of
sites. It is easier to add sites, so let us ensure that $\Phi$ has too few of
them, by considering $\Phi$ as a subset of $\Ba_{\overline{n}}$, with 
\begin{equation}\label{vace}
\overline{n} = \lfloor
 n \big( 1 + \al^{-1} \rho m_1^{1-d} + d\rho \al^{-1} \ell^{1-d}
 \big)^{1/d} \rfloor + 1.
\end{equation}
 Then \req{vacc} ensures that $\vert \Phi
 \vert \leq \al {\overline{n}}^d$, while \req{vaxy} implies that
\begin{equation}\label{vacd}
 \big\vert \Phi \big\vert \geq \al \overline{n}^d \, -  \, ( m_1^{1-d} +
 \ell^{1-d} + \eps  ) O ({\overline{n}}^d),
\end{equation}
where the constant implicit in $O(\overline{n}^d)$ has no dependence
on $m_1$, $\ell$ and $\eps$.
Let $\al' \in (0,1)$ satisfy 
\begin{equation}\label{altemp}
\al' > \max \{\al,1-\al \}.
\end{equation}
By \req{thetfact} and 
Lemma \ref{lempva},
we may fix $\la \in \mathbb{R}$ 
so that, for each high enough value of $n$, the largest connected
component $\WW_{\overline{n}}$ of $\la$-white sites in $\Ba_{\overline{n}}$
satisfies $\vert \WW_{\overline{n}} \vert \geq \al' \overline{n}^d$.  By
\req{vacd} and \req{altemp}, note that $\Phi \cap \WW_{\overline{n}} \not= \emptyset$, if $\ell$ and $m_1$ are chosen 
to be high enough, and $\eps$ to be small. From the fact that $\vert \Phi
 \vert \leq \al {\overline{n}}^d$, which is ensured by \req{vacc}, we
 see that $\vert \WW_{\overline{n}} \setminus \Phi \vert \geq (\al' -
 \al)\overline{n}^d$, for high values of $n$. Given that 
\begin{equation}\label{quuse}
 \lfloor \al {\overline{n}}^d \rfloor - \vert \Phi \vert \leq (m_1^{1-d} +
 \ell^{1-d} + \eps  )  O ( {\overline{n}}^d ) \leq (\al' -
 \al)\overline{n}^d \, ,
\end{equation}
we may thus find a set  $X \subseteq \WW_{\overline{n}}$ 
for which $\vert X \vert =  \lfloor \al {\overline{n}}^d \rfloor - \vert \Phi \vert$, $\Phi \cap X =
\emptyset$ and $\Phi \cup X$ is connected. (The first inequality in
\req{quuse}
follows from 
\req{vacd}). Note that $\vert \Phi \cup X \vert = \lfloor \al {\overline{n}}^d
\rfloor$ and that
\[
 S \big( \Phi \cup X \big) = S(\Phi) + S(X) \geq S(\Phi) - \la  
 ( m_1^{1-d} + \ell^{1-d} + \eps  )  O({\overline{n}}^d).
\]
By \req{vvc} and \req{vace}, for large values of $n$, the dominant term
in the expression for $S(\Phi)$ in terms of $\overline{n}$ is that in ${\overline{n}}^d$,
with coefficient differing from \req{aftt} by at most
$O(\al^{-1}m_1^{1-d} + d\rho \al \ell^{1-d})$. Provided then that
$m_1$ and $\ell$ are chosen to be high enough, we find that 
\[
 \liminf_{n \to \infty}{ \frac{S(\Phi_n \cup X_n)}{{\overline{n}}^d}} > y.
\]
Given that $\Phi_n \cup X_n \subseteq \Ba_{\overline{n}}$, this implies that
$\liminf k^{-d}\tilde{G}_k(\al) > y$, along that subsequence of $k \in \mathbb{N}$
of the form $k = \overline{n}$ for some $n \in \mathbb{N}$.
Given that the gaps between members of this subsequence are uniformly
bounded, we conclude that 
$\liminf_{k \in \mathbb{N}} k^{-d}\tilde{G}_k(\al) > y$ by considering $\Phi_n \cup
X_n \subseteq \Ba_{\overline{n} + C}$ for each $C$ lying in a finite interval
$[0,C_1]$. $\Box$

\noindent{\bf Proof of Lemma \ref{pthmfive}:}
Firstly, we show that
\begin{equation}\label{cr}
\limsup_{\alpha \downarrow 0}{\frac{\tilde{G}(\alpha)}{\alpha}} \leq N .
\end{equation}
For any $\alpha \in (0,1)$, $\tilde{G}(\alpha)$ is non-random
by Lemma \ref{pthmtwo}. 
Thus, if \req{cr} fails, we may find a pair of small values of
$\eps > 0$ and $\alpha$ such that, 
for each sufficiently high $n$, there exists $\Psi_n \in
\LA_{\Ba_n}$ for which $\vert \Psi_n \vert = \lfloor \al n^d \rfloor$
and $S(\Psi_n) > (N+\eps)\al n^d$. 
Let $\tau_n$
denote a path of length at most $d n$ that contains the origin and
satisfies $\Psi_n \cap \tau_n = \emptyset$, with $\tau_n \subseteq
\Ba_n$ and $\Psi_n \cup \tau_n$ connected. Note that 
\begin{equation}\label{navq}
N_{\vert \Psi_n \cup \tau_n \vert} \geq S \big( \Psi_n \cup \tau_n \big) =
S(\Psi_n) + S(\tau_n) \geq (N+\eps)\al n^d - \la d n,
\end{equation} 
if $\la \in \mathbb{R}$ is chosen so that $\PPP(X_0 \geq - \la) = 1$. Note
also that
\begin{equation}\label{navs}
 \vert \Psi_n \cup \tau_n \vert \leq \vert \Psi_n \vert + \vert \tau_n
 \vert \leq \al n^d + d n.
\end{equation}
From \req{navq} and \req{navs}, it follows that
\[
 \liminf_n \vert \Psi_n \cup \tau_n \vert^{-1}
N_{\vert \Psi_n \cup \tau_n \vert} \geq N + \eps. 
\]
This contradicts the definition of $N$ and establishes \req{cr}.

We complete the proof of the lemma by showing that, for any $\delta_1
>0$, there exists $\delta_2 > 0$ such that $\al \in (0,\delta_2)$
implies that $\tilde{G}(\al) > N(1 - \delta_1)\al$. The following result
follows immediately from \cite[Lemma 3.6]{GLAthree}.
\begin{lemma}\label{lemfourf}
Suppose that $F$ satisfies condition \ref{martins}, and that
$X_{\orig}$ is bounded below almost surely. 
For any $\delta > 0$ and for all sufficiently large $\nu \in \mathbb{N}$, 
we have that
$$
\lim_{m \to \infty}{q(m,\nu,\delta)} = 1 , 
$$ 
where
\begin{eqnarray}
q(m,\nu,\delta) &  = & \mathbb{P} \Big( \exists \ \textrm{animal} \ \xi \subseteq
\{0,\ldots,m \} \times \{ -2 \nu, \ldots, 2 \nu \}^{d-1}: \nonumber \\
   & & \quad  (\orig,2\nu \bold \tau),(m,2\nu \bold \tau) \in \xi \, \, \textrm{for
   each $\bold \tau \in \{-1,1 \}^d$, and} \, \,  S(\xi) \geq (1 - \delta)N
   \vert \xi \vert \Big). \, \Box \nonumber
\end{eqnarray}
\end{lemma}
Let $\theta^b(p) =
\PPP_p(\vert C(\orig) \vert = \infty)$ denote
the density of the infinite cluster for the bond
percolation model of parameter $p$ on $\Zd$ (so that here, $C(\orig) =
\{ \bv \in \Zd: \orig \leftrightarrow \bv \}$ means the collection of
sites accessible to the origin by a path of open edges). By \req{thetfact}, or rather, the comment that follows
it, we may choose $\eps \in (0,1)$ so that $\theta^b (1-\eps) > 1/2$.
By Lemma \ref{lemfourf}, we may find $\nu \in \mathbb{N}$ and $M_1
\in \mathbb{N}$ such that $m \geq M_1$ implies that
$q_{m,\nu,{\delta_1}/2} > 1 - \eps$. For $m \geq M_1$ fixed and to each
$\bold a \in \Zd$, we associate the cube $\phi_{\bold a} =
\{-2\nu,\ldots, 2 \nu \}^d + (m + 2 + 4\nu)\bold a$. To each pair $\bold
a_1,\bold a_2 \in \Zd$ of adjacent sites,
let  $\psi_{\bold a_1,\bold a_2}$ 
denote the collection of sites that are not
elements of $\phi_{\bold a_1}
\cup \phi_{\bold a_2}$, but are a convex
combination of a site in $\phi_{\bold a_1}$ and one in  $\phi_{\bold
a_2}$. Note that  $\psi_{\bold a_1,\bold a_2}$ is
an isometric copy of 
$\{0,\ldots,m\} \times \{ -2\nu,\ldots,2\nu \}^{d-1}$.
We define a
bond percolation $P$ on $\Zd$ of parameter $q_{m,\nu,{\delta_1}/2}$ by
declaring the nearest neighbour edge $(\bold a_1,\bold a_2)$ to be
open if $\psi_{\bold a_1,\bold a_2}$ contains a lattice animal
$\xi_{\bold a_1,\bold a_2}$ satisfying $S(\xi) > (1-{\delta_1}/2) N \vert
\xi \vert$ that contains each of the corners of $\psi_{\bold a_1,\bold
a_2}$. We set $\ell = m + 2 + 4 \nu$, and, for $n \in \mathbb{N}$, write $n
= F \ell + r$, with $r \in \{ 0,\ldots, \ell - 1\}$. We write
$P_{F,C}^b$ for the largest connected component of sites for the bond
percolation $P$ in the box
$\{ C,\ldots, F-1-C \}^d$, 
where $C \in \mathbb{N}$ is given. 
Lemma \ref{lempva} is 
valid for bond percolation with the same proof. Applying this lemma, 
and the inequality $\theta^b
(q_{m,\nu,{\delta}_1}) > 1/2$, for any fixed $C > 0$, there
exists $F_0$ such that, for $F \geq F_0$, $\vert P_{F,C}^b \vert \geq
\frac{1}{2} F^d$. To any connected set $\chi = \{ (\bold a_{1r}, \bold
a_{2r}) : r \in \Gamma \}$ of edges in $\Zd$, we associate the lattice
animal
\begin{equation}\label{rumg}
 \zeta_{\chi} = \bigcup_{r \in \Gamma}{ \Big( \xi_{(\bold a_{1r},\bold
 a_{2r})} \, \cup \,  
C[\bold a_{1r},\bold a_{2r}] \Big) }, 
\end{equation}
where 
\[
C[\bold a_{1r},\bold a_{2r}]  = \bigcup_{i=1}^{2} 
 \big\{ 2 \nu \bold j + \ell \bold a_{ir} : \bold j \in
\{-1,1\}^d, \,  2 \nu \bold j + \ell \bold a_{ir} \, \textrm{is adjacent to a
corner of $\psi_{\bold a_1,\bold a_2}$} \big\}.
\]
The sets $C[\bold a_{1r},\bold a_{2r}]$ in \req{rumg} comprise a
subset of 
corners of the
various cubes $\phi_{\bold a_{ir}}$ in $\zeta_{\chi}$ that ensure that
the set $\zeta_{\chi}$ is connected. For any choice of $\chi$
such that $\bold a_{1r},\bold a_{2r} \in \Ba_{F,1}$ for each $r \in
\Gamma$, we have that $\zeta_{\chi} \subseteq \Ba_n$:
for example, if
$\bold a \in \Ba_{F,1}$ and $\bold j \in \{-1,1\}^d$, then
$\ell_{\infty}(2\nu \bold j + \ell \bold a, \Ba_n^c) \geq
\ell_{\infty}(\phi_{\bold a}, \Ba_n^c) \geq  
\ell_{\infty}(\phi_{\bold a}, \Ba_{F\ell}^c) \geq  
\ell_{\infty}(\ell \bold a, \Ba_{F\ell}^c) - 4 \nu \geq m + 2  >
0$; and $\xi_{(\bold a_{1r},\bold a_{2r})} \subseteq conv(\phi_{\bold
a_{1r}}, \phi_{\bold a_{2r}})$. 
Note that, for
each $\chi \subseteq P^{b}_{F,1}$,
\begin{equation}\label{advar}
 S(\zeta_{\chi}) \geq (1 - {\delta_1}/2) N \sum_{r \in \Gamma} \vert
 \xi_{(\bold a_{1r}, \bold a_{2,r})} \vert \, - \, 2^{d+1} \la \vert
 \chi \vert,
\end{equation}
where $\la$ satisfies $\PPP(X_0 \geq - \la) = 1$.
From the inequalities $\sum_{r \in \Gamma}{\vert \xi_{(\bold a_{1r},\bold a_{2r})}
\vert } \geq (m+1) \vert \chi \vert$ and $\vert \zeta_{\chi} \vert \leq \sum_{r
\in \Gamma}{\vert \xi_{(\bold a_{1r}, \bold a_{2r})} \vert} \,
+ \, 2^{d+1} \vert \chi \vert$, as well as \req{advar}, it follows that
\[
 S(\zeta_{\chi})  \geq 
 \Big( (1 - {\delta_1}/2)N - \frac{2^{d+1}\la}{m+1} \Big) \big( 1 -
 2^{d+1} (m+1)^{-1} \big) \vert \zeta_{\chi} \vert. \nonumber   
\] 
Thus, provided that $M_1$ is chosen to be high enough,
\begin{equation}\label{vatn}
S(\zeta_{\chi}) \geq ( 1 - \delta_1 ) N \vert \zeta_{\chi} \vert,
\end{equation}
since $m \geq M_1$.
Note also that
\begin{eqnarray}
 & & \big\vert \xi_{P_{F,1}^b} \big\vert \geq (m+1) \vert P_{F,1}^b \vert \geq
 \frac{m+1}{2} F^d \label{vatp} \\
 &\geq & \frac{m+1}{2(m+2+4\nu)^d} \big( n - (m + 2+4\nu)
 \big)^d = \frac{(m+1) n^d}{2 (m + 2 + 4\nu)^d}  + O(n^{d-1}). \nonumber 
\end{eqnarray}
For each 
\begin{equation}\label{vatr}
\al \in \Big( 0, \frac{m+1}{2(m+2+4\nu)^d} \Big),
\end{equation}
we may, for high values of $n$,
locate a connected $\chi \subseteq P_{F,1}^b$, for which 
\begin{equation}\label{vatv}
 0 \leq \lfloor \al n^d \rfloor - \vert \zeta_{\chi} \vert \leq (m + 1)(4 \nu
 + 1)^d + 2^{d+1},
\end{equation}
by using \req{vatp} and the fact that the removal of one edge from
$\chi$ changes the size of $\zeta_{\chi}$ by at most the term on the
right-hand-side of \req{vatv}.
We choose $E$ to be any set of $\lfloor \al n^d \rfloor - \vert
\zeta_{\chi} \vert$ sites in $\Ba_n$ for which $\zeta_{\chi} \cap E =
\emptyset$ and $\zeta_{\chi} \cup E$ is connected. We find that 
\begin{equation}\label{vatk}
 \vert \zeta_{\chi} \cup E \vert = \lfloor \al n^d \rfloor 
\end{equation}
and 
\begin{equation}\label{utnyk}
 S ( \zeta_{\chi} \cup E ) \geq S(\zeta_{\chi}) - \la \vert E \vert \geq
 (1 - \delta_1) N \lfloor \al n^d \rfloor - \big( \la + (1 -
 \delta_1) N \big) \big( (m+1)(4 \nu + 1)^d + 2^{d+1} \big),  
\end{equation}
the second inequality valid by \req{vatn} and \req{vatv}. 
From $\zeta_{\chi}\cup E \subseteq \Ba_n$, \req{vatk} and \req{utnyk}, it follows that
\[
 \tilde{G}(\al) \geq \liminf_n n^{-d} S(\zeta_{\chi} \cup E) \geq (1 -
 \delta_1) \al N,
\] 
for each $\al$ satisfying \req{vatr}, as required. $\Box$ \\
We have completed the proof of Proposition \ref{thm-conc}.
We will now obtain from the proposition  
another proof of the inequality $G \leq L N$, valid only when
$X_{\orig}$ is bounded below, but without the need for the hypothesis
$N > 0$. (As mentioned in the Introduction, we learn from this that $L
= 0$ when $N < 0$, in the case where $X_{\orig}$ is bounded below.)

Firstly, we prove that the functions $\tilde{G}_n : [0,1] \to [0,\infty)$ satisfy the following
uniform one-sided Lipschitz condition:
\begin{lemma}\label{lemul}
Let $F$ be a distribution satisfying \req{martins}, whose
support lies in the interval $(\mu,\infty)$. Then, for any $n \in
\mathbb{N}$ and any $\alpha_1,\alpha_2 \in [0,1]$ satisfying $\alpha_1
< \alpha_2$, we have that, almost surely,
\begin{equation}\label{counone}
 n^{-d} \tilde{G}_n(\alpha_1) - n^{-d} \tilde{G}_n(\alpha_2) \leq  - \mu
 \big( \alpha_2 - \alpha_1 \big)   + \vert \mu \vert n^{-d}.
\end{equation}
\end{lemma}
\noindent
{\bf Proof:} Let $\gamma \in B_n$ be a lattice animal satisfying
$S(\gamma) = \tilde{G}_n(\alpha_1)$ and $\vert \gamma \vert = \lfloor n^d
\alpha_1 \rfloor$. We may locate
a set $T \subseteq B_n$ satisfying $\vert T \vert = \lfloor n^d
\alpha_2 \rfloor - \lfloor n^d \alpha_1 \rfloor$, $T \cap \gamma =
\emptyset$, with $T \cup \gamma$ connected. Then,
\begin{eqnarray}
  \tilde{G}_n( \alpha_2 ) & \geq & S \big( T \cup \gamma \big) = S(\gamma)
 + S(T) \nonumber \\
 & \geq & \tilde{G}_n(\alpha_1) + \mu \big( \lfloor n^d
\alpha_2 \rfloor - \lfloor n^d \alpha_1 \rfloor \big) \nonumber \\
 & \geq & \tilde{G}_n(\alpha_1) +  n^d \mu  \big( 
\alpha_2 - \alpha_1 \big) - \vert \mu \vert , \label{countwo} 
\end{eqnarray}
where $S(T)$ was bounded below by using the fact that the weight
assigned to each site is at least $\mu$. $\Box$

\noindent{\bf Proof of Corollary \ref{webcor}:}
By the definition of the quantity $L_n$ and the function $\tilde{G}_n: [0,1]
\to [0,\infty)$, $G_n = \tilde{G}_n \big( n^{-d} L_n \big)$. By Theorem
\ref{pthmthree}, $G = \lim_{n \to \infty}{n^{-d} G_n}$ almost surely. These two
statements imply that 
\begin{equation}\label{labone}
 G = \lim_{n \to \infty}{ n^{-d} \tilde{G}_n \Big( \frac{L_n}{n^d} \Big)}
\end{equation}
almost surely.
In the case where $L<1$, we find that, for any $\epsilon \in (0,1-L)$
and for any $n \in \mathbb{N}$,
\begin{eqnarray}
  n^{-d} G_n \big( n^{-d} L_n \big)  
& = & n^{-d} \big[ \tilde{G}_n \big( n^{-d} L_n \big) - \tilde{G}_n \big( L + \epsilon \big) \big]
     +  n^{-d} \tilde{G}_n \big( L + \epsilon \big) 
 \nonumber \\
 & \leq &  \vert \mu \vert \Big[ L + \epsilon - n^{-d}L_n \Big] +
 \vert \mu \vert \frac{1}{n^d}  +   n^{-d} \tilde{G}_n \big( L + \epsilon \big), \nonumber 
\end{eqnarray}
where Lemma \ref{lemul} was applied in the inequality (note that the
quantity $L$ may be random, so that we are choosing $\eps$ to lie on a
random interval). Thus,
by Lemma \ref{pthmfour} and the definition \req{ldef} of $L$,
\begin{displaymath}
 \limsup_{n}{n^{-d} \tilde{G}_n \big( n^{-d} L_n \big)} \leq  \vert \mu \vert \epsilon
 + \tilde{G} \big( L + \epsilon \big) 
\end{displaymath}
for almost all realizations of $\{X_{\bv}: \bv \in \Zd\}$ for
$L < 1$.
By the continuity of the function $\tilde{G}$ on the interval $(0,1)$
(which is implied by Lemma \ref{pthmtwo}), it follows that
\begin{equation}\label{newone}
 \limsup_{n}{n^{-d} \tilde{G}_n \big( n^{-d} L_n \big)} \leq \tilde{G} \big( L \big),
\end{equation}
when $L < 1$.
From the concavity of $\tilde{G}: (0,1) \to [0,\infty)$ and Lemma 
\ref{pthmfive}, we deduce that, almost surely on the set $\{L<1\}$,
\begin{equation}\label{newtwo}
 \tilde{G} (L) \leq L \limsup_{\alpha \downarrow 0}{\frac{\tilde{G}(\alpha)}{\alpha}}
 = LN.
\end{equation}
It follows from \req{labone}, \req{newone} and \req{newtwo} that $G
\leq LN$ almost surely, in the case where $L < 1$. In the case where
$L=1$, we must demonstrate that $G \leq N$. Let $\gamma \in \Ba_n$ be
a greedy lattice animal in $B_n$. Let $\tau$ denote a path
containing the origin 
that satisfies $\ga \cap \tau = \emptyset$, $\ga \cup \tau \in
\LA_{\Ba_n}$ and $\vert \tau \vert \leq dn$.

We see that
$$
 N_{\vert \gamma \vert + \vert \tau \vert} \geq S(\gamma \cup \tau) =
 S(\gamma) + S(\tau) \geq G_n \, + \, \min \{0,\mu\} nd.
$$
Since $\vert
\gamma \vert + \vert \tau \vert \leq n^d$, we find that
$$
\limsup \frac{N_n}{n} \geq  G.
$$
Given that $N = 
\lim_n n^{-1} N_n$ is almost surely constant by \cite[Theorem
2.1]{GLAthree}, we have obtained $N \geq G$. $\Box$

\noindent{\bf Acknowledgments:} I am grateful to Amir Dembo for
introducing the subject of greedy lattice animals to me, and for
numerous helpful discussions during the course of this project. I
thank him in particular for postulating Theorem \ref{pthmsix} and its
proof by means of Proposition \ref{thm-conc}. Thanks to Peter Teichner
for help in proving Lemma \ref{lemtop}, to G\'abor Pete and Yuval
Peres for useful comments on a draft version of the paper.


\begin{thebibliography}{1}


\bibitem{antpis}
P. Antal and A.Pisztora.
\newblock On the chemical distance for supercritical Bernoulli
percolation. 
\newblock{\em Ann. Probab.}, 24: 1036-1048, 1996. 

\bibitem{GLAone}
J.Theodore Cox, Alberto Gandolfi, Philip~S. Griffin, and Harry Kesten.
\newblock Greedy lattice animals I: upper bounds.
\newblock {\em Annals of Applied Probability}, 3(4):1151--1169, 1993.

\bibitem{GLAthree}
Amir Dembo, Alberto Gandolfi, and Harry Kesten.
\newblock Greedy lattice animals: negative values and unconstrained maxima.
\newblock {\em Ann. Probab.}, 29(1):205--241, 2001.

\bibitem{DP}
Jean-Dominique Deuschel, and Agoston Pisztora.
\newblock Surface order deviations for high-density percolation.
\newblock {\em Prob. Th. rel. Fields} 104:467-482, 1996.



\bibitem{Durrett}
Richard Durrett.
\newblock Probability: theory and examples, second edition, Duxbury,
1996

\bibitem{GLAtwo}
Alberto Gandolfi and Harry Kesten.
\newblock Greedy lattice animals II: linear growth.
\newblock {\em Ann. Appl. Probab.}, 4(1):76--107, 1994

\bibitem{Gr}
Geoffrey Grimmett.
\newblock Percolation, second edition, Springer-Verlag, 1999.

\bibitem{Hughes}
Barry D. Hughes
\newblock Random walks and random environments, volume 2: random
environments, Oxford, 1996.

\bibitem{Kesten}
Harry Kesten.
\newblock Aspects of first passage percolation.
\newblock {\em Lecture Notes in Math.}, 1180:125--264, 1986.


\bibitem{keszha}
Harry Kesten and Yu Zhang.
\newblock Strict inequalities for some critical exponents in
              two-dimensional percolation.
\newblock {\em J. Stat. Phys.}, 46(5-6):1031-1055, 1987.

\bibitem{hirsch}
Morris W. Hirsch
\newblock Differential topology, Springer-Verlag, 1976.

\bibitem{Newman}
C.~Douglas Howard and Charles~M. Newman.
\newblock From greedy lattice animals to {E}uclidean first-passage percolation.
\newblock In {\em Perplexing problems in probability}, pages 107--119.
  Birkh\"auser Boston, Boston, MA, 1999.

\bibitem{Leader}
Imre Leader.
\newblock Discrete isoperimetric inequalities. 
\newblock{\em Proc. Sympos. Appl. Math.}, 44: 57-80, 1991.

\bibitem{Lee}
Sungchul Lee.
\newblock An inequality for greedy lattice animals.
\newblock {\em Ann. Appl. Probab.} 3(4):1170--1188, 1993.

\bibitem{Leetwo}
Sungchul Lee.
\newblock The continuity of ${M}$ and ${N}$ in greedy lattice animals.
\newblock {\em J. Theoret. Probab.}, 10(1):87--100, 1997.

\bibitem{Leeone}
Sungchul Lee.
\newblock The power laws of ${M}$ and ${N}$ in greedy lattice animals.
\newblock {\em Stoch. Proc. Appl.}, 69(2):275--287, 1997.

\bibitem{lss}
Thomas~M. Liggett, Roberto~H. Schonmann, and Alan~M. Stacey.
\newblock Domination by product measures.
\newblock {\em Ann. Probab.}, 25(1):71--95, 1997.

\bibitem{loomwhit}
L.H. Loomis and H. Whitney.
\newblock An inequality related to the isoperimetric inequality.
\newblock{\em Bull. Amer. Math. Soc.}, 55:961--962, 1949. 

\bibitem{Martin}
James~B. Martin.
\newblock Linear growth for greedy lattice animals.
\newblock {\em Stoch. Proc. Appl.}, 98(1):43--66, 2002.

\bibitem{matremy}
Pierre Mathieu and Elisabeth Remy.
\newblock Isoperimetry and heat kernel decay on percolation clusters.
\newblock{\em Ann. Probab.}, 32: 100--128, 2004.


\bibitem{man}
Mikhail~V. Menshikov and Sergie~A. Zuyev.
\newblock Models of $\rho$-percolation.
\newblock {\em Petrozavodsk conference on probabilistic models in discrete
	       mathematics}, 337--347, 1992.


\end{thebibliography}
\end{document}